\documentclass[11pt]{amsart}
\usepackage{epsfig}
\newtheorem{theo}{Theorem}
\newtheorem{prop}{Proposition}
\newtheorem*{conj}{Conjecture}
\newtheorem{lemma}{Lemma}
\newtheorem{defi}{Definition}

\newcommand\R{\mathbb{R}}
\newcommand\C{\mathbb{C}}
\newcommand\Z{\mathbb{Z}}

\newcommand\CP{\mathbb{CP}}
\newcommand\dbar{\bar{\partial}}
\begin{document}
\title[A degree doubling formula]
{A degree doubling formula for braid monodromies and Lefschetz pencils}
\author{Denis Auroux}
\address{Department of Mathematics, M.I.T., Cambridge, MA 02139}
\email{auroux@math.mit.edu}
\author{Ludmil Katzarkov}
\address{UC Irvine, Irvine, CA 92612}
\email{lkatzark@math.uci.edu}
\thanks{The first author was partially supported by a Sloan research
fellowship and NSF Grant DMS-0244844.
The second author was partially supported by a Sloan research
fellowship and NSF Grant DMS-9875383.}
\maketitle
\tableofcontents

\section{Introduction}
It was shown in \cite{A2} that every compact symplectic 4-manifold
$(X,\omega)$ can be realized as an approximately holomorphic branched 
covering of $\CP^2$ whose branch curve is a symplectic curve in $\CP^2$ 
with cusps and nodes as only singularities (however the nodes may have 
reversed orientation). Such a covering is obtained by constructing a 
suitable triple of sections of the line bundle $L^{\otimes k}$, where 
$L$ is a line bundle obtained from the symplectic form (its Chern class
is given by $c_1(L)=\frac{1}{2\pi}[\omega]$ when this class is integral),
and where $k$ is a large enough integer.
Moreover, it was shown in \cite{AK} that the braid monodromy techniques
introduced by Moishezon and Teicher in algebraic geometry (see e.g. 
\cite{MOI81,MOI,Teicher})  can be used in this situation to derive, for each
large enough value of the degree $k$, monodromy invariants which completely
describe the symplectic 4-manifold $(X,\omega)$ up to symplectomorphism.
These invariants are also related to those constructed by Donaldson and
arising from the monodromy of symplectic Lefschetz pencils
\cite{D2}, which also are defined only for large values of $k$.

The monodromy invariants arising from branched coverings or symplectic
Lefschetz pencils give, in principle at least, complete information
about the topology of a symplectic manifold~; for example, it is expected 
that they can be used to symplectically tell apart certain pairs of mutually
homeomorphic algebraic surfaces of general type, such as the Horikawa 
manifolds, which no other currently available symplectic invariant can 
distinguish. However, their practical usefulness is immensely limited by 
the difficulties involved in their calculation, even though the computations
by Moishezon, Teicher and Robb of the braid monodromies for certain simple 
types of algebraic surfaces ($\CP^2$, $\CP^1\times \CP^1$, complete 
intersections) \cite{Teicher} give some reason to be hopeful (see also
\cite{ADKY} for other examples); moreover, the difficulty of comparing
two braid group or mapping class group factorizations up to Hurwitz 
equivalence is a major obstacle.

One of the main technical problems arising
in this program is that the monodromy only becomes a symplectic invariant 
when the degree is large enough, which makes it necessary to
handle whole sequences of braid factorizations. Even when the entire
sequence can be obtained directly out of Moishezon-Teicher style
calculations \cite{Teicher,ADKY}, it is very difficult to understand how 
to extract meaningful information out of the monodromy data, due to the 
lack of a clear relationship between the monodromies arising for different
values of the twisting parameter $k$.

The aim of this paper is to describe an explicit formula relating the
braid monodromy invariants obtained for a given degree $k$ to those
obtained for the degree $2k$. The interest of such a formula is obvious
from the above considerations, especially as direct computations of braid
monodromy often become intractable for large degrees. We also give a similar
formula for the monodromy of symplectic Lefschetz pencils; this formula,
which may have even more applications than that for braid monodromies,
answers a question first considered by Donaldson and for which a partial 
(non-explicit) result has been obtained by Smith \cite{smith}. Although
the formula for pencils is much simpler than that for branched coverings,
the currently available technology for monodromy calculations seems
insufficient to allow a direct proof.
\medskip

The techniques introduced in this paper suggest a
wide range of applications.
First of all, calculations similar to those in this
paper appear in any situation involving iterated branched coverings~;
for example, the invariants defined by Moishezon and Teicher should become
effectively computable for a much larger class of algebraic surfaces (see
e.g.\ \S 7 of \cite{ADKY}). 

An obvious class of applications is to study the properties of
high-degree branched coverings or Lefschetz pencils. For example, Smith
has shown using a degree doubling argument that any compact symplectic 
4-manifold admits a symplectic Lefschetz pencil without reducible singular
fibers \cite{smith}. Although this specific result can be obtained just
from the universality property of degree doubling rather than from the
actual formula, other applications require a more detailed understanding
of the degree doubling process.

More importantly, the degree doubling formula provides
precise information on the behavior of various monodromy-related
invariants as the parameter $k$ increases. For example, it is in principle
possible to describe how the fundamental group of the complement of the 
branch curve, or more generally any other invariant directly related to the 
monodromy group of the branched covering or Lefschetz pencil, depends
on the parameter $k$. It is likely that the conjectures formulated in
\cite{ADKY} can be approached from this perspective. 

In a similar direction, it is reasonable to expect the degree doubling formula 
to yield a stability result for the ``directed Fukaya categories'' introduced 
by Seidel as invariants of Lefschetz pencils~\cite{Seidel}. Unlike the direct
calculation methods following Moishezon and Teicher, our formula makes it 
immediately apparent how Lagrangian spheres lying in standard position 
inside the degree $k$ pencil automatically lie inside the degree $2k$ 
pencil; the explicit description of the additional vanishing cycles makes 
it possible to hope that, under certain assumptions, the degree $2k$ pencil
can be shown to contain no other such spheres. 

Yet another question to which
our result may give an answer is that of whether every branched covering
over $\CP^2$ (or every symplectic Lefschetz pencil) is ``of Donaldson type''
(see the remark at the end of \S 1.2). 

Finally, extensions to 
higher-dimensional settings of the stabilization procedure described here 
are theoretically possible, even though it remains uncertain whether it is
actually possible to carry out the calculations.
\medskip

The remainder of this section is devoted to an overview of braid monodromy
invariants (\S 1.1), followed by a sketch of our approach to
the degree doubling process
and a statement of the main results (\S 1.2 and \S 1.3). 

\subsection{Braid monodromy invariants}

We start by recalling the notations and results (see \cite{AK} or
\cite{Agokova} for details).
Let $f:X\to\CP^2$ be an approximately holomorphic branched covering map as
in \cite{A2} and \cite{AK}~: its topology is mostly described by
that of the branch curve $D\subset\CP^2$, which is symplectic and 
approximately holomorphic. The only singularities of $D$ are double 
points (with either orientation) and cusps (with the complex orientation 
only)~; the branching is of order $2$ at every smooth point of $D$.
Fix a generic projection $\pi:\CP^2-\{pt\}\to\CP^1$ whose pole 
does not belong to $D$. We can assume that $D$ is transverse to the fibers 
of $\pi$ everywhere except at a finite set of non-degenerate tangency points,
where a local model is $x^2=y$ with projection to the $x$ component~; 
moreover, we can also assume
that all the special points of $D$ (tangencies and singular points) lie in 
distinct fibers of $\pi$, and that none of them lies in the fiber above
the point at infinity in $\CP^1$.

The idea introduced by Moishezon in the case of a complex curve is that, 
restricting oneself to the preimage of the
affine subset $\C\subset\CP^1$, the monodromy of $\pi_{|D}$ around its 
critical levels can be used to define a map from $\pi_1(\C-\mathrm{crit})$ 
with values in the braid group $B_d$ on $d=\mathrm{deg}\, D$ strings, called
braid monodromy (see e.g.\ \cite{MOI81})~; this monodromy is encoded by a
factorization of the central element $\Delta_d^2$ of the braid group $B_d$.
Namely, the monodromy around the point at infinity in $\CP^1$, which is given 
by the central braid $\Delta_d^2$, decomposes as the product of the monodromies
around the critical levels of the projection to $\CP^1$, each of these being
conjugate to a power of a half-twist. This construction naturally depends on
the choice of an ordered set of generating loops for the free group
$\pi_1(\C-\mathrm{crit})$.

The same techniques extend almost
immediately to the symplectic setting, and the resulting braid
factorizations are of the form

$$\Delta_{d}^{2} = \prod_{j} (Q_{j}^{-1}\,X_1^{r_j}\,Q_{j}),$$
where $X_1$ is the first standard generator of $B_{d}$ (a positive
half-twist), $Q_j$ are arbitrary braids and $r_{j}\in\{-2, 1, 2, 3\}$. 

The case $r_j=1$ corresponds to a {\it tangency point}, where
the curve $D$ is smooth and tangent to the fiber of the projection $\pi$~;
the case $r_j=2$ corresponds to a {\it nodal point} of $D$~; the case
$r_j=-2$ is the mirror image of the previous one, and corresponds to a
{\it negative self-intersection} of $D$ (this is the only type of point 
which does not occur in the algebraic case)~; and finally the case $r_j=3$
corresponds to a {\it cusp} singularity of $D$.

The above-described braid factorization completely determines the topology
of the curve $D$. However it is well-defined only up to the following two 
algebraic operations.
A {\it Hurwitz move} is the replacement of two consecutive factors
$A$ and $B$ by $ABA^{-1}$ and $A$ respectively (we will say that the factor
$A$ has been ``moved to the right''; the opposite move, which amounts to
replacing $A$ and $B$ by $B$ and $B^{-1}AB$ respectively, will be referred
to as ``moving $B$ to the left''). Another possibility is {\it global
conjugation}, i.e. conjugating all factors simultaneously by a given braid.
A Hurwitz move amounts to an elementary change in the choice of generating
loops for the free group $\pi_1(\C-\mathrm{crit})$, while a global
conjugation amounts to a change of trivialization of the reference fiber
of $\pi_{|D}$.
Two factorizations represent the same curve $D$ if and only if they are
Hurwitz and conjugation equivalent.
\medskip

To recover a map $X\to\CP^2$ from the monodromy invariants we also need a 
{\it geometric monodromy representation}. Let $D\subset\CP^2$ be a curve 
of degree $d$ with cusps and nodes (possibly negative), and let $\C\subset
\CP^2$ be a fiber of the projection $\pi:\CP^2-\{pt\}\to\CP^1$ which
intersects $D$ in $d$ distinct points $q_1,\dots,q_d$. Then, the
inclusion of $\C-\{q_1,\dots,q_d\}$ into $\CP^2-D$ induces a
surjective homomorphism on the fundamental groups. Small loops 
$\gamma_1,\dots,\gamma_d$ around 
$q_1,\dots,q_d$ in $\C$ generate $\pi_1(\CP^2-D)$, with relations
coming from the cusps, nodes and tangency points of $D$. These $d$
loops will be called {\it geometric generators} of $\pi_1(\CP^2-D)$.

Recall that there exists a natural right action of $B_d$ on the free group 
$F_d=\pi_1(\C-\{q_1,\dots,q_d\})$~; denote this action by $*$, and
recall the following definition \cite{MOI}~:

\begin{defi} A geometric monodromy representation associated to a
curve $D\subset\CP^2$ is a surjective group homomorphism
$\theta$ from the free group $F_d=\pi_1(\C-\{q_1,\dots,q_d\})$ to the 
symmetric group $S_n$ of order $n$, such that the $\theta(\gamma_i)$ 
are transpositions $($thus also the $\theta(\gamma_i*Q_j))$ and\medskip

$\theta(\gamma_1\dots\gamma_d)=1,$

$\theta(\gamma_1*Q_j)=\theta(\gamma_2*Q_j)$ if $r_j=1$,

$\theta(\gamma_1*Q_j)$ and $\theta(\gamma_2*Q_j)$ are distinct and
commute if $r_j=\pm 2$,

$\theta(\gamma_1*Q_j)$ and $\theta(\gamma_2*Q_j)$ do not commute 
if $r_j=3$.
\end{defi}

In this definition, $n$ corresponds to the number of sheets of the
covering $X\to\CP^2$~; the various conditions imposed on
$\theta(\gamma_i * Q_j)$ express the natural requirements
that the map $\theta:F_d\to S_n$ should factor through the group
$\pi_1(\CP^2-D)$ and that the branching phenomena should occur in
disjoint sheets of the covering for a node and in adjacent sheets 
for a cusp. The surjectivity of $\theta$ corresponds to the
connectedness of the 4-manifold $X$; more precisely, the image
of $\theta$ is a subgroup of $S_n$ generated by transpositions and
acting transitively on $\{1,\dots,n\}$, which implies surjectivity.

Operations such as Hurwitz moves and global conjugations
should be considered simultaneously on the level of braid factorizations
and on that of the corresponding geometric monodromy representations~:
a Hurwitz move does not affect the geometric monodromy representation,
but when performing a global conjugation by a braid $Q$ it is necessary
to compose $\theta$ with the automorphism of $F_d$ induced by $Q$.
\medskip

In the symplectic case the curve $D$ can have negative nodes, and as a
consequence the uniqueness result obtained in \cite{A2} only holds up to
cancellation of pairs of nodes. An additional possibility is therefore
a {\it pair cancellation} move in the braid factorization, where two 
consecutive factors which are the exact inverse of each other are removed 
from the factorization. The converse move (a {\it pair creation}) is 
also allowed, but only when it is compatible with the geometric monodromy
representation~: adding $(Q^{-1}\,X_{1}^{-2}\,Q).(Q^{-1}\,X_{1}^{2}\,Q)$
somewhere in the braid factorization is only legal if $\theta(\gamma_1 *Q)$
and $\theta(\gamma_2 *Q)$ are commuting disjoint transpositions.

\begin{defi} Two braid factorizations (and the
corresponding geometric monodromy representations) are {\em m-equivalent} if
there exists a sequence of operations which turn one into the other, each 
operation being either a global conjugation, a Hurwitz move, or a pair 
cancellation or creation. 
\end{defi} 

We now summarize the main results of \cite{AK}~:

\begin{theo}[\cite{AK}] The compact symplectic 4-manifold $X$ is
uniquely characterized by the sequence of braid factorizations and 
geometric monodromy representations corresponding to the approximately
holomorphic coverings of $\CP^2$ canonically obtained from sections of
$L^{\otimes k}$ for $k\gg 0$, up to m-equivalence.
\end{theo}

It was also shown in \cite{AK} that conversely, given a (cuspidal negative)
braid factorization and a geometric monodromy representation, one can
recover in a canonical way a symplectic 4-manifold (up to symplectomorphism).

\subsection{The degree doubling process}

We now turn to the topic at hand, namely the phenomena that occur when
the degree $k$ is changed to $2k$. 

In all the following, we will assume
that $k$ is large enough for the uniqueness properties of Theorem 1 to
hold (if the considered coverings happen to be algebraic this assumption
is unnecessary). This makes it possible to choose the most convenient
process for constructing the branch curve for degree $2k$ while ensuring that
the resulting branch curve is indeed equivalent to the canonical one.
As observed in \cite{AK}, one especially interesting way to obtain the
covering map $f_{2k}:X\to\CP^2$ is to start with the covering map
$f_k:X\to\CP^2$ and compose it with the Veronese covering
$V_2:\CP^2\to\CP^2$ given by three generic homogeneous polynomials of
degree 2 (this is a 4:1 covering whose branch curve has degree 6, see
below). The map $V_2\circ f_k$ is clearly an approximately holomorphic 
covering given by sections of $L^{\otimes 2k}$, and its branch curve is
the union of the image by $V_2$ of the branch curve $D_k$ of $f_k$ and
$n=\deg f_k$ copies of the branch curve $C_2$ of $V_2$ (the
branch curve $C_2$ is present with multiplicity $n$ because branching occurs
at every preimage by $f_k$ of a branch point of $V_2$). However at every
point where $V_2(D_k)$ intersects $C_2$ the map $V_2\circ f_k$ presents a
non-generic singular behavior~: e.g., composing the branched coverings
$(x,y)\mapsto (x^2,y)$ and $(x,y)\mapsto (x,y^2)$ yields the singular map
$(x,y)\mapsto (x^2,y^2)$, which needs to be perturbed in order to obtain
a generic behavior. Further small perturbations are required in order to 
separate the multiple copies of $C_2$~; nevertheless, $f_{2k}$ is obtained
as a small perturbation of $V_2\circ f_k$ and its branch curve $D_{2k}$ is
obtained as a small perturbation of $V_2(D_k)\cup n\,C_2$.

For all large enough values of $k$, the approximate holomorphicity and 
transversality properties of the above-described perturbation of 
$V_2\circ f_k$ make it subject to the uniqueness results in \cite{A2}
and \cite{AK}~: the coverings constructed directly and those obtained 
by composition with $V_2$ and perturbation therefore become
isotopic. So, for all large values of $k$ we can indeed hope to compute
the braid factorization of $f_{2k}$ by this method. \medskip

Also observe that a
generic isotopy (1-parameter deformation family) of the curve $D_k$ 
behaves ``nicely'' with respect to the chosen Veronese covering $V_2$, and
therefore yields a generic isotopy of the curve $V_2(D_k)$. Since
generic isotopies do not modify braid factorizations (up to Hurwitz and
conjugation equivalences in the algebraic category, or up to $m$-equivalence
in the symplectic category), we are allowed to perform a generic isotopy
on the curve $D_k$ to place it in the most convenient position with respect
to the ramification curve of $V_2$, and this will not affect the end result.

An important consequence of this observation is that the $k\to 2k$ formula
we are looking for is {\it universal} in the sense that it does not depend
on the branch curve $D_k$ itself but only on its degree $d$ and on the
degree $n$ of the covering $f_k$. Indeed, an isotopy can be used to make
sure that all the special points of $D_k$ (cusps, nodes and tangencies)
lie in a small ball $B\subset\CP^2$ located far away from $V_2^{-1}(C_2)$,
and that $D_k$ looks like a union of $d$ lines outside of the ball $B$.
For example, we can take $V_2$ to be a small perturbation of the non-generic
quadratic map $V_2^0:(x:y:z)\mapsto (x^2:y^2:z^2)$, for which the 
ramification curve consists of three lines (the coordinate axes), and we 
can use a linear transformation to contract all the special points of $D$ 
(tangencies, nodes, cusps) into an arbitrarily small ball $B$ centered at
the point $(1:1:1)$.

With this setup, the contribution of $D_{2k}\cap V_2(B)$
to the braid monodromy is the same as that of $D_k\cap B$, and the braid 
monodromy coming from $D_{2k}\cap (\CP^2-V_2(B))$ does not depend on the 
curve $D_k$ but only on its degree and on the geometric monodromy 
representation $\theta$. The braid factorization corresponding to $f_{2k}$
is therefore of the form $$F_k\cdot \mathbf{U}_{d,n,\theta},$$
where $F_k$ is the braid factorization for $f_k$ (after a suitable
embedding of $B_d$ into the larger braid group $B_{\bar{d}}$ corresponding to
$D_{2k}$) and $\mathbf{U}_{d,n,\theta}$ is a word in $B_{\bar{d}}$ depending only on
$d$, $n$ and $\theta$ ($\bar{d}=2d+6n=\deg D_{2k}$).

From the above considerations, the strategy for obtaining the formula 
giving the braid factorization for $D_{2k}$ in terms of the braid 
factorization for $D_k$ is the following. First one needs to understand
the braid factorizations corresponding to the two curves $V_2(D_k)$ and
$C_2$ taken separately. More specifically, the braid factorization for
$V_2(D_k)$ is obtained from that for $D_k$ via a ``folding'' formula
describing the effect of the quadratic map $V_2$; the braid
factorization for $C_2$ (and consequently for $n$ copies of $C_2$) is
obtained by degenerating it to the branch curve of the map $V_2^0$, which
consists of three lines forming a triangle, giving rise to three
similar-looking
contributions from their mutual intersections. Next, one has to study the 
phenomena that arise near the intersections of $C_2$ with $V_2(D_k)$~; these
again give rise to three similar contributions (one for each line in the
branch curve of $V_2^0$). Finally 
more calculations are required in order to combine these ingredients into
a formula for $D_{2k}$. The main result is the following (see \S 3 for
notations):

\begin{theo}
Let $f_k:(X^4,\omega)\to\CP^2$ be
an approximately holomorphic branched covering given by three sections
of $L^{\otimes k}$. Denote by $D_k$ the branch curve of $f_k$, and let
$d=\deg D_k$ and $n=\deg f_k$. Denote by $F_k$ the braid factorization 
corresponding to $D_k$, and assume that $d\le n(n-1)$. Then, with the 
notations of \S 3, the braid factorization corresponding to the branch 
curve $D_{2k}$ of $f_{2k}$ is given up to $m$-equivalence by a formula 
of the following type, provided that $k$ is large enough:
\begin{equation}\label{eq:k2kbasic}
\Delta_{2d+6n}^2=\mathbf{T}_d\cdot \iota(F_k)\cdot \mathbf{I}^\alpha_{d,n,\theta}
\cdot \mathbf{I}^\beta_{d,n,\theta}\cdot\mathbf{V}^{\alpha\beta}_n
\cdot \mathbf{V}^{\alpha\gamma}_n\cdot \mathbf{V}^{\beta\gamma}_n
\cdot \mathbf{I}^\gamma_{d,n,\theta}.
\end{equation}
\end{theo}

In this formula, $F_k$ is viewed as a factorization in $B_{2d+6n}$
using a certain natural embedding $\iota:B_d\hookrightarrow B_{2d+6n}$.
The other terms correspond to universal contributions (depending only
on $d$, $n$ and $\theta$): more precisely, $\mathbf{T}_d$ arises from
the folding of $D_k$ by the quadratic map $V_2$, while
$\mathbf{V}^{\alpha\beta}_n$, $\mathbf{V}^{\alpha\gamma}_n$,
$\mathbf{V}^{\beta\gamma}_n$ arise from the braid monodromy of $n$ parallel
copies of the curve $C_2$, and $\mathbf{I}^\alpha_{d,n,\theta}$,
$\mathbf{I}^\beta_{d,n,\theta}$, $\mathbf{I}^\gamma_{d,n,\theta}$ correspond
to the intersections of $V_2(D_k)$ with $C_2$. The individual factors
in each of these expressions are described  in \S \ref{ss:braidthm}.

The proof of Theorem 2 is carried out in Sections 2 and 3 of this paper:
the strategy of proof outlined above is carefully justified in \S 2
(cf.\ in particular Propositions 1 and 2);
general properties of the braid group and notations are
introduced in \S \ref{ss:generalities}~; \S \ref{ss:folding} 
describes the folding formula which gives
the braid factorization for $V_2(D_k)$~; the braid factorization
of the branch curve $C_2$ of $V_2$ is computed in \S \ref{ss:v2}~; the local
perturbation procedure to be performed near the intersections of $C_2$
with $V_2(D_k)$ is described in \S \ref{ss:3pt}~; \S \ref{ss:assembling}
deals with the assembling procedure that yields the braid factorization for 
$D_{2k}$ from the previous ingredients~; finally, the calculation is 
completed in \S \ref{ss:braidthm}.\medskip

{\bf Remark.} More generally, this procedure applies to any situation
involving {\it iterated} branched coverings~: given two approximately
holomorphic branched covering maps $f$ and $g$, the composed map 
$h=g\circ f$ has a non-generic behavior at each of the intersection points
of the branch curves of $f$ and $g$~; however, the perturbation procedure
described in \S \ref{ss:3pt} also applies to this situation, and calculations 
similar to those of Section 3 can be used to compute the braid monodromy of a
``generic'' perturbation $\tilde{h}$ of $h$.

Also observe that, in the case of complex surfaces, the manner in which we
perturb iterated coverings, even though it is not holomorphic, is very 
similar and in a sense equivalent to the corresponding construction in
complex geometry. In particular, even though our computations are always
performed up to $m$-equivalence (allowing cancellations of pairs of nodes),
in the case of complex manifolds a formula very similar to 
(\ref{eq:k2kbasic}) holds up to Hurwitz and conjugation equivalence 
(without node cancellations). The only issue requiring particular attention
is the manner in which the multiple copies of the curve
$C_2$ are perturbed away from each other (see \S \ref{ss:braidthm}), where
the most natural choice in the approximately holomorphic context may
be slightly different from a holomorphic perturbation; still, evidence
suggests that in practice Theorem 2 does hold up to Hurwitz and
conjugation equivalence for most complex surfaces (see the end of \S
\ref{ss:braidthm}).
\medskip

{\bf Remark.} The branched coverings constructed in \cite{AK} and the
symplectic Lefschetz pencils constructed by Donaldson enjoy transversality
properties which intuitively ought to make their topology very special among 
all possible coverings or pencils. It is therefore interesting to ask for
criteria indicating whether a given covering map (or Leschetz pencil) is 
``of Donaldson type''; more precisely, the question is to decide whether,
after stabilizing by repeatedly applying the degree doubling formula, the 
monodromy data of the given covering map $X\to\CP^2$ eventually 
coincides with the invariants of $X$ given by Theorem 1. 
This question can be reformulated in two equivalent ways
(similar statements about Lefschetz pencils can also be considered)~:

1. Given two sets of monodromy invariants representing branched coverings of 
$\CP^2$ with the same total space up to symplectomorphism, do they always 
become $m$-equivalent to each other by repeatedly applying the degree 
doubling formula~?

2. Is the set of all compact symplectic 4-manifolds with integral symplectic
class up to scaling of the symplectic form in bijection with the set of all
possible braid factorizations and geometric monodromy representations up to
$m$-equivalence and stabilization by degree doubling ?

\subsection{Degree doubling for symplectic Lefschetz pencils}

A direct application of the degree doubling formula for braid monodromies
is a similar formula for the monodromy of the symplectic Lefschetz pencils
constructed by Donaldson \cite{D3}. Indeed, recall from \cite{D3} that every 
compact symplectic 4-manifold admits a structure of Lefschetz pencil
determined by two sections of $L^{\otimes k}$ for large enough $k$.
The monodromy of such a Lefschetz pencil is described by a word in
the mapping class group of a Riemann surface. As explained in \cite{AK},
Lefschetz pencils and branched coverings are very closely related to each
other, and the monodromy of the Lefschetz pencil can be computed explicitly
from the braid factorization and the geometric monodromy representation
describing the covering. 

More precisely, the geometric monodromy
representation $\theta$ determines a group homomorphism $\theta_*$ from 
a subgroup $B_d^0(\theta)$ of $B_d$ to the mapping class group $M_g$ of 
a Riemann surface of genus $g=1-n+(d/2)$~; the braid monodromy is contained
in $B_d^0(\theta)$, and the monodromy of the Lefschetz pencil is obtained
by composing the braid monodromy with $\theta_*$. It was shown in
\S 5 of \cite{AK} that the nodes and cusps of the branch curve do not
contribute to the monodromy of the Lefschetz pencil (the corresponding
braids lie in the kernel of $\theta_*$), while the half-twists corresponding
to the tangency points of the branch curve yield Dehn twists in $M_g$.

Using this description, we derive in Section 4 a degree
doubling formula for Lefschetz pencils. The relation between braid
groups and mapping class groups of Riemann surfaces with boundary components
is described in more detail in \S \ref{ss:mapg}, and the degree doubling 
formula is obtained in \S \ref{ss:mapgthm}.
\medskip

\noindent {\bf Acknowledgements.} We are very grateful to S. Donaldson and
M. Gromov for their constant attention to this work. The second author would
also like to thank IHES for the extremely pleasant working conditions.

\section{Stably quasiholomorphic coverings}

\subsection{Quasiholomorphic coverings and braided curves}
We now describe in more detail the geometric properties of the covering
maps and branch curves that we will be considering.

\begin{defi}
A real $2$-dimensional singular submanifold $D\subset\CP^2$ is a {\em
braided curve} if it satisfies the following properties~: $(1)$ the only 
singular points of $D$ are cusps (with positive orientation) and transverse
double points (with either orientation)~; $(2)$ the point $(0:0:1)$ does not
belong to $D$~; $(3)$ the fibers of the projection $\pi:(x:y:z)\mapsto (x:y)$
are everywhere transverse to $D$, except at a finite set of nondegenerate
tangency points where a local model for $D$ in orientation-preserving
coordinates is $z_2^2=z_1$~; $(4)$ the cusps, nodes and tangency points are
all distinct and lie in different fibers of $\pi$.
\end{defi}

This notion is a topological analogue of the notion of quasiholomorphic
curve as described in \cite{AK}. In fact, a singular curve in $\CP^2$ 
can be described by a braid factorization with factors of degree $1$, 
$\pm 2$, and $3$ if and only if it is braided. As observed in \cite{AK},
every braided curve is isotopic to a symplectic curve, as follows
immediately from applying the transformation $(x:y:z)\mapsto
(x:y:\epsilon z)$, with $\epsilon$ sufficiently small.
However, the branch curves obtained from
asymptotically holomorphic families of branched coverings satisfy much
more restrictive geometric assumptions.

More precisely, recall that the notion of
quasiholomorphicity only makes sense for a sequence of branch curves 
obtained for increasing values of the degree $k$, and that the resulting 
geometric estimates improve when $k$ increases. The geometric properties
that follow immediately from the definitions and arguments in \cite{A2}
and \cite{AK} are the following. Recall that $(X,\omega)$ is endowed with
a compatible almost-complex structure $J$ and the corresponding metric $g$,
and that we rescale this metric to work with the metric $g_k=k\,g$.

\begin{defi}
A sequence of sections $s_k$ of complex vector bundles $E_k$ over $X$
(endowed with Hermitian metrics and connections) is {\em asymptotically 
holomorphic} if there exist constants $C_j$ independent of $k$ such that
$|\nabla^j s_k|_{g_k}\le C_j$ and $|\nabla^{j-1}\dbar s_k|_{g_k}\le C_j
k^{-1/2}$ for all $j$. 

The sections $s_k$ are {\em uniformly transverse to $0$} if there
exists a constant $\gamma>0$ such that, at every point $x\in X$ where
$|s_k(x)|\le\gamma$, the covariant derivative $\nabla s_k(x)$ is surjective 
and has a right inverse of norm less than $\gamma^{-1}$ w.r.t. $g_k$
(we then say that $s_k$ is $\gamma$-transverse to $0$).

If the sections $s_k$ are asymptotically holomorphic and uniformly transverse
to $0$ then for large $k$ their zero sets are smooth asymptotically
holomorphic symplectic submanifolds.
\end{defi}

\begin{defi}
A sequence of branched covering maps $f_k:X\to\CP^2$ determined by
asymptotically holomorphic sections $s_k=(s_k^0,s_k^1,s_k^2)$ of
$\C^3\otimes L^{\otimes k}$ for $k\gg 0$ is {\em quasiholomorphic} if there 
exist constants $C_j$, $\gamma$, $\delta$ independent of $k$,
almost-complex structures $\tilde{J}_k$ on $X$, and finite subsets
$\mathcal{F}_k\subset X$,
such that the following properties hold (using $\tilde{J}_k$ to define
the $\dbar$ operator)~: 

$(0)$ $|\nabla^j(\tilde{J}_k-J)|_{g_k}\le C_j k^{-1/2}$
for every $j\ge 0$~; $\tilde{J}_k=J$ outside of the 
$2\delta$-neighborhood of $\mathcal{F}_k$~; 
$\tilde{J}_k$ is integrable over the $\delta$-neighborhood of
$\mathcal{F}_k$~;

$(1)$ the norm of $s_k$ is everywhere bounded from below by $\gamma$~;
as a consequence, $|\nabla^j f_k|_{g_k}\le C_j$ and 
$|\nabla^{j-1}\dbar f_k|_{g_k}\le C_j k^{-1/2}$ for all $j$~; 

$(2)$ $|\nabla f_k(x)|_{g_k}\ge \gamma$ at every point $x\in X$~; 

$(3)$ the $(2,0)$-Jacobian $\mathrm{Jac}(f_k)=\det\partial f_k$ is
$\gamma$-transverse to $0$~; in particular it
vanishes transversely along a smooth symplectic curve $R_k\subset X$
(the ramification curve).

$(3')$ the restriction of $\dbar f_k$ to
$\mathrm{Ker}\,\partial f_k$  vanishes at every point of $R_k$~;

$(4)$ the quantity $\partial(f_{k|R_k})$, which can be seen as a section
of a line bundle over $R_k$, is $\gamma$-transverse to $0$ and vanishes at 
a finite subset $\mathcal{C}_k\subset\mathcal{F}_k$ 
(the cusp points of $f_k$)~; in particular
$f_k(R_k)=D_k$ is an immersed symplectic curve away from the image of 
$\mathcal{C}_k$~;

$(5)$ $f_k$ is
$\tilde{J}_k$-holomorphic over the $\delta$-neighborhood of $\mathcal{F}_k$~;

$(6)$ the section $(s_k^0,s_k^1)$ of $\C^2\otimes L^{\otimes k}$ is
$\gamma$-transverse to $0$~; as a consequence $D_k$ remains away from the
point $(0\!:\!0\!:\!1)$~;

$(7)$ letting $\phi_k=\pi\circ f_k:R_k\to\CP^1$, the quantity
$\partial(\phi_{k|R_k})$ is $\gamma$-transverse to $0$ over $R_k$, 
and it vanishes
over the union of $\mathcal{C}_k$ with a finite set $\mathcal{T}_k$ 
(the tangency points of $D_k$)~; moreover, $\dbar f_k=0$ at every point of
$\mathcal{T}_k$~;

$(8)$ the projection $f_k:R_k\to D_k$ is injective outside the singular
points of $D_k$, and the branch curve $D_k$ is braided.
\end{defi}

The main result of \cite{AK} is the existence, for large enough values of
$k$, of quasiholomorphic covering maps $X\to\CP^2$ determined by sections
of $\C^3\otimes L^{\otimes k}$, canonical up to isotopy. The braid monodromy
invariants corresponding to these coverings are those mentioned in 
Theorem 1.

\subsection{Stably quasiholomorphic coverings}

We wish to construct and study branched covering maps which, in addition
to being quasiholomorphic, behave nicely when composed with a quadratic
holomorphic map from $\CP^2$ to itself. For this purpose, we extend in
the following way the notions defined in the previous sections~:

\begin{defi}
We say that the image $D\subset \CP^2$ of a smooth curve $R$ by a map $f$
is {\em locally braided} if 
there exists a finite number of open subsets $U_j\subset R$, whose union
is $R$, such that for all $j$ the image $f(U_j)\subset D$ is a braided 
curve in $\CP^2$.\end{defi}

In other words, a locally braided curve is similar to a braided curve except
that it is merely immersed outside its cusps, without any
self-transversality property~; although the cusps and tangencies of a
locally braided curve are still nondegenerate and well-defined, phenomena
such as self-tangencies might occur. For example, if the definition of
a quasiholomorphic covering is relaxed by removing condition $(8)$, the
branch curve $D_k$ is only locally braided.

Although a locally braided branch curve does not have a well-defined braid
monodromy, an arbitrarily small perturbation ensures self-transversality
and yields a braided curve~; it is easy to check that the braid monodromies
of all possible resulting curves are m-equivalent, as the only
phenomenon which can occur in a generic 1-parameter family is the 
cancellation of pairs of double points.

\begin{defi}
A sequence of branched covering maps $f_k:X\to\CP^2$ determined by
asymptotically holomorphic sections $s_k=(s_k^0,s_k^1,s_k^2)$ of
$\C^3\otimes L^{\otimes k}$ for $k\gg 0$ is {\em stably quasiholomorphic} 
if, with the same notations as in Definition $5$, the following 
properties hold~: 

$(1)$ the covering maps $f_k$ are quasiholomorphic~;

$(2)$ the sections $s_k^0$, $s_k^1$ and $s_k^2$ of $L^{\otimes k}$ are 
$\gamma$-transverse to $0$~;

$(3)$ the sections $(s_k^0,s_k^1)$, $(s_k^0,s_k^2)$ and $(s_k^1,s_k^2)$ of 
$\C^2\otimes L^{\otimes k}$ are $\gamma$-transverse to $0$~;

$(4)$ let $\pi^0$, $\pi^1$ and $\pi^2$ be the projections $(x:y:z)\mapsto
(y:z)$, $(x:y:z)\mapsto (x:z)$ and $(x:y:z)\mapsto (x:y)$ respectively, and
define $\phi_k^i=\pi^i\circ f_k$~; the quantity
$\partial((\phi_k^i)_{|(s_k^i)^{-1}(0)})$ is $\gamma$-transverse to $0$ 
over $(s_k^i)^{-1}(0)$ for $i=0,1,2$~;

$(5)$ the quantity $|\partial\phi_k^i|_{g_k}$ is bounded from below by 
$\gamma$ over $(s_k^i)^{-1}(0)$~;

$(6)$ $\mathcal{F}_k=\mathcal{C}_k\cup\mathcal{T}_k\cup\mathcal{I}_k$,
where $\mathcal{T}_k$ is the set of tangency points and $\mathcal{I}_k$ is 
the set of points of $R_k$ where one of the three sections $s_k^i$ 
vanishes.
\end{defi}

We have the following extension of the main results of \cite{A2} and
\cite{AK}, which will be proved in \S 2.3~:

\begin{prop}\label{prop:stqh}
For all large values of $k$, there exist asymptotically holomorphic
sections $s_k$ of $\C^3\otimes L^{\otimes k}$ such that the corresponding
projective maps $f_k:X\to\CP^2$ are stably quasiholomorphic coverings. 
Moreover, for large $k$ the topology of these covering maps is canonical up 
to isotopy and cancellations of pairs of nodes in the branch curve.
\end{prop}

More precisely, the uniqueness statement means that, given two sequences
of stably quasiholomorphic coverings, it is possible for large $k$ to find
an interpolating 1-parameter family of covering maps, all of which are
stably quasiholomorphic, except for finitely many parameter values where
a cancellation or creation of a pair of nodes occurs in the branch curve.

The following result will be used in \S 3.2 to compute the braid monodromy
of the folded branch curve $V'_2(D_k)$:

\begin{lemma}
Consider the two maps $V_2^0:(x:y:z)\mapsto (x^2:y^2:z^2)$ and 
\mbox{$\psi_a:(x:y:z)\mapsto (x:ay+(1-a)x:az+(1-a)x)$} from $\CP^2$ to
itself, and let $f_k$ be a sequence of stably
quasiholomorphic covering maps with branch curves $D_k\subset\CP^2$.
Then the curves $V_2^0(D_k)$ are locally braided. Moreover, if we assume
that $(0:1:1)\not\in D_k$ and that none of the nodes of $D_k$ lies on the
line $L_0=\{(0:y:z)\}$, then for all sufficiently small non-zero values of 
$a\in\C$
the curves $V_2^0(\psi_a(D_k))$ are locally braided and isotopic to
$V_2^0(D_k)$ through locally braided curves.

Furthermore, these properties remain true if $V_2^0$ is replaced by any generic
holomorphic quadratic map $V'_2$ from $\CP^2$ to itself which differs from 
$V_2^0$ by less than $\gamma'$, for some constant $\gamma'$ independent of 
$k$.
\end{lemma}

\proof
The ramification curve of $V_2^0$ consists of three
lines $L_0=\{(0:y:z)\}$, $L_1=\{(x:0:z)\}$ and $L_2=\{(x:y:0)\}$. Moreover,
$V_2^0$ maps each fiber of $\pi$ to a fiber of $\pi$.
Therefore, let $C\subset\CP^2$ be a locally braided curve satisfying
the following properties:
$(a)$ $C$ is transverse to the lines $L_0,L_1,L_2$ and
avoids their intersection points;
$(b)$ the cusps and tangency points of $C$ lie away from $L_0,L_1,L_2$;
$(c)$ at any point $p\in C\cap L_i$,
the curve $C$ is transverse to 
the fiber of the projection $\pi^i$ through $p$; $(d)$ the curve $C$ is
holomorphic near its tangency points and near its intersections with
$L_0\cup L_1$. Then we conclude that
$V_2^0(C)$ is locally braided and holomorphic near its tangency points.

Indeed, conditions $(a)-(c)$ imply that the restriction of $V_2^0$
to $C$ is an immersion, because $V_2^0$ is a local diffeomorphism away
from $L_i$, and $C$ is transverse to the kernel of the differential of
$V_2^0$ at its intersection points with $L_i$. Moreover, $(a)$ also
implies that $V_2^0(C)$ avoids the point $(0:0:1)$.
The cusps of $V_2^0(C)$ are exactly the images of those of
$C$, while the tangency points of $V_2^0(C)$ are of two types: on one hand,
the images of the tangency points of $C$, and on the other hand, the images
of the intersection points of $C$ with either $L_0$ or $L_1$.
Property $(d)$ implies that $V_2^0(C)$ is holomorphic near its tangency
points, and because $C$ is locally braided and transverse to $L_0\cup L_1$, 
these tangencies are non-degenerate, which implies that $V_2^0(C)$ is
locally braided.

We now check that, as a consequence of Definition 7, the curves $D_k$
satisfy properties $(a)-(d)$. Indeed, property (3) of Definition 7 implies
that $f_k$ is a local diffeomorphism wherever two of the components of
$s_k$ are very small; therefore $D_k$ avoids the intersection points of
$L_0,L_1,L_2$. Moreover, property (2) of Definition 7 implies that
$(s_k^i)^{-1}(0)$ is smooth and $D_k$ is transverse to $L_i$ for $i=0,1,2$;
so condition $(a)$ is satisfied. This transversality requirement also
implies that the tangency points of $D_k$ do not lie on $L_0$ or $L_1$;
in the case of $L_2$ we appeal to property (5) of Definition 7 to reach the
same conclusion (recall that by definition $\partial\phi_k^2$ vanishes at
tangency points). Furthermore, property (4) of Definition 7 means that
the restriction of $\phi_k^i$ to $(s_k^i)^{-1}(0)$ has non-degenerate
critical points, which implies that the intersection multiplicity
of $R_k$ with $(s_k^i)^{-1}(0)$ at such a point is always 1 and prevents 
a cusp of $D_k$ from lying on $L_i$. Therefore $(b)$ holds.

Condition $(c)$ is a direct consequence of property (5) of Definition 7,
observing that the points where $D_k$ is tangent to the fiber of $\pi^i$ 
are precisely the critical points of $\phi_k^i$. Finally, condition $(d)$
follows immediately from property (6) of Definition 7. Therefore $D_k$
satisfies $(a)-(d)$, which implies that $V_2^0(D_k)$ is locally braided
and holomorphic near its tangency points.

We now consider the curve $V_2^0(\psi_a(D_k))$. Observe that, when $a\to 0$,
the linear map $\psi_a$ fixes the points of $L_0$ and collapses all other
points towards $p_0=(1:1:1)$. Moreover, $\psi_a$ maps each fiber of $\pi$ to
a fiber of $\pi$. If we assume that the nodes of $D_k$ don't lie on $L_0$,
then for sufficiently small values of $a$ the curve $\psi_a(D_k)$ becomes
arbitrarily close to a union of $d=\deg D_k$ lines, each joining a
point of $D_k\cap L_0$ to $p_0$. The requirement $(0:1:1)\not\in D_k$
ensures that none of these lines is a fiber of the projection $\pi^0$.
The cusps and tangency points of
$\psi_a(D_k)$ are the images of those of $D_k$ and hence all lie in a small
ball centered at $p_0$; moreover the holomorphicity of $D_k$ near the points
of $D_k\cap L_0$ implies that $\psi_a(D_k)$ is holomorphic outside of a
small ball centered at $p_0$.  
Therefore $\psi_a(D_k)$ satisfies the
conditions $(a)-(d)$ listed above, and $V_2^0(\psi_a(D_k))$ is locally
braided and holomorphic near its tangency points for all sufficiently 
small values of $a$.

Observing that properties $(a)-(c)$ are open conditions, one easily checks
that, if the behavior of the curve $D_k$ is generic (which can be ensured
by a small perturbation), then the curves $\psi_a(D_k)$ (or
small perturbations thereof) satisfy
$(a)-(c)$ for all but a discrete set of values of $a$. Therefore, observing
that $\psi_1=\mathrm{Id}$ and choosing a suitable path $a(t)$,
there exists an isotopy between $D_k$ and $\psi_a(D_k)$ through braided
curves satisfying conditions $(a)-(c)$. Although the
possible lack of holomorphicity of $\psi_{a(t)}(D_k)$ near its intersections
with $L_1$ may prevent $(d)$ from holding, this specific
requirement is actually not needed to ensure that $V_2^0(\psi_{a(t)}(D_k))$ 
is locally braided. Therefore, $V_2^0(\psi_a(D_k))$ is isotopic 
to $V_2^0(D_k)$ through locally braided curves.

Finally, we consider a holomorphic quadratic map $V'_2$ sufficiently
close to $V_2^0$. Our main observation is that the curves $V'_2(D_k)$ and
$V_2^0(D_k)$ are $C^1$-close to each other. Therefore, because $V_2^0(D_k)$
is locally braided and holomorphic near its tangency points (which are
all non-degenerate), the curve $V'_2(D_k)$ is also locally braided; indeed,
if $V'_2$ is sufficiently close to $V_2^0$ then every point where
$V'_2(D_k)$ fails to be transverse to the fibers of $\pi$ necessarily lies
close to a tangency point of $V_2^0(D_k)$. Furthermore, choosing a continuous
deformation of $V_2^0$ into $V'_2$, it is clear that $V_2^0(D_k)$ and
$V'_2(D_k)$ are isotopic to each other among locally braided curves.

The reason why we can obtain a uniform estimate $\gamma'$ on the maximum
admissible value of $\|V'_2-V_2^0\|_{C^1}$ is the existence of uniform
estimates on the geometry of $D_k$. Indeed, by carefully keeping track
of the uniform estimates given by Definitions 5 and 7, it is possible
to derive uniform lower bounds for all geometrically relevant quantities,
such as the distance from $D_k$ to the intersection points of the lines $L_i$,
the transversality angle at the intersections of $D_k$ with $L_i$,
the distance between $L_i$ and the cusps and tangency
points of $D_k$, the second derivative of $\pi_{|D_k}$ at the tangency points 
of $D_k$ and its first derivative away from these points, ... This yields
uniform estimates on the geometry of $V_2^0(D_k)$ near its tangency points
and implies that the property of being locally braided remains valid up to
a certain size of perturbation of $V_2^0$ which can be estimated explicitly
in terms of the various bounds.

Moreover, recalling from above the behavior of $\psi_a$ for small values of
$a$, we can similarly show that if $a$ is sufficiently small then 
$V'_2(\psi_a(D_k))$ is locally braided and isotopic to $V_2^0(\psi_a(D_k))$
through locally braided curves; one simply needs to choose $V'_2$ generic
in order to ensure that the images by $V'_2$ of the lines joining $p_0$ to
the points of $D_k\cap L_0$ are smooth conics.

We conclude in particular that the images by $V_2^0$ and $V'_2$ of $D_k$ and
$\psi_a(D_k)$ are all mutually isotopic among locally braided curves, and
their braid monodromies are m-equivalent to each other.
\endproof

The following observation plays a crucial role in our strategy to prove
Theorem 2: given a generic holomorphic quadratic map $V'_2$ close to $V_2^0$,
the composed maps $V'_2\circ f_k$ already satisfy most of the properties
expected of quasiholomorphic coverings except at the points where the
branch curve of $f_k$ intersects that of $V'_2$.

\begin{prop}
Let $f_k$ be a family of stably quasiholomorphic coverings, and let $V'_2$
be a generic holomorphic quadratic map close to $V_2^0$. Then, given any
fixed constant $d_0>0$, there exist constants $C_j$, $\gamma$, $\delta$
independent of $k$ (but depending on $V'_2$ and on $d_0$) such that the 
composed maps 
$f'_{2k}=V'_2\circ f_k$ satisfy all the properties of Definition 5, except
for properties $(3')$ and $(8)$, at every 
point of $X$ whose $g_k$-distance to $\mathcal{I}'_k=R_k\cap f_k^{-1}(R'_2)$ 
is larger than $d_0$ ($R_k$ and $R'_2$ are the ramification curves of $f_k$ 
and $V'_2$ respectively).
\end{prop}

\proof
The projective map $f'_{2k}=V'_2\circ f_k$ is defined by a section $Q(s_k)$
of $\C^3\otimes L^{\otimes 2k}$, each of its three components being a 
quadratic expression $Q_i(s_k)$ ($0\le i\le 2$) in the three sections 
defining $f_k$. It is therefore easy to show that the sections $Q(s_k)$ are 
asymptotically holomorphic.

Because the projective map $V'_2$ induced by the polynomials 
$Q_i$ is well-defined, the inequality $|Q(s)|\ge c\,|s|^2$ holds for some
constant $c>0$. Therefore, the existence of a uniform lower bound on 
$|s_k|$ at every point of $X$ implies that of a uniform lower bound on 
$|Q(s_k)|$, and so property (1) of Definition 5 is satisfied everywhere.

As observed above, by property (2)
of Definition 7 the branch curve of $f_k$ is uniformly transverse to the
ramification curve
of $V_2^0$ and hence to that of $V'_2$. Therefore, if a point $x\in X$ lies
close both to $R_k$ and to $f_k^{-1}(R'_2)$ then it always lies close to
a point of $\mathcal{I}'_k$.

Property (2) of quasiholomorphic coverings follows from the observation
that, since the differentials of $f_k$ and $V'_2$ both have complex rank at 
least $1$ everywhere, $\nabla f'_{2k}(x)$ can only be small if the Jacobians
of $f_k$ at $x$ and of $V'_2$ at $f_k(x)$ are both small. These
quantities vanish transversely ($f_k$ is quasiholomorphic and $V'_2$ is
generic), so $x$ must lie close to both branch curves, and hence, by the
above observation, close to $\mathcal{I'}_k$ (closer than $d_0$ if 
$|\nabla f'_{2k}(x)|$ is assumed small enough). 
In fact, $|\nabla f'_{2k}|$ remains bounded away from
$0$ even near $\mathcal{I}'_k$, because, as observed in the proof of Lemma 1,
property (5) of Definition 7 implies that $V_2^0$ (and hence also $V'_2$)
restricts to the branch curve of $f_k$ as an immersion.

We now turn to the third property. The $(2,0)$-Jacobian of $f'_{2k}$ is
given by $\mathrm{Jac}(f'_2k)=\mathrm{Jac}(f_k)\cdot f_k^*\mathrm{Jac}(V'_2)$.
It can only be small when one of the two terms in the product is small,
i.e.\ near one of the two branch curves. Moreover, $f_k^*\mathrm{Jac}(V'_2)$
is bounded away from zero everywhere except near $f_k^{-1}(R'_2)$, so the 
transverse vanishing of $\mathrm{Jac}(f_k)$ implies that of 
$\mathrm{Jac}(f'_{2k})$ at these points. Similarly $\mathrm{Jac}(f_k)$ is
bounded from below everywhere except near $R_k$, so the transverse vanishing
of $f_k^*\mathrm{Jac}(V'_2)$ implies the desired property at these points.
As a consequence the transversality to $0$ of $\mathrm{Jac}(f'_{2k})$ holds
everywhere except near $\mathcal{I}'_k$ (note that the obtained 
transversality estimate has to be decreased when $d_0$ becomes smaller).

We now look at property (4). Away from $\mathcal{I}'_k$ the branch curve of
$f'_{2k}$ consists of two separate components, $R_k$ and $f_k^{-1}(R'_2)$,
so we work separately on each component. On $R_k-\mathcal{I}'_k$, we
know that $\partial(f_{k|R_k})$ is uniformly transverse to $0$, and because
$\mathcal{I}'_k$ has been removed the complex linear map $\nabla V'_2$ is an
isomorphism at every point of the image, with norm bounded from below
(the constant depends on $d_0$). Composing $\partial(f_{k|R_k})$
with $\nabla V'_2$, we obtain that $\partial(f'_{2k|R_k})$ is also uniformly 
transverse to $0$ at all points of $R_k$ at distance more than $d_0$ from 
$\mathcal{I}'_k$ (again, the constant depends on $d_0$). The argument works
similarly on $f_k^{-1}(R'_2)-\mathcal{I}'_k$~: $\partial f_k$ is an
isomorphism with norm bounded from below (the constant depends on $d_0$),
and because $V'_2$ has been chosen generic the quantity
$\nabla(V'_{2|R'_2})$ vanishes transversely, so
$\partial(f'_{2k|f_k^{-1}(R'_2)})$ is uniformly transverse to $0$ at all
points of $f_k^{-1}(R'_2)$ at distance more than $d_0$ from $\mathcal{I}'_k$.

Observe by the way that all cusp points of $f_k$ and of $V'_2$ lie away
from $\mathcal{I}'_k$. Indeed, for the cusp points of $f_k$ it follows from
property (4) in Definition 7 that they lie away from the branch curve of
$V_2^0$ and hence from that of $V'_2$, as observed in the proof of Lemma 1.
On the other hand, is easy to see that the cusp points of $V'_2$ all lie 
close to one the three singular points of $V_2^0$, while property (3) in 
Definition 7 implies that the branch curve of $f_k$ remains far away from 
these points.

Property (5) is easy to check: since compatible almost-complex 
structures on $X$ are sections of a bundle with contractible fiber, it is 
sufficient to work locally near a cusp point. The points we have to consider 
are either cusp points of $f_k$ or the preimages by $f_k$ of those of $V'_2$.
In the first case, it is sufficient to choose the same almost-complex
structure $\tilde{J}_k$ as for $f_k$, because $V'_2$ is holomorphic. In the
second case, consider the pull-back $f_k^*\mathbb{J}_0$ of the standard 
complex structure of $\CP^2$ via the map $f_k$. Since all cusp points of 
$V'_2$ lie far from the branch curve of $f_k$, the differential of $f_k$ 
is locally an isomorphism and satisfies a uniform lower bound. Therefore the
asymptotic holomorphicity of the sections defining $s_k$ is enough to
ensure that $f_k^*\mathbb{J}_0$ differs from $J$ by at most $O(k^{-1/2})$ in
any $C^r$ norm. A standard argument involving a smooth cut-off function can be
used in order to define a smooth almost-complex structure which coincides
with $f_k^*\mathbb{J}_0$ near the cusp point and with $J$ outside a small
ball.

We now turn to property (6). Consider a point $x\in X$ where the first two
sections defining $f'_{2k}$, namely $Q_0(s_k)$ and $Q_1(s_k)$, are both very
small. Because the quadratic map $V'_2$ is close to $V_2^0$, and
because the only preimage of $(0:0:1)$ by $V_2^0$ is $(0:0:1)$ itself,
the quantities $s_k^0(x)$ and $s_k^1(x)$ are also small. So, if we assume
that $|V'_2-V_2^0|$ is sufficiently small, the uniform transversality
property of $(s_k^0,s_k^1)$ provides a lower bound on $\mathrm{Jac}(f_k)(x)$.
On the other hand, if $V'_2$ is chosen generic, then its branch
curve avoids the point $(1:0:0)$ by a certain distance $\rho>0$. Therefore,
if $Q_0(s_k)$ and $Q_1(s_k)$ are sufficiently small,
then $f'_{2k}(x)$ lies at distance at least $\rho/2$ from the branch curve
of $V'_2$. and we can obtain a uniform lower bound (depending on $\rho$
only) on the Jacobian of $V'_2$ at $f_k(x)$. It follows that
$\mathrm{Jac}(f'_{2k})(x)=\mathrm{Jac}(f_k)(x)
\,\mathrm{Jac}(V'_2)(f_k(x))$ is bounded from below by a fixed constant
independently of $k$. Because of the $C^1$ bounds on $Q_i(s_k)$, 
we conclude that the covariant derivative of
$(Q_0(s_k),Q_1(s_k))$ at $x$ is surjective and bounded from below by a
uniform constant. So property (6) holds.

We finally look at property (7), which actually is equivalent to the
requirement that the branch curve be locally braided. Most 
of the work has already been done in the proof of Lemma 1. More precisely, 
after removing the intersection $\mathcal{I}'_k$, the branch curve of 
$f'_{2k}$ splits into the two components $R_k$ and $f_k^{-1}(R'_2)$, 
and we consider them separately. The critical points of $\psi^0_k=(\pi\circ
V_2^0\circ f_k)_{|R_k}$ and $\psi'_k=(\pi\circ f'_{2k})_{|R_k}$ correspond
to the cusps and tangency points of $V_2^0(D_k)$ and $V'_2(D_k)$,
respectively. Therefore, we have seen in the proof of Lemma 1 
that all the critical points of $\psi_k^0$, and hence those of $\psi'_k$,
are non-degenerate, with a uniform estimate; moreover, 
they all lie in a neighborhood of 
$\mathcal{C}_k\cup\mathcal{T}_k\cup \mathcal{I}_k$, 
which implies that $f_k$ is locally holomorphic with respect to a suitable 
almost-complex structure.

We now look at the component $f_k^{-1}(R'_2)$ away from the points of
$\mathcal{I}'_k$~: since $f_k$ is a local diffeomorphism at all such points,
the expected uniform transversality of $\partial(\pi\circ f'_{2k})$ is equivalent
to the same property for $\partial(\pi\circ V'_2)$ restricted to $R'_2$.
However it is easy to check that such a transversality property holds as 
soon as $V'_2$ is chosen generic (actually, as soon as $V'_2(R'_2)$ is 
locally braided). Of course the transversality estimate on
$\partial(\pi\circ f'_{2k})$ depends on the distance $d_0$, because a lower
bound on $\partial f_k$ is used when lifting the transversality property
from $\pi\circ V'_2$ to $\pi\circ f'_{2k}$.
Also observe that the holomorphicity of $V'_2$ implies that
the differential of $\pi\circ V'_{2|R'_2}$ vanishes completely 
at the tangency points of the branch curve of $V'_2$ (these are genuine
tangencies); this clearly implies the same property for $\pi\circ f'_{2k}$
at the tangency points coming from $f_k^{-1}(R'_2)$. This concludes the
proof.
\endproof

Proposition 2 implies that we can proceed in the following way to construct
quasiholomorphic coverings given by sections of $L^{\otimes 2k}$ for large
$k$~: first construct stably quasiholomorphic coverings $f_k$ as given
by Proposition~\ref{prop:stqh}~; then, define $f'_{2k}=V'_2\circ f_k$ for
a generic perturbation $V'_2$ of $V_2^0$~; and finally perturb $f'_{2k}$
in order to get quasiholomorphic coverings.

Following the arguments in \cite{A2} and \cite{AK} (see also \cite{Agokova}
and the argument in \S 2.3 below), we can make the following
observations concerning the process by which the maps $f'_{2k}$ are perturbed 
and made quasiholomorphic. The first step of the construction of
quasiholomorphic coverings is to ensure that all the required uniform
transversality properties are satisfied over all of $X$. This process is
a purely local iterative construction, so that when one starts with
$f'_{2k}$ it is sufficient to perturb the given sections of $L^{\otimes 2k}$
near the points of $\mathcal{I}'_k$, or equivalently near the points of
$\mathcal{I}_k$~; the required perturbation can be chosen smaller than
any fixed given constant (independent of $k$), so that it does not 
significantly affect the topology of $f'_{2k}$ 
away from the points of $\mathcal{I}'_k$. The next step in order to
construct quasiholomorphic coverings is to ensure property $(5)$ of
Definition~5 at the cusp points as well as the last requirement of 
property $(7)$ at the tangency points~; since the necessary perturbation 
is bounded by a fixed multiple of $k^{-1/2}$, it has no effect whatsoever
on braid monodromy outside of a fixed small
neighborhood of $\mathcal{I}'_k$.

At this point in the construction, the branch curves are 
already locally braided and therefore have well-defined braid monodromies 
up to $m$-equi\-va\-lence~; ensuring the remaining conditions $(3')$ and $(8)$ 
has no effect on the monodromy data. More precisely, the self-transversality
of the branch curves (condition $(8)$) is obtained by an arbitrarily small
perturbation, which is precisely how one defines the braid factorization
associated to a locally braided curve. Meanwhile, condition $(3')$
is obtained by a perturbation process which does not affect the branch 
curve (see \cite{AK}).
Finally, notice that, once the covering maps $f'_{2k}$ are
perturbed and made quasiholomorphic, the braid monodromy invariants
associated to them must coincide with those associated to $f_{2k}$, at
least provided that $k$ is large enough~: this is a direct consequence
of the uniqueness result of \cite{AK}.

As a consequence of these observations, by computing the braid factorization
corresponding to the branch curve of $f'_{2k}$ (very singular, with 
components of large multiplicity), a great step towards computing the
braid factorization for $f_{2k}$ is already accomplished~:
the only remaining task is to understand the effect on braid factorizations
of the perturbation performed near the points of $\mathcal{I}'_k$.
This justifies the strategy of proof used in~\S 3.

\subsection{Proof of Proposition \ref{prop:stqh}}

Proposition \ref{prop:stqh} can be proved using the same techniques as in
\cite{A2} and \cite{AK} (see also \cite{Agokova})~; however, the result
of \cite{Aseoul} can be used to greatly simplify the argument. Observe
that the properties expected of $s_k$ are of two types~: on one hand,
uniform transversality properties, which are open conditions on the 
holomorphic part of the jet of $s_k$, and on the other hand, compatibility
properties, involving the vanishing of certain antiholomorphic derivatives 
along the branch curve. The proof therefore consists of two parts.
In the first part, successive perturbations of $s_k$ are performed in order 
to achieve the various required transversality properties~; each perturbation
is chosen small enough in order to preserve the previously obtained
transversality properties. In the second part, $s_k$ is perturbed along
the curve $R_k$ by at most a fixed multiple of $k^{-1/2}$ in order to
obtain the compatibility conditions.

The first part of the argument can be either carried out as in \cite{A2}
and \cite{AK}, or more efficiently by using the result of \cite{Aseoul}
in the following manner.

Let $E_k=\C^3\otimes L^{\otimes k}$, and consider the holomorphic jet bundles
$\mathcal{J}^2E_k=E_k\oplus T^*X^{(1,0)}\otimes E_k\oplus
(T^*X^{(1,0)})^{\otimes 2}_{\mathrm{sym}}\otimes E_k$. We define the
holomorphic $2$-jet $j^2s$ of a section $s\in\Gamma(E_k)$ as $(s,\partial
s,\partial(\partial s)_\mathrm{sym})$, discarding the antiholomorphic
terms or the antisymmetric part of $\partial\partial s$ (these terms
are bounded by $O(k^{-1/2})$ for asymptotically holomorphic sections). 
Recall from \cite{Aseoul} the notion of finite Whitney quasi-stratification 
of a jet bundle~:

\begin{defi}
Let $(A,\prec)$ be a finite set carrying a binary relation without cycles 
(i.e., $a_1\prec\dots\prec a_p\Rightarrow a_p\not\prec a_1$).
A finite Whitney quasi-stratification of $\mathcal{J}^2 E_k$ indexed by $A$
is a collection $(S^a)_{a\in A}$ of smooth submanifolds of $\mathcal{J}^2
E_k$, transverse to the fibers, not necessarily mutually disjoint, with the 
following properties~: $(1)$ 
$\partial S^a=\overline{S^a}-S^a\subseteq \bigcup_{b\prec a} S^b$~;
$(2)$ given any point $p\in\partial S^a$, there exists $b\prec a$ such
that $p\in S^b$ and such that either $S^b\subset\partial S^a$ and the 
Whitney regularity condition is satisfied at all points of $S^b$,
or $p\not\in\Theta_{S^b}$, where $\Theta_{S^b}\subset S^b$ 
is the set of points where the $2$-jet of a section of $E_k$ can intersect 
$S^b$ transversely $($in particular $\Theta_{S^b}=\emptyset$ whenever
$\mathrm{codim}_\C\,S^b>2)$.
\end{defi}

As in \cite{Aseoul},
say that a sequence of finite Whitney quasi-stratifications $\mathcal{S}_k$ 
of $\mathcal{J}^2E_k$ is asymptotically holomorphic if all the strata are
approximately holomorphic submanifolds of $\mathcal{J}^2 E_k$, with uniform
bounds on the curvature of the strata and on their transversality to the 
fibers of $\mathcal{J}^2 E_k$.

It was shown in \cite{Aseoul} that, given asymptotically holomorphic finite
Whitney quasi-stratifications $\mathcal{S}_k$ of $\mathcal{J}^2 E_k$, it
is always possible for large enough $k$ to construct asymptotically
holomorphic sections of $E_k$ whose $2$-jets are uniformly transverse to
the strata of $\mathcal{S}_k$~; moreover, these sections can be chosen 
arbitrarily close to any given asymptotically holomorphic sections of
$E_k$. The result also holds for one-parameter families of
sections, which implies that the constructed sections are, for large $k$,
canonical up to isotopy. 

Using local
approximately holomorphic sections of $L^{\otimes k}$ and coordinates over $X$,
the fibers of $\mathcal{J}^2 E_k$ can be identified with the space
$\mathcal{J}^2_{2,3}$ of jets of holomorphic maps from
$\C^2$ to $\C^3$. It was observed in \cite{Aseoul} that, if a sequence
of finite Whitney quasi-stratifications of $\mathcal{J}^2 E_k$ is such that
by this process
the restrictions of $\mathcal{S}_k$ to the fibers of $\mathcal{J}^2 E_k$ 
are all identified with a fixed given finite Whitney quasi-stratification
of $\mathcal{J}^2_{2,3}$ by complex submanifolds, then the
quasi-stratifications $\mathcal{S}_k$ are asymptotically holomorphic.

We define finite Whitney quasi-stratifications of $\mathcal{J}^2E_k$ in the 
following way. Consider the symmetric holomorphic part $j^2s(x)$ of the 
$2$-jet of a section $s=(s^0,s^1,s^2)\in\Gamma(E_k)$ at a point $x\in X$~; 
if $s(x)\neq 0$, denote by $f$ the corresponding $\CP^2$-valued map, and by 
$\phi^i$ ($i\in\{0,1,2\}$) its projections to $\CP^1$ along coordinate axes 
if they are well-defined. Finally, if 
$\mathrm{Jac}\,f(x)=\wedge^2\partial f(x)=0$ 
and $\partial\mathrm{Jac}\,f(x)_\mathrm{sym}=(\partial\partial
f(x))_\mathrm{sym}\wedge\partial f(x)\neq 0$, call $R_x$ the kernel of the
$(1,0)$-form $\partial\mathrm{Jac}\,f(x)_\mathrm{sym}$~;
one easily checks that $R_x$ is well defined in terms of $j^2 s$ only and 
that it differs from the tangent space at $x$ to the ramification curve 
of $f$ by at most $O(k^{-1/2})$.
We define the following submanifolds of $\mathcal{J}^2 E_k$ (in the
last two definitions, $\{i,j,k\}=\{0,1,2\}$)~:
\begin{alignat*}{2}
Z=&\{j^2s(x),\ s(x)=0\}\qquad&
(\mathrm{codim.}\ 3)&\\
Z_{ij}=&\{j^2s(x),\ s^i(x)=s^j(x)=0\}\qquad&
(\mathrm{codim.}\ 2)&\\
Z_{i}=&\{j^2s(x),\ s^i(x)=0\}\qquad\qquad\qquad\qquad\qquad\qquad&
(\mathrm{codim.}\ 1)&\\
\Sigma^2=&\{j^2s(x),\ s(x)\neq 0,\ \partial f(x)=0\}\qquad&
(\mathrm{codim.}\ 4)&\\
\Sigma^1=&\{j^2s(x)\not\in Z,\ \partial f(x)\neq 0,\ \mathrm{Jac}\,f(x)=0\}\ &
(\mathrm{codim.}\ 1)&\\
\Sigma^1_s=&\{j^2s(x)\in\Sigma^1,\ \partial\mathrm{Jac}\,f(x)_\mathrm{sym}=0\}\ &
(\mathrm{codim.}\ 3)&\\
\Sigma^{1,1}=&\{j^2s(x)\in\Sigma^1-\Sigma^1_s,\ \partial f(x)_{|R_x}=0\}\ &
(\mathrm{codim.}\ 2)&\\
\Sigma^1_t=&\{j^2s(x)\in\Sigma^1-Z_{01},\ \partial \phi^2(x)=0\}\ &
(\mathrm{codim.}\ 2)&\\
\Sigma^{1,1}_t=&\Sigma^{1,1}\cap\Sigma^1_t\ &
(\mathrm{codim.}\ 3)&\\
S_i=&\{j^2s(x)\in Z_i-Z_{jk},\ \partial \phi^i(x)=0\}\ &
(\mathrm{codim.}\ 3)&\\
S'_i=&\{j^2s(x)\in Z_i-Z_{jk},\ \partial s^i(x)\neq 0,\ 
\partial \phi^i(x)_{|\mathrm{Ker}\,\partial s^i(x)}=0\}\ &
(\mathrm{codim.}\ 2)&
\end{alignat*}

One easily checks that all these subsets are smooth submanifolds of
$\mathcal{J}^2 E_k$. Moreover, $Z$, $Z_i$ and $Z_{ij}$ are closed~;
$\partial \Sigma^2\subseteq Z$~; $\partial \Sigma^1$ and 
$\partial \Sigma^1_s$ are contained in $\Sigma^2\cup Z$~;
$\partial \Sigma^{1,1}\subseteq \Sigma^1_s\cup\Sigma^2\cup Z$~;
$\partial \Sigma^1_t\subseteq \Sigma^2\cup Z\cup
(Z_{01}-\Theta_{Z_{01}})$~; $\partial \Sigma^{1,1}_t\subseteq \Sigma^1_s
\cup\Sigma^2\cup Z\cup (Z_{01}-\Theta_{Z_{01}})$~;
$\partial S_i\subseteq (Z_{jk}-\Theta_{Z_{jk}})$~;
$\partial S'_i\subseteq (Z_{jk}-\Theta_{Z_{jk}})\cup (Z_i-\Theta_{Z_i})$.
Therefore, these submanifolds define quasi-stratifications $\mathcal{S}_k$ 
of $\mathcal{J}^2E_k$. 
Note that, because $\Sigma^1_s=\Sigma^1-\Theta_{\Sigma^1}$, 
the stratum $\Sigma^1_s$ can in fact be eliminated from this description.
Moreover, if one uses local approximately holomorphic coordinates and
asymptotically holomorphic sections of $L^{\otimes k}$ to trivialize 
$\mathcal{J}^2E_k$, it is easy to see that the resulting picture is the
same above every point of $X$~: the submanifolds in $\mathcal{S}_k$ are 
identified with holomorphic submanifolds of $\mathcal{J}^2_{2,3}$ defined
by the same equations. Therefore, by \cite{Aseoul} the quasi-stratifications
$\mathcal{S}_k$ are asymptotically holomorphic.

It is easy to see that conditions $(1)$, $(2)$, $(3)$, $(4)$ and $(6)$ of
Definition 5 are equivalent to the uniform transversality of $j^2 s_k$ to
$Z$, $\Sigma^2$, $\Sigma^1$, $\Sigma^{1,1}$ and $Z_{01}$, respectively. 
Similarly, conditions $(2)$ and $(3)$ of Definition 7 correspond to the
uniform transversality of $j^2 s_k$ to $Z_i$ and $Z_{ij}$ respectively.
Observing that $\partial(\phi_{k|R_k})$ can only vanish at a point $x\in
R_k$ if either $\partial\phi_k(x)=0$ or $\partial(f_{k|R_k})$ vanishes
at $x$, we can rephrase condition $(7)$ of Definition 5 in terms of
uniform transversality to the singular submanifold of $\mathcal{J}^2 E_k$
consisting of the union of $\Sigma^{1,1}$ (cusp points) and $\Sigma^1_t$
(tangencies), intersecting regularly along $\Sigma^{1,1}_t$ (``vertical''
cusp points). Therefore, it is equivalent to the uniform transversality
of $j^2 s_k$ to $\Sigma^{1,1}$, $\Sigma^1_t$, and $\Sigma^{1,1}_t$.
Finally, conditions $(4)$ and $(5)$ of Definition 7 correspond to the
uniform transversality of $j^2 s_k$ to $S'_i$ and $S_i$ respectively.

So, the uniform transversality of $j^2 s_k$ to the quasi-stratifications
$\mathcal{S}_k$, as given by the main result of \cite{Aseoul} provided
that $k$ is large enough, is equivalent
to the various transversality requirements listed in Definitions 5 and 7.
Moreover, the sections of $\C^3\otimes L^{\otimes k}$ constructed in this
manner are canonical up to isotopy, as follows from Theorem 3.2 of 
\cite{Aseoul}~: given any two sequences of such sections, it is possible
for large enough $k$ to find one-parameter families of sections of 
$\C^3\otimes L^{\otimes k}$ interpolating between them and enjoying the 
same uniform transversality properties for all parameter values.

We now turn to the second part of the argument, namely obtaining the
other required properties by perturbing the sections $s_k$ by at most
$O(k^{-1/2})$, which clearly affects neither holomorphicity nor
transversality properties. The argument is exactly the same as in
\cite{AK}~; the only difference is that the set $\mathcal{F}_k$ of points 
where the map $f_k$ must made holomorphic with respect to a slightly perturbed
almost-complex structure is now slightly larger~: one now sets
$\mathcal{F}_k=\mathcal{C}_k\cup\mathcal{T}_k\cup\mathcal{I}_k$ instead
of $\mathcal{F}_k=\mathcal{C}_k$.

As in \cite{A2} and \cite{AK}, one first chooses suitable almost-complex
structures $\tilde{J}_k$ differing from $J$ by $O(k^{-1/2})$ and integrable
near the finite set $\mathcal{F}_k$. It is then possible to
perturb $f_k$ near these points in order to obtain
condition $(5)$ of Definition 5, by the same argument as in \S 4.1
of \cite{A2}. Next, a generic small perturbation 
yields the self-transversality of $D_k$ (property $(8)$ of Definition~5).
Finally, a suitable perturbation of $f_k$, supported near $R_k$ and
vanishing near the points of $\mathcal{F}_k$, yields property $(3')$ of 
Definition 5 along the branch curve, without modifying $R_k$ and $D_k$,
and therefore without affecting the previously obtained compatibility 
properties. As shown in \cite{AK} these various constructions can be
performed in one-parameter families, except for property $(8)$ of Definition
5 where cancellations of pairs of nodes must be allowed~; this yields the
desired result of uniqueness up to isotopy, and completes the proof of
Proposition 1.

\section{The degree doubling formula for braid monodromies}

\subsection{Generalities about the braid group}\label{ss:generalities}
We begin by recalling general definitions and notations concerning the braid
group on $d$ strings. Consider a set $P=\{p_1,\ldots,p_d\}$ of $d$ points in
the plane, and recall that $B_d=\pi_0\,\mathrm{Diff}^+_c(\R^2,P)$ is by
definition the group of equivalence classes of compactly supported
orientation-preserving diffeomorphisms of the plane which leave invariant
the set $P$, where two diffeomorphisms are equivalent if and only if they
induce the same automorphism of $\pi_1(\R^2-P)$. 
Equivalently $B_d$ can be considered as the fundamental group
of the configuration space of $d$ points in the plane~: a braid corresponds
to a motion of the points $p_1,\ldots,p_d$ such that they remain distinct
at all times and eventually return to their original positions (but
possibly in a different order) up to homotopy. An important subgroup of
$B_d$ is the group of pure braids $P_d$ (the braids which preserve each
of the points $p_1,\ldots,p_d$ individually)~; it is clear that 
$B_d/P_d$ is the symmetric group $S_d$.

We will place the points $p_1,\ldots,p_d$ in that order on the real axis,
and denote by $X_i$ the positive (counterclockwise) half-twist along the
line segment joining $p_i$ to $p_{i+1}$, for each $1\le i\le d-1$. It is a
classical fact that $B_d$ is generated by the $d-1$ half-twists $X_i$, and
that the relations between them are $X_iX_j=X_jX_i$ whenever $|i-j|>1$ and
$X_iX_{i+1}X_i= X_{i+1}X_iX_{i+1}$. The center of the braid group is
generated by the element $\Delta_{d}^{2}=(X_{1}\ldots X_{d-1})^{d}$,
which corresponds to rotating everything by $2\pi$.

We will be especially interested in the half-twists
$$Z_{ij}=X_{j-1}\cdot\ldots\cdot X_{i+1}\cdot X_i\cdot X_{i+1}^{-1}\cdot
\ldots\cdot X_{j-1}^{-1}\qquad (1\le i<j\le d).$$
The braid $Z_{ij}$ is a positive half-twist along a path joining the points
$p_i$ and $p_j$ and passing {\it above} all the points inbetween~:

\begin{center}
\setlength{\unitlength}{0.25in}
\begin{picture}(9,2)(-1,-0.7)
\multiput(0,0)(1,0){3}{\circle*{0.1}}
\multiput(5,0)(1,0){3}{\circle*{0.1}}
\multiput(3,0)(0.2,0){6}{\circle*{0.02}}
\put(0,-0.2){\makebox(0,0)[ct]{1}}
\put(1,-0.2){\makebox(0,0)[ct]{$i$}}
\put(6,-0.2){\makebox(0,0)[ct]{$j$}}
\put(7,-0.2){\makebox(0,0)[ct]{$d$}}
\qbezier(1,0)(1.6,0.6)(3.5,0.6)
\qbezier(6,0)(5.4,0.6)(3.5,0.6)
\put(3.5,0.8){\makebox(0,0)[cb]{$Z_{ij}$}}
\end{picture}
\end{center}

Note in particular that $Z_{i,i+1}=X_i$ and that $Z_{ij}$ commutes with
$Z_{kl}$ whenever $i<j<k<l$ or $i<k<l<j$. Other useful relations are
$Z_{ij}Z_{ik}=Z_{ik}Z_{jk}=Z_{jk}Z_{ij}$ whenever $i<j<k$ (these three
expressions differ by a Hurwitz move).
\medskip

The following factorization of $\Delta^2$ as a product of half-twists 
corresponds to the braid monodromy of a smooth curve of degree $d$ in 
$\CP^2$ (see \cite{MOI81})~:
$$\Delta_{d}^{2}=(X_{1}\ldots X_{d-1})^{d}.$$
Another important factorization is
\begin{equation}\label{eq:d_lines}
\Delta_d^2 = \prod_{i=1}^{d-1}\prod_{j=i+1}^{d} Z_{ij}^2 =
\prod_{i=2}^{d}\prod_{j=1}^{i-1} Z_{ji}^2
\end{equation}
(these two expressions are clearly Hurwitz equivalent). This factorization 
corresponds to the braid monodromy of a union of $d$ lines in generic 
position (see \cite{MOI81}).
\medskip

We now turn to geometric monodromy representations. Consider the branch
curve $D$ of an $n$-sheeted branched covering over $\CP^2$, and fix
geometric generators $\gamma_1,\dots,\gamma_d$ of $\pi_1(\CP^2-D)$ 
(small loops going around the $d=\deg D$ intersection points of $D$ with 
a given generic fiber of the projection $\pi$). It is then possible to
define as in \S 1.1 the geometric monodromy representation 
$\theta:F_d\to S_n$ associated to the covering. As observed in \S 1.1,
the fact that the product $\gamma_1\cdot\dots\cdot\gamma_d$ is trivial in 
$\pi_1(\CP^2-D)$ implies that the product of the $d$ transpositions 
$\theta(\gamma_1),\dots,\theta(\gamma_d)$ in $S_n$ is also trivial, and
the connectedness of the considered covering of $\CP^2$ implies that these
transpositions act transitively on $\{1,\dots,n\}$ and hence generate $S_n$.

It is a well-known fact that any two factorizations of the identity element
in $S_n$ as a product of the same number of transpositions generating $S_n$
are equivalent by a succession of Hurwitz moves (this can be seen e.g.\ by
comparing the two corresponding $n$-sheeted simple branched covers of 
$\CP^1$). Therefore, after a
suitable reordering of the sheets of the covering $\pi:D\to\CP^1$ 
(which amounts to a global conjugation of the braid factorization), one 
may freely assume that the permutations $\theta(\gamma_i)$ are equal to
certain predetermined transpositions. Our choice of transpositions 
in the case of the branch curve of $f_k$ will be made explicit in 
\S \ref{ss:braidthm}.

\subsection{The folding process}\label{ss:folding}
We now compute the braid monodromy of the curve $V'_2(D_k)$, where $D_k$ is
the branch curve of one of the stably quasiholomorphic maps $f_k$ given
by Proposition 1 and $V'_2$ is a generic perturbation of $V^0_2$ as in 
\S 2.2. The idea is to use Lemma 1 to reduce oneself to the easy case where
$D_k$ is a union of $d=\deg D_k$ lines through a point in $\CP^2$. In that 
case, $V'_2(D_k)$ becomes a union of $d$ conics through a point, and its
braid monodromy can be computed explicitly. The result is the following:

\begin{prop}\label{prop:folding}
The braid factorization corresponding to the curve $V'_2(D_k)$ is given by
the formula
\begin{equation}\label{eq:folding0}
\Delta_{2d}^2=\Bigl(\prod_{i=1}^{d-1}\prod_{j=i+1}^d Z_{i'j'}^2\Bigr)\cdot
\prod_{i=1}^d Z_{ii'}\cdot F_k\cdot
\Bigl(\prod_{i=1}^{d-1}\prod_{j=i+1}^d Z_{i'j'}^2\Bigr)^2\cdot 
\prod_{i=1}^d Z_{ii'},
\end{equation}
or equivalently
\begin{equation}\label{eq:folding}
\Delta_{2d}^2=\prod_{i=1}^d \hat{Z}_{ii'}\cdot F_k\cdot
\Bigl(\prod_{i=1}^{d-1}\prod_{j=i+1}^d Z_{i'j'}^2\Bigr)^3\cdot 
\prod_{i=1}^d Z_{ii'},
\end{equation}
where $F_k$ is the image of the braid factorization for $D_k$ via the
embedding of the braid group $B_d$ in $B_{2d}$ obtained by considering a ball
containing only the first $d$ points, and $\hat{Z}_{ii'}$ is a half-twist 
along the following path~:
\begin{center}
\setlength{\unitlength}{0.25in}
\begin{picture}(11.5,2)(-1,-1)
\put(0,0){\circle*{0.1}}
\put(2,0){\circle*{0.1}}
\put(4,0){\circle*{0.1}}
\put(5,0){\circle*{0.1}}
\put(7,0){\circle*{0.1}}
\put(9,0){\circle*{0.1}}
\multiput(0.5,0)(0.2,0){6}{\circle*{0.02}}
\multiput(2.5,0)(0.2,0){6}{\circle*{0.02}}
\multiput(5.5,0)(0.2,0){6}{\circle*{0.02}}
\multiput(7.5,0)(0.2,0){6}{\circle*{0.02}}
\put(0,-0.2){\makebox(0,0)[ct]{$1$}}
\put(2,-0.2){\makebox(0,0)[ct]{$i$}}
\put(4,-0.2){\makebox(0,0)[ct]{$d$}}
\put(5,-0.2){\makebox(0,0)[ct]{$1'$}}
\put(7,-0.2){\makebox(0,0)[ct]{$i'$}}
\put(9,-0.2){\makebox(0,0)[ct]{$d'$}}
\qbezier[150](2,0)(2,1)(4,1)
\put(4,1){\line(1,0){4}}
\qbezier[150](8,1)(10,1)(10,0)
\qbezier[150](10,0)(10,-1)(8,-1)
\put(6,-1){\line(1,0){2}}
\qbezier[110](6,-1)(4.5,-1)(4.5,-0.2)
\qbezier[100](4.5,-0.2)(4.5,0.6)(5.7,0.6)
\qbezier[100](5.7,0.6)(7,0.6)(7,0)
\end{picture}
\end{center}
\end{prop}

Equation (\ref{eq:folding}) is an identity in the braid group $B_{2d}$ acting
on $2d$ points labelled $1,\dots,d,1',\dots,d'$ (each pair $i,i'$
corresponds to one of the $d$ conics).

Consider as in Lemma 1 the linear contraction map
$\psi_a:(x\!:\!y\!:\!z)\mapsto
(x\!:\!ay+(1-a)x\!:\!az+(1-a)x)$. When $a$ converges to $0$, the images of
all the points outside of the line $L_0:\{x=0\}$ converge towards the point
$p_0=(1\!:\!1\!:\!1)$. Since $\psi_a$ maps fibers of $\pi$ to fibers of 
$\pi$, the curves $\psi_a(D_k)$ are braided for all values of $a$.
Moreover, $\psi_a$ restricts to the line $L_0:\{x=0\}$ as the identity, and
when $a\to 0$ the image of any line intersecting $L_0$ transversely at a point 
$p=(0\!:\!y\!:\!z)$ converges to the line through 
$p$ and $p_0$.

By an arbitrarily small perturbation, and without losing the other
properties of $D_k$, we can easily assume that the point $(0\!:\!1\!:\!1)$
does not belong to $D_k$, and that none of the nodes of $D_k$ lies on $L_0$.
Therefore, by Lemma~1 the curve $V_2^0(\psi_a(D_k))$ is locally braided
for sufficiently small $a\neq 0$, and isotopic to $V'_2(D_k)$ through
locally braided curves. This implies that the braid factorizations for 
$V_2^0(\psi_a(D_k))$ and for $V'_2(D_k)$ are $m$-equivalent (in fact, 
when $D_k$ is a complex curve the isotopy can be carried out inside the
complex category, so in that case the braid factorizations are even 
Hurwitz and conjugation equivalent).

When $a$ is sufficiently close to $0$, outside
of a small ball centered at $p_0$ the curve $\psi_a(D_k)$ is arbitrarily
close to the union of $d=\deg D_k$ lines joining the points of $D_k\cap L_0$
with $p_0$, and by construction the images by $V_2^0$ of these
$d$ lines are distinct non-degenerate conics in $\CP^2$. Moreover, the
restriction of $V_2^0$ to a neighborhood of $p_0$ is
a diffeomorphism mapping fibers of $\pi$ to fibers of $\pi$.
Therefore, the braid factorization of $V_2^0(\psi_a(D_k))$, or equivalently
that of $V'_2(D_k)$ 
can be obtained by plugging the braid factorization of $D_k$ into
the formula for the braid monodromy of a union of $d$ conics
passing through the point $p_0$, i.e.\ by deleting a neighborhood of $p_0$
from this configuration and replacing it with a braided curve isotopic to
the affine part of $D_k$ (suitably rescaled into a small ball).

As a first step, we therefore need to compute the braid monodromy of a
union of $d$ conics passing through $p_0$. Observe that any 
configuration of $d$ non-degenerate conics in $\CP^2$ intersecting each
other transversely at $p_0$ gives rise to a well-defined braid 
factorization as soon as none of the conics passes 
through the pole of the projection $\pi$~: any such configuration is a
locally braided curve, and can be perturbed into a braided curve (a union
of conics in general position) by an arbitrarily small perturbation.
The connectedness of the space of configurations of conics implies
that, up to Hurwitz and conjugation equivalence, it does not actually 
matter which conics are used for the computation of the braid monodromy.

Following Moishezon, the calculation can be carried out by simultaneously 
``degenerating'' all the conics to pairs of lines, i.e.\ by considering a
limit configuration where each of the conics is very close to a union of two
lines \cite{MOI}. However, for the purpose of proving Theorem 2 it is more 
efficient to perform a direct calculation using a specific
configuration of conics. We consider $d$ conics with real coefficients,
intersecting at the point $p_0$, and with their other mutual intersections
lying close to three given points $p_1$, $p_2$, $p_3$, as in the following
diagram (representing the intersection of the configuration with
$\R^2\subset\C^2$, with the fibers of $\pi$ corresponding to vertical lines).

\begin{center}
\setlength{\unitlength}{0.25in}
\begin{picture}(8,4.5)(-4,-0.5)
\qbezier[300](0.4,0)(-3.1,-0.5)(-3.6,1)
\qbezier[300](-3.6,1)(-4.1,2.5)(-0.6,3)
\qbezier[300](-0.6,3)(2.9,3.6)(3.4,2)
\qbezier[300](3.4,2)(3.9,0.5)(0.4,0)
\qbezier[300](-0.2,0)(3.3,-0.5)(3.8,1)
\qbezier[300](3.8,1)(4.3,2.5)(0.8,3)
\qbezier[300](0.8,3)(-2.7,3.6)(-3.2,2)
\qbezier[300](-3.2,2)(-3.7,0.5)(-0.2,0)
\qbezier[300](-3.4,1.5)(-3.4,3.12)(0.1,3.14)
\qbezier[300](-3.4,1.5)(-3.4,-0.045)(0.1,-0.045)
\qbezier[300](3.6,1.5)(3.6,3.12)(0.1,3.14)
\qbezier[300](3.6,1.5)(3.6,-0.045)(0.1,-0.045)
\multiput(-1.2,-0.5)(0,0.25){18}{\line(0,1){0.1}}
\put(-1.8,0.5){$d$}
\put(-1.8,-0.55){$1$}
\put(-1.8,2.2){$1'$}
\put(-1.8,3.3){$d'$}
\put(0.1,-0.045){\circle*{0.1}}
\put(0.15,-0.5){$p_0$}
\put(-3.4,2.7){$p_1$}
\put(0.05,2.6){$p_2$}
\put(3.3,2.7){$p_3$}
\put(-4.12,1.4){$t_1$}
\put(-3.13,1.4){$t_d$}
\put(2.91,1.4){$t'_1$}
\put(4.08,1.4){$t'_d$}
\end{picture}
\end{center}

All the special points are sent to the real axis by the projection $\pi$;
from left to right, there are $d$ tangency points $t_1,\dots,t_d$, followed
by $d(d-1)/2$ nodes near $p_1$, the multiple point $p_0$, 
nodes near $p_2$, nodes near $p_3$, and finally $d$ tangency points
$t'_1,\dots,t'_d$. The base point is chosen on the real axis, immediately to
the right of $\pi(p_1)$; the $d$ conics intersect the reference fiber of
$\pi$ in $2d$ points (all along the real axis in the fiber), labelled 
$1,\dots,d,1',\dots,d'$.

The system of generating loops that we use to define the braid factorization
is given by paths joining the base point to the projections of the various
tangencies and nodes as shown in the following diagram representing the
base of the fibration $\pi$~:

\begin{center}
\setlength{\unitlength}{0.25in}
\begin{picture}(8,2)(0,-1)
\multiput(0,0)(1,0){9}{\circle*{0.1}}
\put(1.8,0){\circle*{0.1}}
\put(2.2,0){\circle*{0.1}}
\put(2,0){\circle{0.6}}
\put(4.8,0){\circle*{0.1}}
\put(5.2,0){\circle*{0.1}}
\put(5,0){\circle{0.6}}
\put(5.8,0){\circle*{0.1}}
\put(6.2,0){\circle*{0.1}}
\put(6,0){\circle{0.6}}
\multiput(0.2,0)(0.1,0){7}{\circle*{0.03}}
\multiput(7.2,0)(0.1,0){7}{\circle*{0.03}}
\put(0,-0.15){\makebox(0,0)[ct]{\small $t_1$}}
\put(1,-0.15){\makebox(0,0)[ct]{\small $t_d$}}
\put(2,-0.35){\makebox(0,0)[ct]{\small $p_1$}}
\put(4,-0.15){\makebox(0,0)[ct]{\small $p_0$}}
\put(5,-0.35){\makebox(0,0)[ct]{\small $p_2$}}
\put(6,-0.35){\makebox(0,0)[ct]{\small $p_3$}}
\put(7,-0.15){\makebox(0,0)[ct]{\small $t'_1$}}
\put(8,-0.15){\makebox(0,0)[ct]{\small $t'_d$}}
\put(2.3,0){\line(1,0){1.7}}
\qbezier[80](3,0)(4,0.5)(4.71,0.15)
\qbezier[120](3,0)(4.5,1)(5.71,0.15)
\qbezier[80](3,0)(3.5,-1)(4.5,-1)
\put(4.5,-1){\line(1,0){3}}
\qbezier[80](7.5,-1)(8.5,-1)(8.5,0)
\qbezier[30](8.5,0)(8.5,0.25)(8.25,0.25)
\qbezier[30](8.25,0.25)(8,0.25)(8,0)
\qbezier[50](8.5,0)(8.5,0.5)(7.75,0.5)
\qbezier[50](7.75,0.5)(7,0.5)(7,0)
\qbezier[80](3,0)(2.5,0.5)(1.5,0.5)
\qbezier[40](1.5,0.5)(1,0.5)(1,0)
\qbezier[70](1.5,0.5)(0.5,0.5)(0,0)
\end{picture}
\end{center}

We order the various generating loops for $\pi_1(\C-\mathrm{crit})$
counterclockwise around the base point, starting with the first of the
arcs joining the base point to the projection of a node near $p_1$, and
ending with the arc joining the base point to $\pi(t_d)$.

The contribution of each node or tangency point to braid monodromy
can be calculated using a two-step process: first, one computes the {\it
local} braid monodromy, i.e.\ the monodromy action on a fiber of $\pi$ very
close to the critical point; this is the power of a half-twist exchanging
two immediately adjacent intersection points of the considered fiber of $\pi$
with the configuration of $d$ conics. Next, the local configuration is brought
back to the chosen fixed reference fiber of $\pi$ along a prescribed arc,
and the desired braid is obtained as the image of the local monodromy
under this ``parallel transport'' operation; for the purpose of the 
calculation, it is often
efficient to perform a suitable homotopy in order to break down the given 
arc into a succession of half-circles centered on other critical values
along the real axis, since parallel transport along such a half-circle can
be explicitly described as a square root of the local monodromy.

The monodromy around the multiple intersection point $p_0$ is easily seen to
be a full twist of a disc containing the $d$ intersection points labelled 
$1,\dots,d$ in the reference fiber of $\pi$; we use the notation
$\Delta_d^2$ for this element of $B_{2d}$, which is actually the image of
the central element $\Delta_d^2\in B_d$ under the natural embedding
$B_d\hookrightarrow B_{2d}$.

In the case of the nodes near $p_1$, the braid monodromy can be computed
directly from the local picture; for a generic choice of the conics, the
intersection points labelled $1',\dots,d'$ behave as in the case of $d$
lines in general position, and their $d(d-1)/2$ intersections give rise
to the braid monodromy factorization $$L'_d=\prod_{i=1}^{d-1}\prod_{j=i+1}^d
Z_{i'j'}^2$$ or any Hurwitz equivalent expression (compare with equation
(\ref{eq:d_lines})).

For the nodes near $p_2$, the local monodromy is the same as in the case
of $p_1$, except that the ordering of the points $1,\dots,d$ is reversed
compared to the reference fiber of $\pi$ (these points are not affected
by the local monodromy anyway). Since parallel transport along 
a half-circle around $p_0$ precisely amounts to a half-rotation of a disc
containing the points labelled $1,\dots,d$, the contribution to braid
monodromy remains given by the same expression $L'_d$ as above.
Near $\pi(p_3)$, the local configuration is the same as for $p_2$ 
up to reversing the ordering of the points $1',\dots,d'$ inside the fibers
of $\pi$; this discrepancy is taken care of by parallel transport along a
half-circle centered at $p_2$, and so the contribution to braid monodromy is
again $L'_d$.

In the case of the tangency point $t_d$, the intersection of the $d$ conics
with the fiber of $\pi$ above a point immediately to the right of $\pi(t_d)$
consists of $2d$ points in the order $1,\dots,d,d',\dots,1'$ on the real
axis, and the local monodromy is a half-twist exchanging the consecutive
points $d,d'$. Parallel transport along a clockwise half-circle around
$\pi(p_1)$ induces a half-rotation of the disc containing $d',\dots,1'$ in the
clockwise direction, and therefore transforms this half-twist into
$Z_{dd'}$. 

More generally, in a fiber immediately to the right of $t_i$, the 
local picture consists of $2i$ points $1,\dots,i,i',\dots,1'$ on the
real axis, while the points $d,\dots,(i+1),(i+1)',\dots,d'$ have moved to 
the pure imaginary axis, and the local monodromy around $t_i$ 
is a half-twist exchanging the consecutive points $i$ and $i'$. Parallel 
transport along a clockwise half-circle around $\pi(t_j)$ for each $j>i$ 
rotates the two points $j$ and $j'$ clockwise by $\frac{\pi}{2}$, which 
eventually yields the following half-twist in a fiber immediately to the 
right of $\pi(t_d)$:
\begin{center}
\setlength{\unitlength}{0.25in}
\begin{picture}(11.5,2)(-1,-1)
\put(0,0){\circle*{0.1}}
\put(2,0){\circle*{0.1}}
\put(4,0){\circle*{0.1}}
\put(5,0){\circle*{0.1}}
\put(7,0){\circle*{0.1}}
\put(9,0){\circle*{0.1}}
\multiput(0.5,0)(0.2,0){6}{\circle*{0.02}}
\multiput(2.5,0)(0.2,0){6}{\circle*{0.02}}
\multiput(5.5,0)(0.2,0){6}{\circle*{0.02}}
\multiput(7.5,0)(0.2,0){6}{\circle*{0.02}}
\put(0,-0.2){\makebox(0,0)[ct]{$1$}}
\put(2,-0.2){\makebox(0,0)[ct]{$i$}}
\put(4,-0.2){\makebox(0,0)[ct]{$d$}}
\put(5,-0.2){\makebox(0,0)[ct]{$d'$}}
\put(7.2,-0.2){\makebox(0,0)[ct]{$i'$}}
\put(9,-0.2){\makebox(0,0)[ct]{$1'$}}
\qbezier[100](2,0)(2,0.6)(3.2,0.6)
\qbezier[110](6,-1)(4.5,-1)(4.5,-0.2)
\qbezier[100](4.5,-0.2)(4.5,0.6)(3.2,0.6)
\qbezier[100](6,-1)(7,-1)(7,0)
\end{picture}
\end{center}
Finally, going around $\pi(p_1)$ we need to perform a clockwise half-rotation
of a disc containing $d',\dots,1'$, which yields the half-twist $Z_{ii'}$
in the reference fiber of $\pi$. Therefore, the contribution of
$t_i$ to the braid factorization is $Z_{ii'}$.

The tangencies $t'_1,\dots,t'_d$ are handled in the exactly the same manner;
the calculations are slightly more tedious because of the more complicated
choices of arcs joining $\pi(t'_i)$ to the base point, but one easily checks
that the braid monodromy around $t'_i$ is again the half-twist $Z_{ii'}$.

Putting the various contributions together in the correct order, we obtain
that the braid monodromy for the chosen configuration of conics can be
expressed by the factorization
\begin{equation}\label{eq:d_conics}\Delta_{2d}^2=
L'_d\cdot \prod_{i=1}^d Z_{ii'}\cdot \Delta_d^2\cdot (L'_d)^2\cdot
\prod_{i=1}^d Z_{ii'}.\end{equation}

As explained at the beginning of this section, in order to get the braid 
factorization for $V'_2(D_k)$ we need to replace in (\ref{eq:d_conics})
the factor $\Delta_d^2$, corresponding to the local monodromy at the
intersection point $p_0$, with the braid factorization corresponding to 
$D_k$, embedded into $B_{2d}$ in the natural way by considering a disc 
containing the points $1,\dots,d$ (see also the remark below). This immediately 
yields the formula (\ref{eq:folding0}).

The equivalent expression (\ref{eq:folding}) is obtained from
(\ref{eq:folding0}) by a sequence of Hurwitz moves, or equivalently,
by a change in the choice of generators for $\pi_1(\C-\mathrm{crit})$. 
Indeed, moving the
factors in the first $L'_d$ to the right across $\prod Z_{ii'}$ and $F_k$
affects these latter factors by a conjugation by the inverse of the product
of all the factors in $L'_d$, i.e.\ by a clockwise full twist of the disc
containing the points $1',\dots,d'$. As a result, $Z_{ii'}$ is transformed
into $\hat{Z}_{ii'}$, while the factors in $F_k$ commute with those in
$L'_d$ and remain unaffected. This completes the proof of Proposition
\ref{prop:folding}.
\medskip

\noindent {\bf Remark.} As observed in \S 1.1, the braid factorization $F_k$
is only defined up to certain algebraic operations, among which global
conjugation by an element of $B_d$. 
At first glance, the expressions obtained from (\ref{eq:folding0}) and 
(\ref{eq:folding}) by replacing $F_k$ with its conjugate $(F_k)_Q$ by some 
braid $Q\in B_d$ appear to be inequivalent to the original unconjugated
ones. Nonetheless, as suggested by the geometric intuition, all possible
choices yield equivalent results for the braid factorization of $V'_2(D_k)$.
More precisely, defining $X_r=Z_{r,r+1}$ and $X'_r=Z_{r',(r+1)'}$ for any
$1\!\le\!r\!\le\!d\!-\!1$, we claim that replacing $F_k$ by $(F_k)_{X_r}$ in the 
r.h.s. of (\ref{eq:folding}) yields an expression which is Hurwitz and 
conjugation equivalent to the original one. This is proved by observing
that the conjugated expressions $(L'_d)_{X'_r}$, 
$(\prod \hat{Z}_{ii'})_{X_r X'_r}$ and $(\prod Z_{ii'})_{X_r X'_r}$ are 
Hurwitz equivalent to the unconjugated ones (checking these identities is
an easy task left to the reader), so that a global conjugation by $X_r X'_r$ 
and a sequence of Hurwitz moves can compensate for the conjugation of $F_k$.

\subsection{The $V_2$ branch curve}\label{ss:v2}
We now compute the braid factorization corresponding to the branch curve
$C_2$ of the quadratic map $V'_2$ (or more generally of any generic
quadratic holomorphic map from $\CP^2$ to itself). Elementary calculations 
show that $C_2$ is a curve of degree $6$ with nine cusps, no nodal points,
and tangent to the fibers of $\pi$ in three points.

The braid factorizations for branch curves of generic polynomial maps from 
$\CP^2$ to itself in any degree have been computed by Moishezon 
\cite{MoiVeronese} (see also \cite{MTVeronese}), using a very technical and 
intricate argument. For the sake of completeness, we provide a direct 
calculation in the degree $2$ case.

\begin{prop}\label{prop:v2}
The braid factorization for the branch curve of $V'_2$ is
given by the formula
\begin{equation}\label{eq:v_2}
\Delta_6^2=\bigl(Z_{13}^3Z_{14}^3Z_{12;(34)}Z_{23}^3\bigr)\cdot
\bigl(Z_{15}^3Z_{16}^3Z_{12;(56)}Z_{25}^3\bigr)\cdot
\bigl(Z_{35}^3Z_{36}^3Z_{34;(56)}Z_{45}^3\bigr),
\end{equation}
where $Z_{ab;(cd)}=(Z^2_{bc}Z^2_{bd})Z_{ab}(Z^2_{bc}Z^2_{bd})^{-1}$ is
a half-twist interchanging $a$ and $b$ along a path that goes around
the points labelled $c$ and $d$.
\end{prop}

Proposition \ref{prop:v2} is proved by studying the effect of a generic
small deformation of the degenerate map $V_2^0:(x:y:z)\mapsto (x^2:y^2:z^2)$
on its branch curve. The ramification locus of $V_2^0$
in the source $\CP^2$ consists of three lines, which map two-to-one to three
lines in the target $\CP^2$~: the branching divisor of $V_2^0$ therefore
consists of three lines with multiplicity $2$ (this behavior is extremely
non-generic). The perturbation of $V_2^0$ into the generic map $V'_2$
in particular affects the local behavior of the branch curve
near the three points where the lines in the branch curve of $V_2^0$
intersect. It also affects the branch curve in a more global manner, since
the multiplicity 2 lines making up the branch divisor of $V_2^0$ are
deformed into a configuration without multiplicities; roughly speaking, away
from the intersection points each line of multiplicity 2 is separated into 
two distinct lines lying close to each other (even though one must keep in 
mind that the curve $C_2$ is irreducible).

In order to avoid the pole of the projection $\pi$, we compose the map
$V_2^0$ with the linear transformation $(x\!:\!y\!:\!z)\mapsto 
(x+\eta z\!:\!y+\eta z\!:\!z)$, for $\eta>0$ small. The resulting branch
divisor still consists of three multiplicity 2 lines, intersecting
the real slice $\R^2\subset\C^2$ in the following manner:

\begin{center}
\setlength{\unitlength}{0.25in}
\begin{picture}(5,3)(-2.5,-0.5)
\put(-1.3,-0.5){\line(1,2){1.5}}
\put(1.3,-0.5){\line(-1,2){1.5}}
\put(-1.5,0){\line(1,0){3}}
\put(-1.4,0.1){\makebox(0,0)[rb]{\small $q_1=(1\!:\!0\!:\!0)$}}
\put(1.4,0.1){\makebox(0,0)[lb]{\small $(0\!:\!1\!:\!0)=q_3$}}
\put(0.2,2.1){\makebox(0,0)[lc]{\small $(\eta\!:\!\eta\!:\!1)=q_2$}}
\put(0.8,1.2){\makebox(0,0)[lc]{\small 56}}
\put(-0.8,1.2){\makebox(0,0)[rc]{\small 12}}
\put(0,-0.1){\makebox(0,0)[ct]{\small 34}}
\end{picture}
\end{center}

On this diagram, the fibers of $\pi$ correspond to vertical lines.
We choose the reference fiber of $\pi$ far to the left on the real axis;
after a generic perturbation, each of the three lines gives rise to two 
intersection points between $C_2$ and the reference fiber of $\pi$, 
for a total of $6$ intersection points, all lying close to the real axis
in the fiber. We label these points from $1$ to $6$ in the natural order
along the real
axis, namely we label $1$ and $2$ the two intersection points corresponding 
to the line $y=0$~; we label $3$ and $4$ those corresponding to $z=0$, and 
finally $5$ and $6$ those corresponding to $x=0$. 

The braid factorization is computed by considering the three intersection
points, which obviously play very similar roles. The first intersection
point $q_1$, for which we study the braid monodromy by considering paths close
to the real axis in the base, involves the double lines $1-2$ and 
$3-4$, the first of which has the greatest slope~; computations in local 
coordinates yield a word in the braid group $B_4$, which needs to be
embedded into $B_6$ simply by considering a disc containing the points
$1,2,3,4$ and centered on the real axis (the ``parallel transport'' 
operation is trivial in this case).

The second intersection point $q_2$ involves $1-2$ and $5-6$, the first of 
which again has the greatest slope~; because the local picture is the same, 
the local computation yields the same word in $B_4$ as for $q_1$.
In the base of the fibration $\pi$, we choose to join $\pi(q_2)$ to the
reference fiber via a path passing above the real axis; one easily checks
that parallel transport along this path (going around $\pi(q_1)$ in the
counterclockwise direction) amounts to a counterclockwise
half-rotation of the points $1-2$ around $3-4$. As a consequence,
we must now use an embedding of $B_4$ into $B_6$ corresponding 
to a domain containing the points $1,2,5,6$ and passing {\it above} the real 
axis near the points $3,4$.
Finally, the third point $q_3$ involving $3-4$ and $5-6$ again
corresponds to the same local picture. We choose to join $\pi(q_3)$ to the
base point by a path passing above the real axis, and one easily checks that,
after parallel transport around $\pi(q_1)$ and $\pi(q_2)$, the relevant
embedding of $B_4$ into $B_6$ is simply that given by a disc containing
the points $3,4,5,6$ and centered on the real axis.

Consider any of the three intersection points $q_1,q_2,q_3$, where we want 
to compute the local contribution to braid monodromy after a small generic 
perturbation. Above such a point, the map $V_2^0$ is given in local affine 
coordinates by
$(x,y)\mapsto (x^2,y^2)$~; we choose to perturb it into the map
$$f:(x,y)\mapsto (x^2+\alpha y,y^2+\beta x),$$ where $\alpha$ and $\beta$
are small nonzero constants. The ramification curve is given by the
vanishing of the Jacobian of $f$, which is $4xy-\alpha\beta$~; the branch
curve of $f$ is therefore parametrized as
$$\left\{\Bigl(x^2+\frac{\alpha^2\beta}{4x},\frac{\alpha^2\beta^2}{16x^2}
+\beta x\Bigr),x\in\C-\{0\}\right\}.$$
We also need to specify the projection map in the local coordinates~: 
it can be assumed to be $(z_1,z_2)\mapsto z_1+\epsilon z_2$ for a small 
nonzero value of the constant $\epsilon$.

With this setup, the branch curve of $f$ presents one tangency point and
three cusps, and the corresponding factorization in $B_4$ can be expressed as
\begin{equation}\label{eq:v2_3pt}
Z_{13}^3\cdot Z_{14}^3\cdot Z_{12;(34)}\cdot Z_{23}^3,
\end{equation}
where $Z_{12;(34)}=(Z^2_{23}Z^2_{24})Z_{12}(Z^2_{23}Z^2_{24})^{-1}$
is the following half-twist~:

\begin{center}
\setlength{\unitlength}{0.25in}
\begin{picture}(5,1.7)(-1,-1)
\put(0,0){\circle*{0.1}}
\put(1,0){\circle*{0.1}}
\put(2,0){\circle*{0.1}}
\put(3,0){\circle*{0.1}}
\put(0,-0.2){\makebox(0,0)[ct]{1}}
\put(0.8,-0.2){\makebox(0,0)[ct]{2}}
\put(2,-0.2){\makebox(0,0)[ct]{3}}
\put(3,-0.2){\makebox(0,0)[ct]{4}}
\qbezier[120](0,0)(0.5,0.7)(2,0.7)
\qbezier[120](2,0.7)(3.7,0.7)(3.7,-0.2)
\qbezier[100](3.7,-0.2)(3.7,-1)(2.5,-1)
\qbezier[100](2.5,-1)(1.5,-1)(1,0)
\end{picture}
\end{center}

One can easily check that the product of the factors in (\ref{eq:v2_3pt})
is equal to $Z_{12}Z_{34}Z_{13}^2Z_{14}^2Z_{23}^2Z_{24}^2$, which amounts
to the double lines $1-2$ and $3-4$ intersecting each other while the two
lines in each double line (1 and 2 on one hand, 3 and 4 on the other hand)
twist by a half-turn around each other~: this is exactly the expected 
contribution (the presence of the half-twists is due to the fact that each 
double line is the image of a $2:1$ covering branched at the singular point).

It is worth mentioning that the expression (\ref{eq:v2_3pt}) 
is Hurwitz equivalent to its conjugates under the action of the 
half-twists $Z_{12}$ or $Z_{34}$ (or any combination of them).
This ``invariance property'' is suggested by the geometric 
intuition, since the two points of each pair $1-2$ or $3-4$ arising from the 
perturbation of a double line play interchangeable roles; in fact, the
diffeomorphisms $Z_{12}$ or $Z_{34}$ of the reference fiber of $\pi$ are
induced by suitable changes in the parameters of the perturbation
(from $\beta$ to $-\beta$ via $(e^{i\theta}\beta)_{0\le\theta\le\pi}$,
and similarly for $\alpha$, respectively). This observation explains why,
although the three embeddings $B_4\to B_6$ described above are in fact
naturally determined only up to certain conjugations, we need not worry
about the lack of canonicality of our choices.

Finally, one obtains the braid factorization for $C_2$ by putting together 
the contributions of the three intersection points $q_1,q_2,q_3$, using the
embeddings $B_4\to B_6$ described above.
The images of (\ref{eq:v2_3pt}) under these embeddings are exactly the
three expressions appearing in the r.h.s. of (\ref{eq:v_2}).
One easily checks that all the special points of $C_2$ are accounted for,
either by using the Pl\"ucker formulas to
show that $C_2$ only has 9 cusps and 3 tangency points, or by verifying 
directly that the product of the factors in the r.h.s. of (\ref{eq:v_2})
is equal to the central element $\Delta_6^2$. This completes the proof
of Proposition~\ref{prop:v2}.\endproof

Let us point out that, although (\ref{eq:v2_3pt}) looks very similar to the formula
obtained by Moishezon for the braid monodromy at what he calls a ``3-point''
\cite{MSegre}, the two geometric situations are very different~: Moishezon's
3-points correspond to a generic projection of a very degenerate algebraic
surface, with locally a covering map of degree 3, while the points we
describe here correspond to a very degenerate projection of a smooth
algebraic surface, with locally a covering map of degree 4. The fact
that two very different geometric descriptions of the curve $C_2$ yield
identical braid factorizations is one of the many remarkable properties of
quadratic maps from $\CP^2$ to itself.%
\medskip

We finish this section by briefly describing the geometric monodromy
representation $\theta_{V_2}:\pi_1(\CP^2-C_2)\to S_4$ corresponding to
the factorization (\ref{eq:v_2}). Each double line in the branch curve of 
$V_2^0$ corresponds to two disjoint transpositions in $S_4$, while the
transpositions corresponding to lines in different double lines are
adjacent. Therefore, after a suitable reordering of the four sheets
of the covering $V'_2$, one may assume that the six geometric generators
$\gamma_1,\dots,\gamma_6$ (small loops going around each of the six
points labelled $1,\dots,6$ in the reference fiber of $\pi$) are mapped to the 
transpositions $(1\ 2)$, $(3\ 4)$, $(1\ 3)$, $(2\ 4)$, $(1\ 4)$ and $(2\ 3)$
respectively. One easily checks that all the braids appearing in the
factorization (\ref{eq:v_2}) satisfy the compatibility relations
stated in the introduction (e.g., for the first factor $Z_{13}^3$, the
transpositions $\theta_{V_2}(\gamma_1)=(1\ 2)$ and 
$\theta_{V_2}(\gamma_3)=(1\ 3)$ are indeed adjacent).

\subsection{Regeneration of the mutual intersections}\label{ss:3pt}
We now describe the contribution to the braid monodromy of $D_{2k}$ of an
intersection point of $V'_2(D_k)$ with $C_2$. As observed in \S 2.2, 
the behavior of the map $f'_{2k}=V'_2\circ f_k$ above such a point 
is not generic, and a perturbation is
needed in order to obtain the generic map $f_{2k}$.
The local description of this perturbation is the following~:

\begin{lemma}
Over a neighborhood of a point where $R_k$ intersects $f_k^{-1}(R'_2)$,
up to an isotopy
of the branch curve among locally braided curves we can assume that
$f'_{2k}$ and $f_{2k}$ are given by the following models in local complex 
coordinates: $f'_{2k}(x,y)=(-x^2+y,-y^2)$, and 
$f_{2k}(x,y)=(-x^2+y,-y^2+\epsilon x)$, where $\epsilon$ is a small 
non-zero constant, $\pi$ being the projection to the first component.
\end{lemma}

\begin{proof}  
Provided that $k$ is large enough and given a point 
$p\in R_k\cap f_k^{-1}(R'_2)$, the argument in \S 3 of \cite{A2} (see also
\cite{AK},\cite{Aseoul}) implies 
that a small perturbation term, localized near $p$, can be added to 
$f'_{2k}$ in order to make it generic and achieve the required transversality 
properties near $p$~;
the other transversality properties of $f'_{2k}$ are not affected if the
perturbation is chosen small enough. Moreover, the one-parameter
construction used in \cite{A2} to prove uniqueness up to isotopy implies that 
the space of admissible perturbations is path connected (once again 
provided that $k$ is large enough). 

Local models for the various maps can be obtained as follows.
First observe that there exist local holomorphic coordinates $(z_1,z_2)$ on 
$\CP^2$ near $f_k(p)$ in which $V'_2$ can be expressed as 
$(z_1,z_2)\mapsto (z_1,-z_2^2)$. Moreover, it was shown
in \cite{A2} that there exist local approximately holomorphic
coordinates $(x,y)$ on $X$ and $(\tilde{z}_1,\tilde{z}_2)$ on $\CP^2$ 
in which $f_k$ is given by $(x,y)\mapsto (x^2,y)$.

Recall that $f_k$ satisfies properties $(2)$ and $(5)$ of Definition~7.
Therefore, provided that $V'_2$ is chosen sufficiently close to $V_2^0$ 
(which is always assumed to be the case), we know two things~: first, by
property $(2)$, the branch curve $D_k=f_k(R_k)$ intersects the ramification
curve $R'_2$ of $V'_2$ transversely~; second, by property $(5)$, the 
tangent space to $D_k$ at $f_k(p)$ does not lie in the kernel of the 
differential of $V'_2$, i.e.\ the image by $V'_2$ of the branch curve 
of $f_k$ is locally immersed. Therefore, $D_k$, given by the equation
$\tilde{z}_1=0$, is transverse at $f_k(p)$ to both axes of the coordinate 
system $(z_1,z_2)$ on $\CP^2$.

A first consequence is that $(\tilde{z}_1,z_2)$ are local approximately
holomorphic coordinates on $\CP^2$~; replacing the coordinate $y$ on $X$ by
$\tilde{y}=f_k^*(z_2)$, we obtain that the expression of $f_k$ in the
local coordinates $(x,\tilde{y})$ and $(\tilde{z}_1,z_2)$ remains
$(x,y)\mapsto (x^2,y)$.

Another consequence is that the coefficients of $\tilde{z}_1$ and $z_2$
in the expression of $z_1$ as a function of $\tilde{z}_1$ and $z_2$ are
both non-zero. Therefore, near the origin
we can write $z_1=\tilde{z}_1\phi(\tilde{z}_1,z_2)+z_2\psi(z_2)+O(k^{-1/2})$, 
where $\phi$ and $\psi$ are non-vanishing holomorphic functions
and the last part corresponds to the antiholomorphic terms.

Working with coordinates $(z_1,z_2)$ on $\CP^2$, the expression of $f_k$
becomes $(x,y)\mapsto (x^2\phi(x^2,y)+y\psi(y)+O(k^{-1/2}),y)$. Performing
the coordinate change $(x,y)\mapsto (ix\phi(x^2,y)^{1/2},y)$ on $X$, we
can reduce the model for $f_k$ to the simpler expression $(x,y)\mapsto
(-x^2+y\psi(y)+O(k^{-1/2}),y)$. Decomposing $\psi$ into even and odd
degree parts, we can write
$$f'_{2k}(x,y)=(-x^2+y\psi_0(y^2)+y^2\psi_1(y^2)+O(k^{-1/2}),-y^2).$$
Composing with the coordinate change $(u,v)\mapsto 
(u+v\psi_1(-v),v\psi_0(-v)^2)$ on $\CP^2$, we reduce to 
$f'_{2k}(x,y)=(-x^2+y\psi_0(y^2)+O(k^{-1/2}),-y^2\psi_0(y^2)^2)$.
Finally, the coordinate change $(x,y)\mapsto (x,y\psi_0(y^2))$ on $X$ yields
the expression $f'_{2k}(x,y)=(-x^2+y+O(k^{-1/2}),-y^2)$. This expression
differs from the desired one only by antiholomorphic terms, which are
bounded by $O(k^{-1/2})$ and therefore can be discarded without affecting 
the local braid monodromy computations.

We know that for large enough $k$ the space of admissible asymptotically 
holomorphic local perturbations of $f'_{2k}$ near $p$ (i.e.\ perturbations
satisfying the required uniform transversality properties) is path
connected. Therefore, we are free to choose the perturbation which suits 
best our purposes~; fixing a constant $\epsilon\neq 0$, we set $f_{2k}$
to be of the form $(x,y)\mapsto (-x^2+y,-y^2+\epsilon x)$. One easily checks
that, provided that the chosen value of $\epsilon$ is bounded from below
independently of $k$, this map locally satisfies all the required
properties.

Concretely, this perturbation of $f'_{2k}$ can be performed in the same
manner as in \cite{A2}, by considering the very localized asymptotically
holomorphic sections $s^\mathrm{ref}_{k,p}$ of $L^{\otimes k}$ with
exponential decay away from $p$ first introduced by Donaldson in \cite{D1}.
It is easy to check that, by adding to one of the sections of
$L^{\otimes 2k}$ defining the covering map $f'_{2k}$ a small multiple of 
$x\cdot s^\mathrm{ref}_{2k,p}$, where $x$ is the first coordinate function on
$X$ near $p$, the map $f'_{2k}$ itself is affected by a perturbation which
coincides at the first order with the desired one. In view of the local
models, this is sufficient to ensure that the branch curve agrees with the
prescribed one up to isotopy among locally braided curves, and hence to
ensure that the braid monodromy is as desired. In fact, replacing the
coefficient in front of $s^\mathrm{ref}_{2k,p}$ by a suitable polynomial of 
higher degree in the coordinates, we can even make the perturbation of 
$f'_{2k}$ coincide with the desired one up to arbitrarily high order.

We finally consider the projection $\pi$ used to define braid monodromy.
Recall that the various hypotheses made on $V'_2$ and $f_k$ ensure that
the branch curve of $V'_2$ remains locally transverse to the fibers of
$\pi$. Furthermore, over a neighborhood of the considered point, the 
tangent space to the branch curve of $f'_{2k}$ in $\CP^2$ remains very 
close to the direction determined by the branch curve of $V'_2$ 
(in our local model, the first coordinate axis)~; an easy calculation 
shows that the same property remains true for $f_{2k}$ (see also below). 
It follows that the local braid monodromy does not depend at all on choice 
of the projection $\pi$ as long as its fibers are locally transverse to the 
first coordinate axis. Therefore, performing if necessary an isotopy among
locally braided curves by means of a suitable rotation, we can safely
assume $\pi$ to be the projection to the first coordinate axis.
\end{proof}

By Lemma 2 we know that the local braid monodromy of $f_{2k}$
can be computed using for $f_{2k}$ the local model
$$(x,y)\mapsto (-x^2+y,-y^2+\epsilon x)$$ where $\epsilon$ is a 
small non-zero constant. The Jacobian of this map is $4xy-\epsilon$, and
its branch curve can be parametrized as
$$\left\{\Bigl(-x^2+\frac{\epsilon}{4x},-\frac{\epsilon^2}{16x^2}+\epsilon x
\Bigr),x\in\C-\{0\}\right\}.$$

The signs have been chosen in such a way that, taking $\epsilon$ along the
positive real axis and taking the base point at a large negative real value
of the first coordinate, the intersection of the branch curve with the 
reference fiber of $\pi$ consists of three points aligned along the real 
axis, the left-most one corresponding to the branch curve of $f_k$ while 
the two others correspond to the branch curve of $V'_2$.

Projecting to the first component (or choosing any other generic projection),
the only remarkable features of the branch curve near the origin are three 
cusps, and the corresponding braid factorization can be expressed as
\begin{equation}\label{eq:regener1}
Z_{12}^3\cdot Z_{13}^3\cdot Z_{12;(3)}^3,
\end{equation}
where the point labelled $1$ corresponds to the branch curve of $f_k$ while
the points labelled $2$ and $3$ correspond to the branch curve of $V'_2$,
and where $Z_{12;(3)}=Z_{23}^2Z_{12}Z_{23}^{-2}$ is a half-twist exchanging
$1$ and $2$ along a path that goes around $3$~:

\begin{center}
\setlength{\unitlength}{0.25in}
\begin{picture}(4,1.7)(-1,-1)
\put(0,0){\circle*{0.1}}
\put(1,0){\circle*{0.1}}
\put(2,0){\circle*{0.1}}
\put(0,-0.2){\makebox(0,0)[ct]{1}}
\put(0.8,-0.2){\makebox(0,0)[ct]{2}}
\put(2,-0.2){\makebox(0,0)[ct]{3}}
\qbezier[100](0,0)(0.5,0.7)(1.5,0.7)
\qbezier[100](1.5,0.7)(2.7,0.7)(2.7,-0.2)
\qbezier[80](2.7,-0.2)(2.7,-1)(2,-1)
\qbezier[80](2,-1)(1.5,-1)(1,0)
\end{picture}
\end{center}

A short calculation in $B_3$ shows that the product of the factors in
(\ref{eq:regener1}) is equal to $Z_{23}(Z_{12}^2Z_{13}^2)^2$, which
amounts to the line labelled $1$ twisting twice around $2$ and $3$ while
these two lines undergo a half-twist. This is consistent with the geometric
intuition, since the branch curve of $f_k$, folded onto itself by $V'_2$, 
hits the branch curve of $V'_2$ in a manner that can be represented by
the following picture:

\begin{center}
\setlength{\unitlength}{0.25in}
\begin{picture}(6,2.5)(-3,-2)
\put(-3,0.15){\line(1,0){2.6}}
\put(-3,-0.15){\line(1,0){2.6}}
\put(0.4,0.15){\line(1,0){2.6}}
\put(0.4,-0.15){\line(1,0){2.6}}
\qbezier[40](-0.4,0.15)(0,0.15)(0,0)
\qbezier[40](-0.4,-0.15)(0,-0.15)(0,0)
\qbezier[40](0.4,0.15)(0,0.15)(0,0)
\qbezier[40](0.4,-0.15)(0,-0.15)(0,0)
\qbezier[200](-2,-2)(-1,-0.15)(0,-0.15)
\qbezier[200](2,-2)(1,-0.15)(0,-0.15)
\qbezier[20](-3,-1.1)(-2.35,-1.1)(-1.7,-1.1)
\qbezier[20](-3,-0.8)(-2.35,-0.8)(-1.7,-0.8)
\qbezier[20](0.4,-1.1)(-0.25,-1.1)(-0.9,-1.1)
\qbezier[20](0.4,-0.8)(-0.25,-0.8)(-0.9,-0.8)
\qbezier[8](-1.7,-1.1)(-1.3,-1.1)(-1.3,-0.95)
\qbezier[8](-1.7,-0.8)(-1.3,-0.8)(-1.3,-0.95)
\qbezier[8](-0.9,-1.1)(-1.3,-1.1)(-1.3,-0.95)
\qbezier[8](-0.9,-0.8)(-1.3,-0.8)(-1.3,-0.95)
\qbezier[20](3,-1)(2.35,-1)(1.7,-1)
\qbezier[20](3,-0.7)(2.35,-0.7)(1.7,-0.7)
\qbezier[20](-0.4,-1)(0.25,-1)(0.9,-1)
\qbezier[20](-0.4,-0.7)(0.25,-0.7)(0.9,-0.7)
\qbezier[8](1.7,-1)(1.3,-1)(1.3,-0.85)
\qbezier[8](1.7,-0.7)(1.3,-0.7)(1.3,-0.85)
\qbezier[8](0.9,-1.0)(1.3,-1.0)(1.3,-0.85)
\qbezier[8](0.9,-0.7)(1.3,-0.7)(1.3,-0.85)
\put(-2.5,-2){1}
\put(-3.4,0.1){3}
\put(-3.4,-0.4){2}
\end{picture}
\end{center}
The line labelled $1$ intersects $2$ and $3$ with multiplicity 2 because
the image of $D_k$ by $V'_2$ is necessarily tangent to the branch curve
of $V'_2$ wherever they intersect~; the lines $2$ and $3$ twist around
each other by a half-turn because they arise as the two sheets of a
2:1 covering branched at the origin (they correspond
to the two preimages by $f_k$ of each point where $V'_2$ is ramified).

It is also worth observing that the expression (\ref{eq:regener1}) is
easily shown to be Hurwitz equivalent to its conjugates under the action
of the group generated by the half-twist $Z_{23}$ exchanging the two points
labelled $2$ and $3$, in agreement with the geometric intuition suggesting
that their roles are interchangeable.

In order to understand how the braid monodromy given in 
(\ref{eq:regener1}) fits in the global picture, 
we now need to explain the labelling of the various components making up
$D_{2k}$ and the corresponding geometric monodromy representation.
\medskip

{\bf Notations.}
As described above the branch curve $D_{2k}$ is obtained by deforming
the union of $V'_2(D_k)$ and $n$ copies of $C_2$. Its degree is therefore
$\bar{d}=2d+6n$. For braid group calculations in $B_{2d+6n}$, the $2d+6n$
intersection points of $D_{2k}$ with the reference fiber of $\pi$ will be
labelled as follows: we assign labels $1,\dots,d$ and $1',\dots,d'$ to the
$2d$ intersection points which correspond to
$V'_2(D_k)$ (in the same manner as in \S \ref{ss:folding}), 
and $i_\alpha,i'_\alpha,i_\beta,i'_\beta,i_\gamma,i'_\gamma$ 
for $1\le i\le n$ to the $6n$ intersection points corresponding to the
$n$ copies of $C_2$. More precisely, recall that the branch curve of 
$V'_2$ is obtained as a perturbation of the branch curve of $V_2^0$, which
consists of three double lines~: therefore the $n$ copies of $C_2$ can be 
thought of as three groups of $2n$ lines. These
three groups correspond to the three subscripts $\alpha$, $\beta$ and
$\gamma$~; for each value of $i$ the two labels $i_\alpha$ and $i'_\alpha$
correspond to the perturbation of a double line in the $i$-th copy of
the branch curve of $V_2^0$.

We will choose the reference fiber of $\pi$ and the configuration of the
branch curve in such a way that the $2d+6n$ intersection points of the
reference fiber with $D_{2k}$ all lie on the real axis, in the
order $1,\dots,d,1',\dots,d',$ $1_\alpha,1'_\alpha,2_\alpha,2'_\alpha,\dots,$
$n_\alpha,n'_\alpha,$ $1_\beta,1'_\beta,\dots,n_\beta,n'_\beta,$ $1_\gamma,
1'_\gamma,\dots,n_\gamma,n'_\gamma$; in fact, we will actually choose a
reference fiber yielding a slightly different configuration of intersection
points, and then conjugate the obtained monodromy by a suitable braid.
In any case, when using $Z_{ij}$ notations it will be understood
that the $2d+6n$ intersection points of $D_{2k}$ with the reference fiber
of the projection $\pi$ are to be placed on the real axis in the above-given
order.

In order to describe the geometric monodromy representation morphism
$\theta_{2k}:\pi_1(\CP^2-D_{2k})\to S_{4n}$, we first need to choose
a set of geometric generators of $\pi_1(\CP^2-D_{2k})$. We choose the base
point for $\pi_1(\CP^2-D_{2k})$ to lie in the reference
fiber of $\pi$, far above the real axis which contains the $2d+6n$
intersection points with $D_{2k}$, and we use a system of $2d+6n$
generating loops, each joining the base point to one of the intersection
points along a straight line, circling once around the intersection point,
and going back to the base point along the same straight line.

The $4n$ sheets of the
covering $f_{2k}$ can be thought of as four groups of $n$ sheets, which
we will label as $i_a,i_b,i_c,i_d$ for $1\le i\le n$. Consider a
situation similar to that of \S \ref{ss:folding}, where most of the branch 
curve of $f_k$ is concentrated into a small ball
far away from the branch curve of $V'_2$~: this results in a picture
where the parts of the branch curve corresponding to $V'_2(D_k)$ connect 
to each other the $n$ 
sheets of a single group ($1_a,\dots,n_a$ for example), while the copies 
of $C_2$ connect the various groups of $n$ sheets to each other. 
In particular, the transpositions
in $S_{4n}$ corresponding to the geometric generators around
$1,\dots,d,1',\dots,d'$ are directly given by the geometric monodromy 
representation $\theta_k$ associated to $D_k$~: for any $1\le r\le d$, 
if $\theta_k$ maps the $r$-th geometric generator to the 
transposition $(ij)$ in $S_n$ then, calling $\gamma_r$ and $\gamma_{r'}$ the 
geometric generators in $\pi_1(\CP^2-D_{2k})$ corresponding to $r$ and $r'$, 
one gets $\theta_{2k}(\gamma_r)=\theta_{2k}(\gamma_{r'})=(i_aj_a)$.
Finally, each of the $n$ copies of $C_2$ connects four sheets to each other,
one in each group of $n$, in the same manner as for $V'_2$ itself~: 
therefore $\theta_{2k}$ maps the geometric
generators around $i_\alpha$, $i'_\alpha$, $i_\beta$, $i'_\beta$, $i_\gamma$ 
and $i'_\gamma$ to 
$(i_ai_b)$, $(i_ci_d)$, $(i_ai_c)$, $(i_bi_d)$, $(i_ai_d)$ and $(i_bi_c)$
respectively, for all $1\le i\le n$.

A suitable choice of geometric configuration and reference fiber of $\pi$
yields a situation in which $\theta_{2k}$ is as described above. Our choice
of configuration will be made explicit in \S \ref{ss:assembling}. 
A different set of geometric choices would lead to a different description of
the braid monodromy and of $\theta_{2k}$, but the final answer always
remains the same up to Hurwitz and conjugation equivalence.\medskip

With this understood,
we now describe the contribution to the braid monodromy of a point where
a piece of $V'_2(D_k)$, say e.g. the portion of conic labelled $r'$ for some
$1\le r\le d$, hits one of the three groups of $2n$ lines making up the
$n$ copies of $C_2$, say e.g. the lines labelled $1_\alpha,1'_\alpha,\dots,
n_\alpha,n'_\alpha$. 

If one just considers the composed map $V'_2\circ f_k$, the $n$ copies of
the branch curve $C_2$ of $V'_2$ all lie in the same position, and
the curve $V'_2(D_k)$ hits them tangently (and therefore with local
intersection multiplicity $+2$). To obtain the generic map
$f_{2k}$ we add a small perturbation, which affects the situation by 
moving the $n$ copies of $C_2$ apart from each other and also by modifying
the intersection of $R_k$ with $f_k^{-1}(R'_2)$ in the manner explained
above. More precisely, $R'_2$ admits $2n-2$ local lifts to $X$ which do not
locally intersect the branch curve of $f_k$ (because they lie in different
sheets of the covering) and thus do not require any
special treatment, while the two other sheets of $f_k$ give rise to
``lifts'' of $R'_2$ intersecting the branch curve of $f_k$ and each other.
Therefore, when computing the braid factorization of $D_{2k}$, we can
locally consider the $n$ copies of $C_2$ as consisting of $2n-2$ parallel
lines, each intersected twice by $V'_2(D_k)$ (giving rise to two nodes),
and two ``lines'' parallel to the others which are hit by $V'_2(D_k)$ in
the manner previously explained.

The geometric monodromy
representation $\theta_{2k}$ maps the geometric generator around $r'$ to
a transposition of the form $(p_aq_a)$, for some $1\le p,q\le n$. 
The two lines hit in a non-trivial manner are those
labelled $p_\alpha$ and $q_\alpha$, which under the
map $\theta_{2k}$ correspond respectively to the transpositions 
$(p_ap_b)$ and $(q_aq_b)$ in $S_{4n}$. The other $2n-2$ 
lines ($i_\alpha$ for
$i\not\in\{p,q\}$ and $i'_\alpha$ for all $i$) lie in different sheets of the
covering and their intersections with $r'$ simply
remain as nodes in the branch curve $D_{2k}$.

Parallel transport of the local configuration along a given arc
in the base of the fibration $\pi$ reveals the important role played by
two specific paths in the reference fiber, namely the path along which the 
point labelled $r'$ approaches the group of $2n$ points
$1_\alpha,\dots,n'_\alpha$ and the path along which two of these $2n$ points
approach each other. To phrase things differently, these two paths determine
an embedded triangle with vertices at $r',p_\alpha,q_\alpha$ in the
reference fiber, which collapses as one moves from the reference fiber 
towards the intersection point. 

We assume the configuration to be such that, after parallel transport
of the local configuration into the reference fiber of $\pi$, the path
along which the point labelled $r'$ approaches the $2n$ other points is the
simplest possible one passing {\it above} the real axis, while the two
points $p_\alpha$ and $q_\alpha$ approach each other along a path
isotopic to an arc contained in the upper half-plane.
Equivalently, inside the reference fiber of $\pi$, we assume that
the embedded triangle with vertices $r'$, $p_\alpha$ and $q_\alpha$ which
collapses as one approaches the considered singular point is the simplest
possible one lying in the upper half-plane.
Whether this is truly the case or whether the formula needs to be adjusted
by a suitable global conjugation will be determined later on, when the
contributions of the various points are put together into a global braid
factorization in $B_{2d+6n}$; if the motions of $r'$, $p_\alpha$ and
$q_\alpha$ are different from the (purely arbitrary) above choice, then
the formula giving the local monodromy will need to be conjugated by a
certain element of $B_{2d+6n}$ (any braid that maps the triangle joining
$r',p_\alpha,q_\alpha$ into the correct position can be used, as they all
yield Hurwitz equivalent factorizations).

\begin{prop}\label{prop:3pt}
The braid monodromy for the intersection of the portion of conic labelled
$r'$ with the $2n$
lines $1_\alpha,1'_\alpha,\dots,n_\alpha,n'_\alpha$ is Hurwitz and
conjugation equivalent to the following factorization~:

\begin{multline}\label{eq:mutual}
\prod_{i=n}^1\left(\acute{Z}^2_{r'i'_\alpha}
\bigl[\acute{Z}^2_{r'i_\alpha}\bigr]_{i\not\in\{p,q\}}\right)
\cdot Z^3_{r'p_\alpha}\cdot Z^3_{r'q_\alpha}\cdot
Z^3_{r'p_\alpha;(q_\alpha)}\cdot\\
\prod_{i=n}^1\left(\grave{Z}^2_{r'i'_\alpha}
\bigl[\grave{Z}^2_{r'i_\alpha}\bigr]_{i\not\in\{p,q\}}\right),
\end{multline}
where $\acute{Z}_{r'\tau}$ and $\grave{Z}_{r'\tau}$ are half-twists along
the following paths:

\begin{center}
\setlength{\unitlength}{0.25in}
\begin{picture}(9,2)(0,-1)
\put(0,0.2){\makebox(0,0)[cb]{$\acute{Z}_{r'\tau}$}}
\put(2,0){\circle*{0.1}}
\put(4,0){\circle*{0.1}}
\put(5,0){\circle*{0.1}}
\put(7,0){\circle*{0.1}}
\multiput(0.5,0)(0.2,0){6}{\circle*{0.02}}
\multiput(2.5,0)(0.2,0){6}{\circle*{0.02}}
\multiput(5.5,0)(0.2,0){6}{\circle*{0.02}}
\multiput(7.5,0)(0.2,0){6}{\circle*{0.02}}
\put(2,-0.2){\makebox(0,0)[ct]{$r'$}}
\put(4,-0.2){\makebox(0,0)[ct]{$d'$}}
\put(5,0.2){\makebox(0,0)[cb]{$1_\alpha$}}
\put(7,0.32){\makebox(0,0)[cb]{$\tau$}}
\qbezier[90](2,0)(2,0.8)(3.25,0.8)
\qbezier[90](3.25,0.8)(4.5,0.8)(4.5,0)
\qbezier[90](4.5,0)(4.5,-0.8)(5.75,-0.8)
\qbezier[90](5.75,-0.8)(7,-0.8)(7,0)
\end{picture}
\quad
\begin{picture}(9,2)(0,-1)
\put(0,0.2){\makebox(0,0)[cb]{$\grave{Z}_{r'\tau}$}}
\put(2,0){\circle*{0.1}}
\put(4,0){\circle*{0.1}}
\put(6,0){\circle*{0.1}}
\put(7,0){\circle*{0.1}}
\multiput(0.5,0)(0.2,0){6}{\circle*{0.02}}
\multiput(2.5,0)(0.2,0){6}{\circle*{0.02}}
\multiput(4.5,0)(0.2,0){6}{\circle*{0.02}}
\multiput(7.5,0)(0.2,0){6}{\circle*{0.02}}
\put(2,-0.2){\makebox(0,0)[ct]{$r'$}}
\put(4,0.2){\makebox(0,0)[cb]{$\tau$}}
\put(7,-0.2){\makebox(0,0)[ct]{$1_\beta$}}
\qbezier[110](2,0)(2,1)(3.5,1)
\put(3.5,1){\line(1,0){1.5}}
\qbezier[110](5,1)(6.5,1)(6.5,0)
\qbezier[90](6.5,0)(6.5,-1)(5.25,-1)
\qbezier[90](5.25,-1)(4,-1)(4,0)
\end{picture}
\end{center}
and $Z_{r'p_\alpha;(q_\alpha)}$ is a half-twist along the path

\begin{center}
\setlength{\unitlength}{0.25in}
\begin{picture}(10,2)(0,-1)
\put(2,0){\circle*{0.1}}
\put(4,0){\circle*{0.1}}
\put(6,0){\circle*{0.1}}
\put(7,0){\circle*{0.1}}
\put(8,0){\circle*{0.1}}
\multiput(0.5,0)(0.2,0){6}{\circle*{0.02}}
\multiput(2.5,0)(0.2,0){6}{\circle*{0.02}}
\multiput(4.5,0)(0.2,0){6}{\circle*{0.02}}
\multiput(8.5,0)(0.2,0){6}{\circle*{0.02}}
\put(2,-0.2){\makebox(0,0)[ct]{$r'$}}
\put(4,-0.2){\makebox(0,0)[ct]{$p_\alpha$}}
\put(7,-0.2){\makebox(0,0)[ct]{$q_\alpha$}}
\qbezier[110](2,0)(2,1)(3.5,1)
\put(3.5,1){\line(1,0){2.5}}
\qbezier[110](6,1)(7.5,1)(7.5,0)
\qbezier[70](7.5,0)(7.5,-1)(7,-1)
\qbezier[70](7,-1)(6.5,-1)(6.5,0)
\qbezier[90](6.5,0)(6.5,0.5)(5.25,0.5)
\qbezier[90](5.25,0.5)(4,0.5)(4,0)
\end{picture}
\end{center}

In $(\ref{eq:mutual})$, the products are to be performed in the
reverse order (first $i=n$, finishing with $i=1$), and the notation
$[\dots]_{i\not\in\{p,q\}}$ means that the enclosed factor is not
present for $i=p$ or $i=q$. 
\end{prop}

\proof
We start by considering a slightly simpler setup where, instead of being in
their normal positions, the points $p_\alpha$ and $q_\alpha$ have been moved 
to the right of the $2n-2$ other points $1_\alpha,\dots,n'_\alpha$ (i.e.,
further along the positive real axis in the reference fiber of $\pi$).
More precisely, we assume that the points $p_\alpha$ and $q_\alpha$ have
been moved into these positions along arcs in the upper half-plane, so that
the point $r'$ still reaches them by passing above the real axis and the
vanishing cycle is the line segment joining $p_\alpha$ to $q_\alpha$. 
The situation is then described by the following picture in $\R^2\subset\C^2$:

\begin{center}
\setlength{\unitlength}{0.25in}
\begin{picture}(6,2.5)(-3,-2)
\put(-3,0.15){\line(1,0){2.6}}
\put(-3,-0.15){\line(1,0){2.6}}
\put(0.4,0.15){\line(1,0){2.6}}
\put(0.4,-0.15){\line(1,0){2.6}}
\qbezier[30](-0.4,0.15)(0,0.15)(0,0)
\qbezier[30](-0.4,-0.15)(0,-0.15)(0,0)
\qbezier[30](0.4,0.15)(0,0.15)(0,0)
\qbezier[30](0.4,-0.15)(0,-0.15)(0,0)
\qbezier[200](-2,-2)(-1,-0.15)(0,-0.15)
\qbezier[200](2,-2)(1,-0.15)(0,-0.15)
\put(-2.5,-1.95){$r'$}
\put(-3.7,-0.3){$p_\alpha$}
\put(-3.7,0.2){$q_\alpha$}
\put(3.2,-1.5){$1_\alpha$}
\put(3.2,-0.6){$n'_\alpha$}
\put(-3,-0.45){\line(1,0){6}}
\multiput(-3,-0.75)(0.2,0){6}{\circle*{0.02}}
\multiput(3,-0.75)(-0.2,0){6}{\circle*{0.02}}
\put(-3,-1.05){\line(1,0){6}}
\put(-3,-1.35){\line(1,0){6}}
\end{picture}
\end{center}

The reference fiber is once again placed to the left of the diagram, and
the vertical direction corresponds to the real axis in the fibers of $\pi$.

Recalling that $V'_2(D_k)$ hits $C_2$ tangently, the expected total 
contribution to the braid monodromy corresponds to $r'$ twisting twice 
around each of the lines $1_\alpha,1'_\alpha,\dots,n_\alpha,n'_\alpha$.
For the reasons explained above, a half-twist between the lines $p_\alpha$ 
and $q_\alpha$ is also to be expected.

In order to compute the braid monodromy, we observe that in the chosen
configuration the singular fibers of $\pi$ all lie along the real axis, and
choose the following system of generating paths in the base of the fibration
$\pi$: the first path connects the base point (far away on the negative
real axis) to the first intersection
of $r'$ with $n'_\alpha$ by passing {\it below} the real axis; the second
one similarly joins the base point to the first intersection of $r'$ with
$n_\alpha$ by passing below the real axis; and so on, going from right to
left, until all $2n-2$ nodes in the left half of the diagram have been
considered. The following three paths join the base point to the three
cusp singularities arising from the perturbation of the singular point in
the middle of the diagram, passing {\it above} the real axis. Finally,
the remaining $2n-2$ paths join the base point to the intersections in the
right half of the diagram, passing above the real axis, and going from left
to right (the first of these paths ends at the second intersection of $r'$
with $n'_\alpha$, the last one ends at the second intersection with
$1_\alpha$). As should always be the case, the paths are ordered
counterclockwise around the base point.

Observing that the conic labelled $r'$ behaves similarly to the graph
of the identity function in the left half of the diagram and similarly
to the graph of $-\mathrm{Id}$ in the right half, one easily obtains the
following expression for the local braid monodromy of our configuration:
\begin{equation}\label{eq:4d_twists}
\prod_{i=n}^1\left(\acute{Z}^2_{r'i'_\alpha}
\bigl[\acute{Z}^2_{r'i_\alpha}\bigr]_{i\not\in\{p,q\}}\right)
\cdot F_{r'p_\alpha q_\alpha}\cdot
\prod_{i=n}^1\left(\grave{Z}^2_{r'i'_\alpha}
\bigl[\grave{Z}^2_{r'i_\alpha}\bigr]_{i\not\in\{p,q\}}\right),
\end{equation}
where the notation $F_{r'p_\alpha q_\alpha}$ represents
an expression similar to (\ref{eq:regener1}), and $\acute{Z}_{r'i_\alpha}$
and $\grave{Z}_{r'i_\alpha}$ are the same half-twists as in the statement of
Proposition~5.

We now bring the two points $p_\alpha$ and $q_\alpha$ back to their
respective positions, moving them along paths passing {\it above} the
real axis. The half-twists $\acute{Z}_{r'\tau}$ and $\grave{Z}_{r'\tau}$
are not affected by this motion; whereas $F_{r'p_\alpha q_\alpha}$ is
changed into $Z^3_{r'p_\alpha}\cdot Z^3_{r'q_\alpha}\cdot
Z^3_{r'p_\alpha;(q_\alpha)}$. Therefore, the expression
(\ref{eq:4d_twists}) turns into (\ref{eq:mutual}).
\endproof

Observe that the conjugates of the expression (\ref{eq:mutual}) by certain
elements of $B_{2n}$ (acting on $1_\alpha,\dots,n'_\alpha$)
are Hurwitz equivalent to (\ref{eq:mutual}).
Indeed, consider the subgroup $B_{2n-2}\times B_2\subset B_{2n}$ of braids 
which globally preserve the triangle formed by $r'$, $p_\alpha$ and 
$q_\alpha$. The factor $B_2$ is generated by the half-twist 
$Z_{p_\alpha q_\alpha}$ interchanging $p_\alpha$ and $q_\alpha$, while
the factor $B_{2n-2}$ is generated by half-twists interchanging 
two of the $2n-2$ other points along a path passing below the real axis. 
Conjugating (\ref{eq:mutual}) by $Z_{p_\alpha q_\alpha}$ simply
amounts to a modification of the three central degree $3$ factors of
(\ref{eq:mutual}) by two Hurwitz moves. Similarly, conjugation by
one of the half-twists generating $B_{2n-2}$ (interchanging two 
consecutive points among the $2n-2$) is equivalent to
two Hurwitz moves, one among the $\acute{Z}_{r'\tau}^2$ factors 
and the other among the $\grave{Z}_{r'\tau}^2$ factors. This is in
agreement with the geometric intuition suggesting that, since all these
conjugations do not affect the triangle joining $r'$, $p_\alpha$ and 
$q_\alpha$, they do not modify the braid monodromy in any significant way.

However, as already pointed out above, conjugating (\ref{eq:mutual}) by an
element of $B_{2n}$ lying outside of $B_{2n-2}\times B_2$ affects
non-trivially the path along which $p_\alpha$ and $q_\alpha$ approach each
other, and therefore yields an expression which is not Hurwitz equivalent to
the original one (this can be seen directly by observing that the product of
all factors in (\ref{eq:mutual}) is modified by the conjugation).


\subsection{The assembling rule}\label{ss:assembling}
We now study how the various elements described above fit together to provide
the braid factorization for $D_{2k}$. We will start by considering, as a
toy model, a curve made up of $d$ conics and three lines, corresponding
to the following diagram (drawn for $d=2$)~:

\begin{center}
\setlength{\unitlength}{0.25in}
\begin{picture}(8,5.2)(-4,-0.5)
\qbezier(0.2,0)(-3.3,-0.5)(-3.8,1)
\qbezier(-3.8,1)(-4.3,2.5)(-0.8,3)
\qbezier(-0.8,3)(2.7,3.5)(3.2,2)
\qbezier(3.2,2)(3.7,0.5)(0.2,0)
\qbezier(0,0)(3.5,-0.5)(4,1)
\qbezier(4,1)(4.5,2.5)(1,3)
\qbezier(1,3)(-2.5,3.5)(-3,2)
\qbezier(-3,2)(-3.5,0.5)(0,0)
\qbezier(-5,1.4)(-2,3)(1,4.6)
\qbezier(5.2,1.4)(2.2,3)(-0.8,4.6)
\put(-4,3.16){\line(1,0){8.2}}
\put(-4.8,1.8){$\alpha$}
\put(4.5,2){$\gamma$}
\put(4.4,3.1){$\beta$}
\put(-1.8,0.5){$2$}
\put(1.5,0.5){$1$}
\put(-1.8,2.2){$1'$}
\put(1.5,2.2){$2'$}
\end{picture}
\end{center}

The $d$ conics play the role of $V'_2(D_k)$, while the three lines
correspond to the $C_2$ part. As usual, the vertical direction corresponds
to the real axis in the fibers of $\pi$, and the reference fiber is to the
left of the diagram~; in the reference fiber the points
are placed on the real axis in the order $1,\dots,d,1',\dots,d',
\alpha,\beta,\gamma$.
Although the space of all configurations of $d$ conics and three lines
tangent to them in $\CP^2$ is connected, thus making all possible choices
equally suitable, the choice of the configuration represented above is
motivated by its remarkable similarity to the configurations chosen in
\S \ref{ss:folding} and \S \ref{ss:v2} for $V'_2(D_k)$ and $C_2$
respectively. In particular, one easily checks that the braid monodromy
for the chosen configuration of the $d$ conics is exactly the one computed
in \S \ref{ss:folding} (equation (\ref{eq:d_conics}) and Proposition
\ref{prop:folding}).

The braid monodromy for this configuration of $d$ conics and three
lines can be computed explicitly in coordinates. However this 
calculation is tedious and not very illuminating, so we first motivate
the formula by deriving it
by a different method~: we start from a situation where the lines
are in general position with respect to the conics, and we follow on the
level of braid factorizations the deformation of such a generic
configuration into the specific desired one. In fact,
keeping track of the deformation amounts to performing a sequence of Hurwitz
moves with the aim of bringing next to each other the two factors
arising from the intersections of each line with each conic~; the
resulting braid factorization contains consecutive identical degree $2$
factors, so that merging the intersections becomes a trivial task.

Alternatively, the reader may jump ahead to the statement of
Proposition \ref{prop:assembling} and the outline of proof given afterwards
for a description of the direct monodromy calculation.

The standard braid factorization assembling formula for the union
of two tranversely intersecting curves of respective degrees $p$ and $q$ is
given by
\begin{equation}\label{eq:assembling0}
\Delta_{p+q}^2=\Delta_p^2\cdot\prod_{i=1}^p\prod_{j=p+1}^{p+q}
Z_{ij}^2\cdot \Delta_q^2,
\end{equation}
where the points are labelled $1,\dots,p$ for the first curve and
$p+1,\dots,p+q$ for the second, and $\Delta_p^2$ and $\Delta_q^2$
stand for the braid factorizations of the two components. The braid groups
$B_p$ and $B_q$ are implicitly embedded into $B_{p+q}$ by considering two
disjoint disks containing the $p$ first points and the $q$ last points 
respectively. The formula (\ref{eq:assembling0}) can be easily checked 
by applying a suitable isotopy to the
two components so that, outside of two mutually disjoint balls, they behave
like respectively $p$ and $q$ mutually transverse lines.

In our case we want the three lines to be tangent to the conics, so we need
to perform Hurwitz moves on this factorization so that the two intersections
of each line with each conic can be brought together. Our starting point,
as given by (\ref{eq:assembling0}) and (\ref{eq:folding}), is the
factorization
\begin{multline}\label{eq:assembling1}
\Delta^2=\Bigl(\prod_{i=1}^d\hat{Z}_{ii'}\cdot L_d\cdot (L'_d)^3\cdot
\prod_{i=1}^d Z_{ii'}\Bigr)\cdot\\ \prod_{i=1}^d \left(Z_{i\alpha}^2
Z_{i\beta}^2Z_{i\gamma}^2\right)\cdot\prod_{i=1}^d \left(Z_{i'\alpha}^2
Z_{i'\beta}^2Z_{i'\gamma}^2\right)\cdot
\left(Z_{\alpha\beta}^2
Z_{\alpha\gamma}^2Z_{\beta\gamma}^2\right).
\end{multline}
Moving the $Z_{ii'}$ factors to the right, one replaces the central
$Z_{i\alpha}^2Z_{i\beta}^2Z_{i\gamma}^2$ terms by
$Z_{i'\alpha}^2Z_{i'\beta}^2Z_{'i\gamma}^2$~; then, moving the
rightmost terms to the left, one obtains the new expression
$$\Delta^2=\Bigl(\prod_{i=1}^d\hat{Z}_{ii'}\cdot L_d\cdot (L'_d)^3
\Bigr)\cdot\Bigl(\prod_{i=1}^d \left(Z_{i'\alpha}^2
Z_{i'\beta}^2Z_{i'\gamma}^2\right)\Bigr)^2\cdot
\left(Z_{\alpha\beta}^2
Z_{\alpha\gamma}^2Z_{\beta\gamma}^2\right)\cdot\prod_{i=1}^d \check{Z}_{ii'}
$$
where $\check{Z}_{ii'}$ is a half-twist along the following path~:

\begin{center}
\setlength{\unitlength}{0.25in}
\begin{picture}(14,2.2)(-1,-1.2)
\put(0,0){\circle*{0.1}}
\put(2,0){\circle*{0.1}}
\put(4,0){\circle*{0.1}}
\put(5,0){\circle*{0.1}}
\put(7,0){\circle*{0.1}}
\put(9,0){\circle*{0.1}}
\put(10,0){\circle*{0.1}}
\put(11,0){\circle*{0.1}}
\put(12,0){\circle*{0.1}}
\multiput(0.5,0)(0.2,0){6}{\circle*{0.02}}
\multiput(2.5,0)(0.2,0){6}{\circle*{0.02}}
\multiput(5.5,0)(0.2,0){6}{\circle*{0.02}}
\multiput(7.5,0)(0.2,0){6}{\circle*{0.02}}
\put(0,-0.2){\makebox(0,0)[ct]{1}}
\put(2,-0.2){\makebox(0,0)[ct]{$i$}}
\put(4,-0.2){\makebox(0,0)[ct]{$d$}}
\put(5.2,-0.2){\makebox(0,0)[ct]{1'}}
\put(7,-0.2){\makebox(0,0)[ct]{$i'$}}
\put(9,-0.2){\makebox(0,0)[ct]{$d'$}}
\put(10,-0.4){\makebox(0,0)[ct]{$\alpha$}}
\put(11,-0.3){\makebox(0,0)[ct]{$\beta$}}
\put(12,-0.4){\makebox(0,0)[ct]{$\gamma$}}
\qbezier(2,0)(2,0.8)(3,0.8)
\qbezier(3,0.8)(4.3,0.8)(4.3,0)
\qbezier(4.3,0)(4.3,-1.2)(5.3,-1.2)
\put(5.3,-1.2){\line(1,0){6}}
\qbezier(11.3,-1.2)(12.5,-1.2)(12.5,0)
\qbezier(12.5,0)(12.5,0.8)(11,0.8)
\qbezier(11,0.8)(9.5,0.8)(9.5,0)
\qbezier(9.5,0)(9.5,-0.9)(8.5,-0.9)
\put(5.5,-0.9){\line(1,0){3}}
\qbezier(5.5,-0.9)(4.6,-0.9)(4.6,0)
\qbezier(4.6,0)(4.6,0.8)(5.8,0.8)
\qbezier(5.8,0.8)(7,0.8)(7,0)
\end{picture}
\end{center}

To shorten notations, we will write this factorization in the form 
\begin{equation}\label{eq:theta}
\Delta^2=\prod\limits_{i=1}^d\hat{Z}_{ii'}\cdot L_d\cdot \Theta\cdot
\prod\limits_{i=1}^d\check{Z}_{ii'},\end{equation}
and work only with the central
part $\Theta$, which geometrically corresponds to the upper half of the
considered diagram. Using the commutativity rules in the central part, one
can rewrite $\Theta$ as
$$\Theta=(L'_d)^3\cdot\Bigl(\prod_{i=1}^d Z_{i'\alpha}^2
\prod_{i=1}^d Z_{i'\beta}^2 \prod_{i=1}^d Z_{i'\gamma}^2\Bigr)^2\cdot
\left(Z_{\alpha\beta}^2 Z_{\alpha\gamma}^2Z_{\beta\gamma}^2\right).$$

Moving the second set of $Z^2_{i'\alpha}$ and $Z^2_{i'\beta}$ factors
to the left, one can rewrite this expression as
$$\Theta=(L'_d)^3\cdot\Bigl(\prod_{i=1}^d Z_{i'\alpha}^2\Bigr)^2\cdot
\prod_{i=1}^d Z_{i'\beta;\{\alpha\}}^2\cdot
\prod_{i=1}^d Z_{i'\beta}^2\cdot
\prod_{i=1}^d Z_{i'\gamma;\{\beta\}}^2\cdot
\prod_{i=1}^d Z_{i'\gamma}^2\cdot
\left(Z_{\alpha\beta}^2 Z_{\alpha\gamma}^2Z_{\beta\gamma}^2\right)$$
where $Z_{i'\beta;\{\alpha\}}$ and $Z_{i'\gamma;\{\beta\}}$ are 
half-twists along the following paths~:

\begin{center}
\setlength{\unitlength}{0.25in}
\begin{picture}(9,3)(3,-2)
\put(4,0){\circle*{0.1}}
\put(5,0){\circle*{0.1}}
\put(7,0){\circle*{0.1}}
\put(9,0){\circle*{0.1}}
\put(10,0){\circle*{0.1}}
\put(11,0){\circle*{0.1}}
\put(12,0){\circle*{0.1}}
\multiput(3.1,0)(0.2,0){3}{\circle*{0.02}}
\multiput(5.5,0)(0.2,0){6}{\circle*{0.02}}
\multiput(7.5,0)(0.2,0){6}{\circle*{0.02}}
\put(4,-0.2){\makebox(0,0)[ct]{$d$}}
\put(5.2,-0.2){\makebox(0,0)[ct]{1'}}
\put(7,-0.2){\makebox(0,0)[ct]{$i'$}}
\put(9,-0.2){\makebox(0,0)[ct]{$d'$}}
\put(10,-0.4){\makebox(0,0)[ct]{$\alpha$}}
\put(11,-0.3){\makebox(0,0)[ct]{$\beta$}}
\put(12,-0.4){\makebox(0,0)[ct]{$\gamma$}}
\qbezier(9.5,1)(11,1)(11,0)
\put(5.5,1){\line(1,0){4}}
\qbezier(5.5,1)(4.3,1)(4.3,0)
\qbezier(4.3,0)(4.3,-1.2)(5.3,-1.2)
\put(5.3,-1.2){\line(1,0){4}}
\qbezier(9.3,-1.2)(10.5,-1.2)(10.5,0)
\qbezier(10.5,0)(10.5,0.5)(10,0.5)
\qbezier(10,0.5)(9.5,0.5)(9.5,0)
\qbezier(9.5,0)(9.5,-0.9)(8.5,-0.9)
\put(5.5,-0.9){\line(1,0){3}}
\qbezier(5.5,-0.9)(4.6,-0.9)(4.6,0)
\qbezier(4.6,0)(4.6,0.7)(5.8,0.7)
\qbezier(5.8,0.7)(7,0.7)(7,0)
\put(7.5,-2){\makebox(0,0)[cb]{$Z_{i'\beta;\{\alpha\}}$}}
\end{picture}
\hfill
\begin{picture}(9,3)(3,-2)
\put(4,0){\circle*{0.1}}
\put(5,0){\circle*{0.1}}
\put(7,0){\circle*{0.1}}
\put(9,0){\circle*{0.1}}
\put(10,0){\circle*{0.1}}
\put(11,0){\circle*{0.1}}
\put(12,0){\circle*{0.1}}
\multiput(3.1,0)(0.2,0){3}{\circle*{0.02}}
\multiput(5.5,0)(0.2,0){6}{\circle*{0.02}}
\multiput(7.5,0)(0.2,0){6}{\circle*{0.02}}
\put(4,-0.2){\makebox(0,0)[ct]{$d$}}
\put(5.2,-0.2){\makebox(0,0)[ct]{1'}}
\put(7,-0.2){\makebox(0,0)[ct]{$i'$}}
\put(9,-0.2){\makebox(0,0)[ct]{$d'$}}
\put(10,-0.4){\makebox(0,0)[ct]{$\alpha$}}
\put(11,-0.3){\makebox(0,0)[ct]{$\beta$}}
\put(12,-0.4){\makebox(0,0)[ct]{$\gamma$}}
\qbezier(10.5,1)(12,1)(12,0)
\put(5.5,1){\line(1,0){5}}
\qbezier(5.5,1)(4.3,1)(4.3,0)
\qbezier(4.3,0)(4.3,-1.2)(5.3,-1.2)
\put(5.3,-1.2){\line(1,0){5}}
\qbezier(10.3,-1.2)(11.5,-1.2)(11.5,0)
\qbezier(11.5,0)(11.5,0.7)(10.5,0.7)
\qbezier(10.5,0.7)(9.5,0.7)(9.5,0)
\qbezier(9.5,0)(9.5,-0.9)(8.5,-0.9)
\put(5.5,-0.9){\line(1,0){3}}
\qbezier(5.5,-0.9)(4.6,-0.9)(4.6,0)
\qbezier(4.6,0)(4.6,0.7)(5.8,0.7)
\qbezier(5.8,0.7)(7,0.7)(7,0)
\put(7.5,-2){\makebox(0,0)[cb]{$Z_{i'\gamma;\{\beta\}}$}}
\end{picture}\end{center}

A succession of Hurwitz moves to the right makes it possible to rewrite
$\Theta$ as
$$(L'_d)^3\cdot\Bigl(\prod_{i=1}^d Z_{i'\alpha}^2\Bigr)^2\cdot
\prod_{i=1}^d Z_{i'\beta;\{\alpha\}}^2\cdot \tilde{Z}^2_{\alpha\beta,0}\cdot
\prod_{i=1}^d Z_{i'\beta}^2\cdot
\prod_{i=1}^d Z_{i'\gamma;\{\beta\}}^2\cdot
\tilde{Z}^2_{\alpha\gamma}\tilde{Z}^2_{\beta\gamma}\cdot
\prod_{i=1}^d Z_{i'\gamma}^2$$
where $\tilde{Z}_{\alpha\beta,0}$, $\tilde{Z}_{\alpha\gamma}$ and
$\tilde{Z}_{\beta\gamma}$ are
half-twists along the following paths~:

\begin{center}
\setlength{\unitlength}{0.25in}
\begin{picture}(5.2,3)(-0.6,-2)
\put(0,0){\circle*{0.1}}
\put(0.8,0){\circle*{0.1}}
\put(2.2,0){\circle*{0.1}}
\put(3,0){\circle*{0.1}}
\put(3.8,0){\circle*{0.1}}
\put(4.6,0){\circle*{0.1}}
\multiput(-0.6,0)(0.15,0){3}{\circle*{0.02}}
\multiput(1.12,0)(0.15,0){6}{\circle*{0.02}}
\put(0,-0.2){\makebox(0,0)[ct]{$d$}}
\put(1,-0.2){\makebox(0,0)[ct]{$1'$}}
\put(2.2,-0.2){\makebox(0,0)[ct]{$d'$}}
\put(3,0.6){\makebox(0,0)[ct]{$\alpha$}}
\put(3.8,-0.3){\makebox(0,0)[ct]{$\beta$}}
\put(4.6,-0.4){\makebox(0,0)[ct]{$\gamma$}}
\qbezier(3,0)(3,-1)(2,-1)
\put(1.4,-1){\line(1,0){0.6}}
\qbezier(1.4,-1)(0.4,-1)(0.4,0)
\qbezier(0.4,0)(0.4,1)(1.4,1)
\put(1.4,1){\line(1,0){1.4}}
\qbezier(2.8,1)(3.8,1)(3.8,0)
\put(2,-2){\makebox(0,0)[cb]{$\tilde{Z}_{\alpha\beta,0}$}}
\end{picture}
\hfil
\begin{picture}(5.2,3)(-0.6,-2)
\put(0,0){\circle*{0.1}}
\put(0.8,0){\circle*{0.1}}
\put(2.2,0){\circle*{0.1}}
\put(3,0){\circle*{0.1}}
\put(3.8,0){\circle*{0.1}}
\put(4.6,0){\circle*{0.1}}
\multiput(-0.6,0)(0.15,0){3}{\circle*{0.02}}
\multiput(1.12,0)(0.15,0){6}{\circle*{0.02}}
\put(0,-0.2){\makebox(0,0)[ct]{$d$}}
\put(1,-0.2){\makebox(0,0)[ct]{$1'$}}
\put(2.2,-0.2){\makebox(0,0)[ct]{$d'$}}
\put(3,0.6){\makebox(0,0)[ct]{$\alpha$}}
\put(3.8,-0.3){\makebox(0,0)[ct]{$\beta$}}
\put(4.6,-0.4){\makebox(0,0)[ct]{$\gamma$}}
\qbezier(3,0)(3,-1)(2,-1)
\put(1.4,-1){\line(1,0){0.6}}
\qbezier(1.4,-1)(0.4,-1)(0.4,0)
\qbezier(0.4,0)(0.4,1)(1.4,1)
\put(1.4,1){\line(1,0){2.2}}
\qbezier(3.6,1)(4.6,1)(4.6,0)
\put(2,-2){\makebox(0,0)[cb]{$\tilde{Z}_{\alpha\gamma}$}}
\end{picture}
\hfil
\begin{picture}(5.2,3)(-0.6,-2)
\put(0,0){\circle*{0.1}}
\put(0.8,0){\circle*{0.1}}
\put(2.2,0){\circle*{0.1}}
\put(3,0){\circle*{0.1}}
\put(3.8,0){\circle*{0.1}}
\put(4.6,0){\circle*{0.1}}
\multiput(-0.6,0)(0.15,0){3}{\circle*{0.02}}
\multiput(1.12,0)(0.15,0){6}{\circle*{0.02}}
\put(0,-0.2){\makebox(0,0)[ct]{$d$}}
\put(1,-0.2){\makebox(0,0)[ct]{$1'$}}
\put(2.2,-0.2){\makebox(0,0)[ct]{$d'$}}
\put(3,-0.4){\makebox(0,0)[ct]{$\alpha$}}
\put(3.8,-0.3){\makebox(0,0)[ct]{$\beta$}}
\put(4.6,-0.4){\makebox(0,0)[ct]{$\gamma$}}
\qbezier(2.6,0)(2.6,0.6)(3.2,0.6)
\qbezier(3.2,0.6)(3.8,0.6)(3.8,0)
\qbezier(2.6,0)(2.6,-1)(1.5,-1)
\qbezier(1.5,-1)(0.4,-1)(0.4,0)
\qbezier(0.4,0)(0.4,1)(1.4,1)
\put(1.4,1){\line(1,0){2.2}}
\qbezier(3.6,1)(4.6,1)(4.6,0)
\put(2,-2){\makebox(0,0)[cb]{$\tilde{Z}_{\beta\gamma}$}}
\end{picture}
\end{center}
Moving $\tilde{Z}^2_{\alpha\beta,0}$, $\tilde{Z}^2_{\alpha\gamma}$ and
$\tilde{Z}^2_{\beta\gamma}$ to the left, one can rewrite $\Theta$ as
$$\Theta=(L'_d)^3\cdot\Bigl(\prod_{i=1}^d Z_{i'\alpha}^2\Bigr)^2\cdot
\tilde{Z}^2_{\alpha\beta,0}\cdot
\Bigl(\prod_{i=1}^d Z_{i'\beta}^2\Bigr)^2\cdot
\tilde{Z}^2_{\alpha\gamma}\tilde{Z}^2_{\beta\gamma}\cdot
\Bigl(\prod_{i=1}^d Z_{i'\gamma}^2\Bigr)^2.$$
Moving the $Z_{i'\beta}^2$ factors to the left, one obtains the new
expression
$$\Theta=(L'_d)^3\cdot\Bigl(\prod_{i=1}^d Z_{i'\alpha}^2\Bigr)^2\cdot
\Bigl(\prod_{i=1}^d Z_{i'\beta}^2\Bigr)^2\cdot
\tilde{Z}^2_{\alpha\beta}\cdot
\tilde{Z}^2_{\alpha\gamma}\tilde{Z}^2_{\beta\gamma}\cdot
\Bigl(\prod_{i=1}^d Z_{i'\gamma}^2\Bigr)^2,$$
where $\tilde{Z}_{\alpha\beta}$ is a half-twist along the path

\begin{center}
\setlength{\unitlength}{0.25in}
\begin{picture}(5.2,2)(-0.6,-1)
\put(0,0){\circle*{0.1}}
\put(0.8,0){\circle*{0.1}}
\put(2.2,0){\circle*{0.1}}
\put(3,0){\circle*{0.1}}
\put(3.8,0){\circle*{0.1}}
\put(4.6,0){\circle*{0.1}}
\multiput(-0.6,0)(0.15,0){3}{\circle*{0.02}}
\multiput(1.12,0)(0.15,0){6}{\circle*{0.02}}
\put(0,-0.2){\makebox(0,0)[ct]{$d$}}
\put(1,-0.2){\makebox(0,0)[ct]{$1'$}}
\put(2.2,-0.2){\makebox(0,0)[ct]{$d'$}}
\put(3.1,0.4){\makebox(0,0)[ct]{$\alpha$}}
\put(3.8,-0.3){\makebox(0,0)[ct]{$\beta$}}
\put(4.6,-0.4){\makebox(0,0)[ct]{$\gamma$}}
\qbezier(2.6,0)(2.6,0.6)(3.2,0.6)
\qbezier(3.2,0.6)(3.8,0.6)(3.8,0)
\qbezier(2.6,0)(2.6,-1)(1.5,-1)
\qbezier(1.5,-1)(0.4,-1)(0.4,0)
\qbezier(0.4,0)(0.4,1)(1.4,1)
\put(1.4,1){\line(1,0){1.8}}
\qbezier(3.2,1)(4.2,1)(4.2,0)
\qbezier(4.2,0)(4.2,-1)(3.6,-1)
\qbezier(3.6,-1)(3,-1)(3,0)
\end{picture}
\end{center}

Observing that each factor $Z_{i'j'}^2$ in $L'_d$ commutes with the
products $\prod Z_{i'\alpha}^2$ and $\prod Z_{i'\beta}^2$
and also with $\tilde{Z}^2_{\alpha\beta}$, $\tilde{Z}^2_{\alpha\gamma}$
and $\tilde{Z}^2_{\beta\gamma}$,
a sequence of Hurwitz moves to the left makes it possible to rewrite
$\Theta$ as
\begin{equation}\label{eq:assembling2}
L'_d\cdot\Bigl(\prod_{i=1}^d Z_{i'\alpha}^2\Bigr)^2\cdot L'_d\cdot
\Bigl(\prod_{i=1}^d Z_{i'\beta}^2\Bigr)^2\cdot
\tilde{Z}^2_{\alpha\beta}\tilde{Z}^2_{\alpha\gamma}\tilde{Z}^2_{\beta\gamma}
\cdot L'_d\cdot \Bigl(\prod_{i=1}^d Z_{i'\gamma}^2\Bigr)^2.
\end{equation}

We now study more in detail the first part of (\ref{eq:assembling2}), namely
$$\Theta_{\alpha}=L'_d\cdot\Bigl(\prod_{i=1}^d Z_{i'\alpha}^2\Bigr)^2=
\prod_{i=1}^{d-1}\prod_{j=i+1}^d Z^2_{i'j'}
\cdot\Bigl(\prod_{i=1}^d Z_{i'\alpha}^2\Bigr)^2.$$
A sequence of Hurwitz moves to the right makes it possible to rewrite this
expression as
$$\Theta_{\alpha}=\prod_{i=1}^{d-1}\prod_{j=i+1}^d Z^2_{i'j'}
\cdot\prod_{i=1}^d \left(Z^2_{i'\alpha}\hat{Z}^2_{i'\alpha}\right),$$
where $\hat{Z}_{i'\alpha}$ is a half-twist along the path

\begin{center}
\setlength{\unitlength}{0.25in}
\begin{picture}(7,2)(-2,-1)
\put(-2,0){\circle*{0.1}}
\put(0,0){\circle*{0.1}}
\put(1,0){\circle*{0.1}}
\put(3,0){\circle*{0.1}}
\put(4,0){\circle*{0.1}}
\multiput(1.5,0)(0.2,0){6}{\circle*{0.02}}
\multiput(-1.5,0)(0.2,0){6}{\circle*{0.02}}
\put(-2,-0.2){\makebox(0,0)[ct]{$1'$}}
\put(0,-0.2){\makebox(0,0)[ct]{$i'$}}
\put(3,-0.2){\makebox(0,0)[ct]{$d'$}}
\put(4,-0.3){\makebox(0,0)[ct]{$\alpha$}}
\qbezier(0,0)(0,1)(2,1)
\qbezier(2,1)(4.7,1)(4.7,-0.1)
\qbezier(4.7,-0.1)(4.7,-1)(3,-1)
\qbezier(3,-1)(0.4,-1)(0.4,-0.2)
\qbezier(0.4,-0.2)(0.4,0.6)(2,0.6)
\qbezier(2,0.6)(4,0.6)(4,0)
\end{picture}
\end{center}
Using commutation relations, more Hurwitz moves yield the identity
$$\Theta_{\alpha}=\prod_{i=1}^{d}\Bigl(\prod_{j=i+1}^d Z^2_{i'j'}\cdot
Z^2_{i'\alpha}\hat{Z}^2_{i'\alpha}\Bigr).$$
Next we move $Z_{i'\alpha}^2$ to the left and obtain
$$\Theta_{\alpha}=
\prod_{i=1}^{d}\Bigl(Z^2_{i'\alpha}\cdot\prod_{j=i+1}^d Z^2_{i'j';(\alpha)}
\cdot\hat{Z}^2_{i'\alpha}\Bigr),$$
where $Z_{i'j';(\alpha)}=Z_{i'\alpha}^{-2}Z_{i'j'}Z_{i'\alpha}^2$ 
is a twist along the path
\begin{center}
\setlength{\unitlength}{0.25in}
\begin{picture}(9,2)(0,-1)
\put(2,0){\circle*{0.1}}
\put(4,0){\circle*{0.1}}
\put(6,0){\circle*{0.1}}
\put(7,0){\circle*{0.1}}
\put(8,0){\circle*{0.1}}
\put(9,0){\circle*{0.1}}
\multiput(0.5,0)(0.2,0){6}{\circle*{0.02}}
\multiput(2.5,0)(0.2,0){6}{\circle*{0.02}}
\multiput(4.5,0)(0.2,0){6}{\circle*{0.02}}
\put(2,-0.2){\makebox(0,0)[ct]{$i'$}}
\put(4,-0.2){\makebox(0,0)[ct]{$j'$}}
\put(6,-0.2){\makebox(0,0)[ct]{$d'$}}
\put(7,-0.4){\makebox(0,0)[ct]{$\alpha$}}
\put(8,-0.2){\makebox(0,0)[ct]{$\beta$}}
\put(9,-0.4){\makebox(0,0)[ct]{$\gamma$}}
\qbezier(2,0)(2,1)(3.5,1)
\qbezier(3.5,1)(4.75,1)(6,1)
\qbezier(6,1)(7.5,1)(7.5,0)
\qbezier(7.5,0)(7.5,-1)(7,-1)
\qbezier(7,-1)(6.5,-1)(6.5,0)
\qbezier(6.5,0)(6.5,0.5)(5.25,0.5)
\qbezier(5.25,0.5)(4,0.5)(4,0)
\end{picture}
\end{center}
Finally, moving the $Z^2_{i'j';(\alpha)}$ factors to the left, one obtains
the identity
$$\Theta_{\alpha}=
\prod_{i=1}^{d}\Bigl((Z^2_{i'\alpha})^2\cdot
\prod_{j=i+1}^d Z^2_{i'j';(\alpha)}\Bigr).$$

Geometrically this expression corresponds to the following picture~:

\begin{center}
\setlength{\unitlength}{0.25in}
\begin{picture}(6,1.8)(-3,-1.6)
\put(-3,0.02){\line(1,0){6}}
\qbezier(-2.3,-1.5)(-1.3,0)(-0.3,0)
\qbezier(1.7,-1.5)(0.7,0)(-0.3,0)
\qbezier(-3,-1.5)(-2,0)(-1,0)
\qbezier(1,-1.5)(0,0)(-1,0)
\qbezier(-1,-1.5)(0,0)(1,0)
\qbezier(3,-1.5)(2,0)(1,0)
\put(-3.5,-0.1){$\alpha$}
\put(-3.4,-1.6){$1'$}
\put(-2.7,-1.6){$2'$}
\put(-1.2,-0.7){...}
\put(-1.4,-1.6){$d'$}
\end{picture}
\end{center}

Proceeding similarly with the pieces involving $\beta$ and $\gamma$ in the
expression (\ref{eq:assembling2}), and letting 
$Z_{i'j';(\beta)}=Z_{i'\beta}^{-2}Z_{i'j'}Z_{i'\beta}^2$ and
$Z_{i'j';(\gamma)}=Z_{i'\gamma}^{-2}Z_{i'j'}Z_{i'\gamma}^2$
(these twists correspond to the same picture as $Z_{i'j';(\alpha)}$ but
going around $\beta$ or $\gamma$ instead of $\alpha$), the factorization
(\ref{eq:assembling2}) rewrites as
\begin{multline*}
\Theta=
\prod_{i=1}^{d}\Bigl((Z^2_{i'\alpha})^2\cdot
\prod_{j=i+1}^d Z^2_{i'j';(\alpha)}\Bigr)\cdot
\prod_{i=1}^{d}\Bigl((Z^2_{i'\beta})^2\cdot
\prod_{j=i+1}^d Z^2_{i'j';(\beta)}\Bigr)\cdot\\
\tilde{Z}^2_{\alpha\beta}
\tilde{Z}^2_{\alpha\gamma}\tilde{Z}^2_{\beta\gamma}\cdot 
\prod_{i=1}^{d}\Bigl((Z^2_{i'\gamma})^2\cdot
\prod_{j=i+1}^d Z^2_{i'j';(\gamma)}\Bigr).
\end{multline*}
We have finally achieved our goal of bringing next to each other the two 
intersections of each conic with each line. Therefore, going back to 
(\ref{eq:theta}), we finally obtain~:

\begin{prop}\label{prop:assembling}
The braid factorization corresponding to the union of $d$ conics and three
lines tangent to them is given by
\begin{multline}\label{eq:assembled}
\Delta^2=\prod\limits_{i=1}^d\hat{Z}_{ii'}\cdot L_d\cdot
\prod_{i=1}^{d}\Bigl(Z^4_{i'\alpha}\cdot
\prod_{j=i+1}^d Z^2_{i'j';(\alpha)}\Bigr)\cdot
\prod_{i=1}^{d}\Bigl(Z^4_{i'\beta}\cdot
\prod_{j=i+1}^d Z^2_{i'j';(\beta)}\Bigr)\cdot\\
\tilde{Z}^2_{\alpha\beta}
\tilde{Z}^2_{\alpha\gamma}\tilde{Z}^2_{\beta\gamma}\cdot 
\prod_{i=1}^{d}\Bigl(Z^4_{i'\gamma}\cdot
\prod_{j=i+1}^d Z^2_{i'j';(\gamma)}\Bigr)
\cdot\prod\limits_{i=1}^d\check{Z}_{ii'}.
\end{multline}
\end{prop}

As explained above, the connectedness of the space of configurations of
mutually tangent conics and lines implies that, for a different choice of
the initial configuration, the braid factorization remains the same up
to Hurwitz equivalence and global conjugation. For instance, using different
expressions as starting points for the above geometric derivation of
(\ref{eq:assembled}) leads to formulas in which the half-twists
$\hat{Z}_{ii'}$ and $\check{Z}_{ii'}$ are replaced by slightly different half-twists with the
same end points; these formulas are equivalent to (\ref{eq:assembled}) up
to suitable Hurwitz moves and conjugations.

For completeness, we now describe how the reader may re-obtain
the formula (\ref{eq:assembled}) by a direct calculation from the diagram
presented at the beginning of this section (we describe the case $d=2$, the
extension to all values of $d$ is trivial). We start again from the 
diagram representing the intersection of the configuration with
$\R^2\subset\C^2$.

\begin{center}
\setlength{\unitlength}{0.25in}
\begin{picture}(8,5)(-4,-0.3)
\qbezier(0.2,0)(-3.3,-0.5)(-3.8,1)
\qbezier(-3.8,1)(-4.3,2.5)(-0.8,3)
\qbezier(-0.8,3)(2.7,3.5)(3.2,2)
\qbezier(3.2,2)(3.7,0.5)(0.2,0)
\qbezier(0,0)(3.5,-0.5)(4,1)
\qbezier(4,1)(4.5,2.5)(1,3)
\qbezier(1,3)(-2.5,3.5)(-3,2)
\qbezier(-3,2)(-3.5,0.5)(0,0)
\qbezier(-5,1.4)(-2,3)(1,4.6)
\qbezier(5.2,1.4)(2.2,3)(-0.8,4.6)
\put(-4,3.16){\line(1,0){8.2}}
\put(-4.8,1.8){$\alpha$}
\put(4.5,2){$\gamma$}
\put(4.4,3.1){$\beta$}
\put(-1.8,0.5){$2$}
\put(1.5,0.5){$1$}
\put(-1.8,2.2){$1'$}
\put(1.5,2.2){$2'$}
\end{picture}
\end{center}

All the special points are sent to the real axis by the projection $\pi$, and
labelling them in the obvious manner they are, from left to right, in the
following order (after slightly deforming the projection in a manner which
clearly doesn't affect the braid factorization)~: $11'$, 
$22'$ (tangencies), $12$, $1'\alpha$, $1'2'$, $2'\alpha$, $\alpha\beta$, 
$2'\beta$, $1'2'$, $1'\beta$, $\alpha\gamma$, $\beta\gamma$,
$1'\gamma$, $1'2'$, $2'\gamma$ (nodes and double nodes), $11'$, $22'$ 
(tangencies). 

The base point is placed on the real axis, immediately to
the right of the first two tangencies (and to the left of all other points).
The intersection with the reference fiber differs from the expected one
by a permutation of the points labelled $1'$ and $2'$ (the points are in the 
order $1,2,2',1',\alpha,\beta,\gamma$)~; this is taken care of by 
conjugating all computed monodromies by a half-twist, namely the point 
labelled $1'$ is brought back to the left of $2'$ by moving it
counterclockwise along a half-circle passing {\it above} $2'$.

The system of generating loops that we use to define the braid factorization
is given by paths joining the base point to the various other points in
the following manner (one easily checks that these paths are ordered
counterclockwise around the base point). The first two paths join the base
point to the points $11'$ and $22'$ on its left, starting below the real
axis and rotating twice clockwise around $11'$ and $22'$ (see diagram
below). The four following paths join the base point to the points $12$, 
$1'\alpha$, $1'2'$ and $2'\alpha$ on its right, passing above the real axis.
The next four paths reach the points $1'\beta$, $1'2'$, $2'\beta$ and
$\alpha\beta$ in that order, starting above the real axis and crossing it
between $1'\beta$ and $\alpha\gamma$ to reach their end points from below,
as shown on the diagram. The following two paths join the base point to
$\alpha\gamma$ and $\beta\gamma$, simply passing above the real axis.
The next three paths have $1'\gamma$, $1'2'$ and $2'\gamma$ as end points,
passing above the real axis but circling once clockwise around the three
points before reaching them. Finally, the last two paths connect the base
point to the two rightmost points $11'$ and $22'$, passing above the real
axis and circling twice clockwise around them. The picture is as follows~:

\begin{center}
\setlength{\unitlength}{0.25in}
\begin{picture}(17.5,2)(-2,-1)
\multiput(-2,0)(1,0){16}{\circle*{0.1}}
\multiput(14.5,0)(1,0){2}{\circle*{0.1}}
\put(-2,-0.5){\makebox(0,0)[cb]{\small $11'$}}
\put(-1,-0.5){\makebox(0,0)[cb]{\small $22'$}}
\put(1,-0.5){\makebox(0,0)[cb]{\small $12$}}
\put(2,-0.5){\makebox(0,0)[cb]{\small $1'\alpha$}}
\put(3,-0.5){\makebox(0,0)[cb]{\small $1'2'$}}
\put(4,-0.5){\makebox(0,0)[cb]{\small $2'\alpha$}}
\put(5,0.1){\makebox(0,0)[cb]{\small $\alpha\beta$}}
\put(6,0.1){\makebox(0,0)[cb]{\small $2'\beta$}}
\put(7,0.1){\makebox(0,0)[cb]{\small $1'2'\vphantom{\beta}$}}
\put(7.8,0.1){\makebox(0,0)[cb]{\small $1'\beta$}}
\put(9,-0.5){\makebox(0,0)[cb]{\small $\alpha\smash{\gamma}$}}
\put(10,-0.5){\makebox(0,0)[cb]{\small $\smash{\beta\gamma}$}}
\put(11,-0.5){\makebox(0,0)[cb]{\small $1'\smash{\gamma}$}}
\put(12,-0.5){\makebox(0,0)[cb]{\small $1'2'$}}
\put(13,-0.5){\makebox(0,0)[cb]{\small $2'\smash{\gamma}$}}
\put(14.5,-0.5){\makebox(0,0)[cb]{\small $11'$}}
\put(15.5,-0.5){\makebox(0,0)[cb]{\small $22'$}}
\qbezier(0,0)(-0.5,-1)(-1.2,-1)
\qbezier(-0.5,0)(-0.5,-0.8)(-1.2,-0.8)
\put(-1.2,-1){\line(-1,0){0.6}}
\put(-1.2,-0.8){\line(-1,0){0.6}}
\qbezier(-1.8,-0.8)(-2.5,-0.8)(-2.5,0)
\qbezier(-1.8,-1)(-2.7,-1)(-2.7,0)
\qbezier(-2.5,0)(-2.5,0.7),(-1.5,0.3)
\qbezier(-2.7,0)(-2.7,0.85)(-1.5,0.85)
\qbezier(-1.5,0.85)(-0.5,0.85)(-0.5,0)
\qbezier(-2,0)(-1.75,0.15)(-1.5,0.3)
\qbezier(-1,0)(-1.25,0.15)(-1.5,0.3)
\put(0,0){\line(1,0){1}}
\qbezier(0,0)(0.5,0.5)(1.5,0.5)
\qbezier(1.5,0.5)(2,0.5)(2,0)
\qbezier(1.5,0.5)(2.5,0.5)(3,0)
\qbezier(1.5,0.5)(3,0.5)(4,0)
\qbezier(0,0)(0.5,1)(1.5,1)
\put(1.5,1){\line(1,0){13.5}}
\qbezier(7.5,1)(8.5,1)(8.5,0)
\qbezier(8.5,1)(9.5,1)(9.5,0.3)
\qbezier(10,0)(9.75,0.15)(9.5,0.3)
\qbezier(9,0)(9.25,0.15)(9.5,0.3)
\qbezier(8.5,0)(8.5,-0.6)(7,-0.6)
\put(7,-0.6){\line(0,1){0.6}}
\qbezier(7,-0.6)(7.5,-0.3)(8,0)
\qbezier(7,-0.6)(6.5,-0.3)(6,0)
\qbezier(7,-0.6)(6,-0.3)(5,0)
\qbezier(12.5,1)(13.5,1)(13.5,0)
\qbezier(13.5,0)(13.5,-1)(12.5,-1)
\put(11.5,-1){\line(1,0){1}}
\qbezier(11.5,-1)(10.4,-1)(10.4,0)
\qbezier(10.4,0)(10.4,0.7)(11.5,0.5)
\qbezier(11.5,0.5)(11.25,0.25)(11,0)
\qbezier(11.5,0.5)(11.75,0.25)(12,0)
\qbezier(11.5,0.5)(12.25,0.25)(13,0)
\qbezier(15,1)(16.2,1)(16.2,0)
\qbezier(16.2,0)(16.2,-1)(15,-1)
\qbezier(15,-1)(13.8,-1)(13.8,0)
\qbezier(13.8,0)(13.8,0.8)(15,0.8)
\qbezier(15,0.8)(16,0.8)(16,0)
\qbezier(16,0)(16,-0.8)(15,-0.8)
\qbezier(15,-0.8)(14,-0.8)(14,0)
\qbezier(14,0)(14,0.6)(15,0.5)
\qbezier(15,0.5)(14.75,0.25)(14.5,0)
\qbezier(15,0.5)(15.25,0.25)(15.5,0)
\end{picture}
\end{center}

The monodromy around each point is computed using the following
observation~: placing oneself along the real axis, close to the image in the
base of one of the special points, the intersection points of the curve with
the fiber of $\pi$ all lie along the real axis (except at the outermost
tangencies where some points have moved off the axis), and the two points
involved in the monodromy lie next to each other. The monodromy then
corresponds to a twist along a line segment between these two points~; more
importantly, restricting oneself to a half-circle around the considered
point in the base amounts to rotating the two points in the fiber around
each other by half the total angle. With this understood, and decomposing
each path into half-circles around the various points, the computations
simply become a tedious task of careful accounting.

After suitably conjugating by a half-twist between $1'$ and $2'$, it turns
out that the braid monodromies along the various given loops are 
exactly the factors appearing in (\ref{eq:assembled}), except in the case 
of the tangency points $11'$ and $22'$ at either extremity. In fact, the 
monodromies around the tangency points differ from $\hat{Z}_{ii'}$ and 
$\check{Z}_{ii'}$ by a conjugation by $Z_{12}^4$ (or more generally the 
square of $\Delta_d^2$ when $d>2$)~; a global conjugation of all factors
by this braid eliminates the discrepancy and yields the desired formula.

(The choices made for the two sets of $d$ tangencies may seem rather
artificial, and indeed other choices would lead in a slightly more direct
manner to equally valid
expressions -- Hurwitz and conjugation equivalent to (\ref{eq:assembled}) --
involving different half-twists instead of 
$\hat{Z}_{ii'}$ and $\check{Z}_{ii'}$. The choices made here are motivated
by consistency with the geometric derivation outlined before the statement
of Proposition~\ref{prop:assembling}, where these half-twists come up as
a consequence of the use of Proposition~\ref{prop:folding} as a starting
point for the calculation.)

\subsection{The degree doubling formula}\label{ss:braidthm}
We finally turn to our main objective, computing the braid factorization 
for $D_{2k}$. Recall from \S 2.2 that the generic covering map $f_{2k}$ can be
obtained as a small perturbation of $f'_{2k}=V'_2\circ f_k$, 
where $V'_2$ is a generic quadratic holomorphic map obtained by 
slightly perturbing $V_2^0:(x:y:z)\mapsto (x^2:y^2:z^2)$.
More precisely, Proposition 2 states that, away from the intersection
points of the two branch curves $R_k$ and $f_k^{-1}(R'_2)$, the map 
$f'_{2k}$ satisfies almost all expected properties, the only problem for 
the definition of braid monodromy invariants being that its branch curve 
is not everywhere transverse to itself~; of course, it is also necessary 
to perturb $f'_{2k}$ near the intersection points in order to obtain a 
generic local model.

Recall that, by the main result of \cite{AK}, $f'_{2k}$ can be made generic 
near the points of $\mathcal{I}'_k=R_k\cap f_k^{-1}(R'_2)$ by adding to it 
small perturbation terms (see also the argument at the end of \S 2.2).
Provided that the perturbations are chosen small enough, the transversality
properties satisfied by $f'_{2k}$ away from these points are not affected.
Moreover, recall that for large $k$ the one-parameter argument proving the 
uniqueness up to isotopy of quasiholomorphic coverings also implies the 
connectedness of the space of admissible perturbations of $f'_{2k}$ near a 
given point of $\mathcal{I}'_k$. Therefore, the 
perturbation of $f'_{2k}$ affects the braid monodromy near each
point of $\mathcal{I}'_k$ exactly as described in \S \ref{ss:3pt}.

It is important to observe that these perturbations of $f'_{2k}$ only 
significantly affect the branch curve near the points of $\mathcal{I}'_k$~:
away from $\mathcal{I}'_k$, the branch curve of the perturbed map remains
$C^1$-close to that of $f'_{2k}$ (the perturbation terms are very small
in comparison with the transversality estimates satisfied by $f'_{2k}$).
Therefore, no unexpected changes can take place in the braid monodromy, 
although some pairs of nodes may be created when self-transversality 
is lost.

Another seemingly crucial point to be understood is the manner in which
the $n$ copies of the branch curve of $V'_2$ are moved into mutually 
transverse positions. Indeed, as explained at the end of \S \ref{ss:3pt}
this information directly determines the contribution to the braid monodromy
of the points of $\mathcal{I}'_k$ by modifying the local configuration
of vanishing cycles. Similarly, the braid monodromy arising near the
points $(1\!:\!0\!:\!0)$, $(0\!:\!1\!:\!0)$ and $(0\!:\!0\!:\!1)$ from the
cusps and tangency points in the $n$ copies of $C_2$ is strongly related
to the local configuration in each group of $2n$ lines. Therefore, our
lack of control over the manner in which each of the three groups 
of $2n$ lines is arranged may seem rather disturbing. 

Fortunately, up to $m$-equivalence this does not affect the final outcome of
the calculations. Indeed, in most places the $2n$ components labelled 
$1_\alpha,\dots,n'_\alpha$ (or similarly the two other groups of $2n$ lines)
all lift into different sheets of the covering $f_{2k}:X\to\CP^2$~; the only
exceptions are near the intersection points of $C_2$ with $V'_2(D_k)$,
where two of the $2n$ curves actually meet each other (e.g., those 
labelled $p_\alpha$ and $q_\alpha$ in \S \ref{ss:3pt}), and similarly near 
the points of intersection between two groups of $2n$ lines, where the two
curves coming from the same copy of $C_2$ (e.g., those labelled
$i_\alpha$ and $i'_\alpha$) also merge. In any case, we are free to move
the various lines across each other, as long as the two distinguished
components are kept together~; in this process, the braid factorization
only changes when pairs of intersections are created or cancelled, which
always amounts to an $m$-equivalence. Observe moreover that all possible
configurations can be deformed into each other in this way~; this follows
e.g.\ from the fact that all the curves under consideration, whether
self-transverse or not, are locally braided. We conclude that up to 
$m$-equivalence the braid monodromy does not depend on the chosen 
configuration.

Another more algebraic way to express the same idea is the following.
As observed at the end of \S \ref{ss:3pt}, the manner in which the local
braid monodromy arising from a point of $\mathcal{I}'_k$ 
depends on the local configuration is a conjugation by an element $Q$ of 
$B_{2n}$ which after multiplication by an element of $B_{2n-2}\times B_2$ can 
easily be assumed to be a pure braid. Denoting by $\Phi$ the factorized
expression corresponding to the standard configuration and by $\Phi_Q$ its
conjugate by the braid $Q$, we have the chain of $m$-equivalences
$$\Phi_Q\sim Q\cdot Q^{-1}\cdot \Phi_Q \sim Q\cdot \Phi\cdot Q^{-1}\sim
Q\cdot Q^{-1}_\Phi\cdot \Phi,$$
where the first operation is a pair creation and the two others are Hurwitz
moves~; therefore, conjugating $\Phi$ by $Q$ is equivalent to inserting the
two factors $Q$ and $Q^{-1}_\Phi$, which are both pure braids in $B_{2n}$. 
A similar phenomenon occurs near the intersection points between two
groups of $2n$ lines~: the choice of a specific configuration amounts to
a conjugation by a pure braid in $B_{2n}\times B_{2n}$, which after a
suitable $m$-equivalence simply amounts to inserting some pure braids into
the factorization. Finally, some intersections between the $2n$ lines also 
occur outside of these points, which means that, independently of the issue of the local
configurations, some pure braids in $B_{2n}$ appear as factors. Collecting
all the pure braids in $B_{2n}$ we have obtained in this description, we get
that the choice of a specific configuration amounts to the choice of a set
of pure braid factors in $B_{2n}$ (or more precisely, three such sets of
factors, one for each of the groups of lines labelled $\alpha$, $\beta$ and 
$\gamma$). The product of these factors is always the same independently
of the chosen configuration, because in the end we only consider
factorizations of $\Delta^2$. The result then follows from the following
observation~: {\em given a pure braid $Q\in B_{2n}$, any two decompositions
of $Q$ into products of positive and negative twists differ from each other
by Hurwitz moves and pair cancellations.} This can be seen
by realizing a factorization of $Q$ as the braid monodromy of a
curve with $2n$ components in $\C^2$ and by observing that any two such
configurations are deformation equivalent (e.g., when
$Q=1$ the components can be unknotted by translating them).
See also \cite{AKS}.
\medskip

As explained in \S \ref{ss:folding}, we can deform the curve $D_k$ so that
its image by $V_2^0$ becomes arbitrarily close to a union of $d$ conics,
at which point the braid factorization for $V_2^0(D_k)$, or equivalently
$V'_2(D_k)$, is given by (\ref{eq:folding}). First consider the 
singular map $V_2^0\circ f_k$, whose branch curve is the union of 
$V_2^0(D_k)$ with three lines (each of which has multiplicity $2n$).
These three lines always intersect $V_2^0(D_k)$ tangently. Therefore, 
after slightly deforming the map $V_2^0$ so that the three lines composing
its branch curve avoid the pole of the projection $\pi$, the braid
factorization for the branch curve of $V_2^0\circ f_k$ is very close to that
given by Proposition \ref{prop:assembling}~; keeping in mind the result of
\S \ref{ss:folding}, the only difference between the braid monodromy for
$V_2^0\circ f_k$ and (\ref{eq:assembled}) is that the $L_d$ term in 
(\ref{eq:assembled}) should be replaced by the braid factorization $F_k$ 
for $D_k$.

The discussion at the beginning of this section gives a description of
the modifications that occur when $V_2^0$ is replaced by $V'_2$ and 
$f'_{2k}$ is perturbed into the generic map $f_{2k}$. In this situation, 
the lines labelled $\alpha$, $\beta$ and $\gamma$ in \S \ref{ss:assembling}
each need to be replaced by a set of $2n$ lines. As we know from our study
of the structure of $f_{2k}$ near the points of $\mathcal{I}'_k$, the 
factors $Z^4_{i'\alpha}$, $Z^4_{i'\beta}$ and $Z^4_{i'\gamma}$ in 
(\ref{eq:assembled}) need to be replaced by expressions similar to
(\ref{eq:mutual})~; as explained above we do not have to worry about the
details of the local configurations.

Moreover, the factors $\tilde{Z}^2_{\alpha\beta}$, 
$\tilde{Z}^2_{\alpha\gamma}$ and $\tilde{Z}^2_{\beta\gamma}$ in
(\ref{eq:assembled}) need to be replaced by the factorizations 
describing the behavior of $n$ copies of $C_2$ near one of the points
where two groups of $2n$ lines intersect each other. The contribution
of each copy of $C_2$ has been computed in \S \ref{ss:v2}, but we must
also take into account the mutual intersections between the various
components. Fortunately, as explained above we do not have to worry about
the exact local configuration, so we can choose one that simplifies
calculations. 

Finally, we also need consider the mutual intersections of 
the $2n$ lines labelled $1_\alpha,\dots,n'_\alpha$ (and similarly in the
two other groups)~; although the possibility of moving the lines across
each other gives a lot of freedom, the manner in which they intersect is 
largely determined by the twisting phenomena arising at the points of
intersection with $V'_2(D_k)$ or with the other groups of $2n$ lines.
Indeed, since the total braid monodromy for the branch curve of $f_{2k}$
has to be $\Delta^2$, the amount of twisting of any two lines around each
other, and more precisely the product of all the degree $\pm 2$ factors
involving $1_\alpha,\dots,n'_\alpha$, is entirely determined by the chosen 
configurations at the intersection points with $V'_2(D_k)$ and the other
groups of $2n$ lines. As observed above, the various 
possible decompositions of this product into degree $\pm 2$ factors are
all $m$-equivalent to each other, so that once again we can choose one
freely (more geometrically, it is quite clear that any two configurations
of the lines that are compatible with the local configurations chosen 
at the intersection points can be deformed into one another and hence 
yield $m$-equivalent results).
\medskip

We now need to explicitly describe the geometric monodromy representation
$\theta:F_d\to S_n$ for $f_k$. Recalling from \S \ref{ss:generalities} 
that all geometric morphisms $\theta:F_d\to S_n$ are
equivalent to each other up to conjugation, we are free to choose the one
most suited to our purposes~; since the choice that we now make is in some
particular cases not the most practical one, we will also explain how to
adapt the formula for a different choice of $\theta$.

Let us assume from now on that $n=\deg f_k$ and $d=\deg D_k$ satisfy the 
inequality $d\le n(n-1)$. This inequality is satisfied in almost all
examples~; in particular, given any symplectic 4-manifold, it is satisfied 
as soon as $k$ is large enough. Consider geometric generators
$\gamma_1,\dots,\gamma_d$ of $\pi_1(\CP^2-D_k)$~: the loops $\gamma_i$ are
contained in the reference fiber of the projection to $\CP^1$, in which,
assuming that the base point and the $d$ intersection points with $D_k$
all lie on the real axis, they join the base point to the $i$-th
intersection point by passing above the real axis, circle once
counterclockwise around the intersection point, and return to the base 
point along the same path.

Performing a suitable global conjugation of 
the braid monodromy of $f_k$ if necessary, we can assume that the 
geometric monodromy representation is such that the transpositions
$\theta(\gamma_1),\dots,\theta(\gamma_d)$ are respectively equal to
the $d$ first terms of the factorized expression
$$\mathrm{Id}=\prod_{i=1}^{n-1}\prod_{j=i+1}^{n} (ij)\,(ij)$$
in the symmetric group $S_n$. This choice is legal because $d$ is even
and $d\ge 2n-2$. For each $1\le i\le n(n-1)$ we define the two indices
$1\le p(i)<q(i)\le n$ such that the $i$-th factor of this expression
in $S_n$ is the transposition $(p(i)q(i))$~; in particular
$\theta(\gamma_i)=(p(i)q(i))$ for all $i\le d$.

We first consider the contribution of the intersection points of 
$V'_2(D_k)$ with $C_2$. Making the same choice of local configurations as 
in \S \ref{ss:3pt}, each factor $Z^4_{i'\alpha}$ in (\ref{eq:assembled})
needs to be replaced by 
\begin{multline}\label{eq:regen_alpha}
\prod_{j=n}^1\left(\acute{Z}^2_{i'j'_\alpha}
\bigl[\acute{Z}^2_{i'j_\alpha}\bigr]_{j\not\in\{p(i),q(i)\}}\right)
\cdot Z^3_{i'p(i)_\alpha}\cdot Z^3_{i'q(i)_\alpha}\cdot
Z^3_{i'p(i)_\alpha;(q(i)_\alpha)}\cdot\\
\smash{\prod_{j=n}^1}\left(\grave{Z}^2_{i'j'_\alpha}
\bigl[\grave{Z}^2_{i'j_\alpha}\bigr]_{j\not\in\{p(i),q(i)\}}\right),
\end{multline}
and similarly for the $Z^4_{i'\beta}$ and $Z^4_{i'\gamma}$ factors.

We next consider the intersections of the $2n$ lines labelled
$1_\alpha,\dots,n'_\alpha$ with the $2n$ lines
labelled $1_\beta,\dots,n'_\beta$. We choose as local configuration a
situation consisting of $n$ identical copies of $C_2$ shifted away from 
each other by generic translations. The amounts by which the various 
copies are translated away from each other are assumed to be much larger 
than the distance between the two lines in a pair (e.g., $i_\alpha$ and
$i'_\alpha$)~; although this configuration can no longer be considered 
as a very small perturbation of $f'_{2k}$, it is quite clear that the
translation process preserves the property of being locally braided, so
that in terms of braid monodromy this configuration is $m$-equivalent to
that obtained by a small perturbation of $f'_{2k}$. This choice of
configuration can be represented on the following diagram~:

\begin{center}
\setlength{\unitlength}{0.25in}
\begin{picture}(8,3.5)(0,0)
\put(0,0){\line(1,1){3.5}}
\put(0.15,0){\line(1,1){3.5}}
\put(1.5,0){\line(1,1){3.5}}
\put(1.65,0){\line(1,1){3.5}}
\put(4,0){\line(1,1){3.5}}
\put(4.15,0){\line(1,1){3.5}}
\put(0,0.5){\line(1,0){8}}
\put(0,0.6){\line(1,0){8}}
\put(0,1.5){\line(1,0){8}}
\put(0,1.6){\line(1,0){8}}
\put(0,3){\line(1,0){8}}
\put(0,3.1){\line(1,0){8}}
\multiput(4.6,2.05)(0.4,0.2){4}{\circle*{0.02}}
\put(0.62,0.55){\circle*{0.2}}
\put(3.12,1.55){\circle*{0.2}}
\put(7.12,3.05){\circle*{0.2}}
\put(-0.2,2.3){\makebox(0,0)[rt]{$\beta$}}
\put(2.5,0.2){\makebox(0,0)[ct]{$\alpha$}}
\end{picture}
\end{center}

In this picture each intersection along the diagonal corresponds to 
a copy of $C_2$, yielding an expression similar to that in 
(\ref{eq:v2_3pt}), while all other intersections occur between different 
copies of $C_2$ and simply yield nodes. However, recall from the
computations in \S \ref{ss:assembling} that, when inserted into the
expression for the global braid monodromy, all local braid monodromy 
contributions need to be conjugated in such a way that the various twists
are performed along paths similar to the one appearing in the definition of 
$\tilde{Z}_{\alpha\beta}$.
Therefore, if we momentarily ignore the specificities of the intersections
along the diagonal, the braid monodromy for nodal intersections between
the two sets of $2n$ lines should be given by
$$\prod_{i=1}^n \prod_{j=1}^n \left(\tilde{Z}^2_{i_\alpha j_\beta}
\tilde{Z}^2_{i_\alpha j'_\beta}
\tilde{Z}^2_{i'_\alpha j_\beta}
\tilde{Z}^2_{i'_\alpha j'_\beta}\right),$$
where for any $\tau\in\{1_\alpha,1'_\alpha,\dots,
n_\alpha,n'_\alpha\}$ and $\upsilon\in\{1_\beta,1'_\beta,\dots,
n_\alpha,n'_\beta\}$ the notation $\tilde{Z}_{\tau\upsilon}$ represents
a half-twist along the path

\begin{center}
\setlength{\unitlength}{0.25in}
\begin{picture}(17.5,2)(0,-1)
\multiput(0.5,0)(0.2,0){6}{\circle*{0.02}}
\put(2,0){\circle*{0.1}}
\put(2,-0.2){\makebox(0,0)[ct]{$d$}}
\put(3,0){\circle*{0.1}}
\put(3,-0.2){\makebox(0,0)[ct]{$1'$}}
\multiput(3.5,0)(0.2,0){6}{\circle*{0.02}}
\put(5,0){\circle*{0.1}}
\put(5,-0.2){\makebox(0,0)[ct]{$d'$}}
\put(6,0){\circle*{0.1}}
\put(6.1,-0.2){\makebox(0,0)[ct]{$1_\alpha$}}
\multiput(6.5,0)(0.2,0){6}{\circle*{0.02}}
\put(8,0){\circle*{0.1}}
\put(8,-0.3){\makebox(0,0)[ct]{$\tau$}}
\multiput(8.5,0)(0.2,0){6}{\circle*{0.02}}
\put(10,0){\circle*{0.1}}
\put(10,-0.1){\makebox(0,0)[ct]{$n'_\alpha$}}
\put(11,0){\circle*{0.1}}
\put(11,-0.2){\makebox(0,0)[ct]{$1_\beta$}}
\multiput(11.5,0)(0.2,0){6}{\circle*{0.02}}
\put(13,0){\circle*{0.1}}
\put(13,-0.3){\makebox(0,0)[ct]{$\upsilon$}}
\multiput(13.5,0)(0.2,0){6}{\circle*{0.02}}
\put(15,0){\circle*{0.1}}
\put(14.9,-0.1){\makebox(0,0)[ct]{$n'_\beta$}}
\put(16,0){\circle*{0.1}}
\put(16,-0.2){\makebox(0,0)[ct]{$1_\gamma$}}
\multiput(16.5,0)(0.2,0){6}{\circle*{0.02}}
\qbezier(8,0)(8,0.5)(6.8,0.5)
\qbezier(6.8,0.5)(5.6,0.5)(5.6,-0.2)
\qbezier(5.6,-0.2)(5.6,-1)(6.6,-1)
\put(6.6,-1){\line(1,0){7.9}}
\qbezier(14.5,-1)(15.5,-1)(15.5,0)
\qbezier(15.5,0)(15.5,1)(14.5,1)
\put(3.5,1){\line(1,0){11}}
\qbezier(3.5,1)(2.5,1)(2.5,0)
\qbezier(2.5,0)(2.5,-1)(3.5,-1)
\put(3.5,-1){\line(1,0){1}}
\qbezier(4.5,-1)(5.2,-1)(5.4,0)
\qbezier(5.4,0)(5.54,0.7)(6.5,0.7)
\put(6.5,0.7){\line(1,0){5.5}}
\qbezier(12,0.7)(13,0.7)(13,0)
\end{picture}
\end{center}

However, according to the calculations performed in \S \ref{ss:v2}, the 
intersections corresponding to $i=j$ consist of three cusps and one tangency 
point set up as in (\ref{eq:v2_3pt}) rather than four nodes. Therefore, the
correct contribution to the braid factorization of $f_{2k}$ is
given by the expression 
\begin{multline}\label{eq:nv2_2}
\prod_{i=1}^n \Bigl(\prod_{j=1}^{i-1} \left(\tilde{Z}^2_{i_\alpha j_\beta}
\tilde{Z}^2_{i_\alpha j'_\beta}
\tilde{Z}^2_{i'_\alpha j_\beta}
\tilde{Z}^2_{i'_\alpha j'_\beta}\right) \cdot
\tilde{Z}^3_{i_\alpha i_\beta}\tilde{Z}^3_{i_\alpha i'_\beta}
\tilde{Z}_{i_\alpha i'_\alpha;(i_\beta i'_\beta)}
\tilde{Z}^3_{i'_\alpha i_\beta}\cdot\\
\prod_{j=i+1}^{n} \left(\tilde{Z}^2_{i_\alpha j_\beta}
\tilde{Z}^2_{i_\alpha j'_\beta}
\tilde{Z}^2_{i'_\alpha j_\beta}
\tilde{Z}^2_{i'_\alpha j'_\beta}\right)\Bigr),
\end{multline}
where $\tilde{Z}_{i_\alpha i'_\alpha;(i_\beta i'_\beta)}=
(\tilde{Z}^2_{i'_\alpha i_\beta}\tilde{Z}^2_{i'_\alpha i'_\beta})
Z_{i_\alpha i'_\alpha}
(\tilde{Z}^2_{i'_\alpha i_\beta}\tilde{Z}^2_{i'_\alpha i'_\beta})^{-1}$ is
a half-twist exchanging $i_\alpha$ and $i'_\alpha$ along a path that goes
around $i_\beta$ and $i'_\beta$ (the $\alpha$ points being connected to the
$\beta$ points along the same type of path described above). 

The factors $\tilde{Z}^2_{\alpha\gamma}$ and $\tilde{Z}^2_{\beta\gamma}$
in (\ref{eq:assembled}) are treated similarly, and give rise to expressions
similar to (\ref{eq:nv2_2}), except that the paths along which the 
$\tilde{Z}^2_{\tau\upsilon}$ factors twist now follow the model of
$\tilde{Z}^2_{\alpha\gamma}$ or $\tilde{Z}^2_{\beta\gamma}$ instead of
$\tilde{Z}^2_{\alpha\beta}$.

Our choice of local configuration for the $\alpha-\beta$ intersection is 
rather arbitrary~; however, a different choice would only affect the braid 
factorization by conjugation by a pure braid in $B_{2n}\times B_{2n}$ (each 
factor acting on one group of lines, while the path along which 
the groups are connected to each other necessarily remains that of 
$\tilde{Z}_{\alpha\beta}$). By the argument at the beginning of this section,
such a conjugation amounts up to $m$-equivalence to inserting some pure 
braid factors in $B_{2n}\times B_{2n}$ into the global braid monodromy, 
which has been shown not to affect the outcome of the computations, so that 
we can safely ignore this issue.
\medskip

We now look at the remaining nodal intersections between the $2n$ lines 
$1_\alpha,1'_\alpha,\dots,n_\alpha,n'_\alpha$. The product of all these
contributions to the braid monodromy is determined in the
following manner by the previously chosen configurations at intersection 
points with $V'_2(D_k)$ and with the other groups of $2n$ lines.
If we consider only the relative motions of the $2n$ points labelled
$1_\alpha,\dots,n'_\alpha$ induced by the various braids in the
factorization, it is easy to check from the above formulas that the tangent 
intersection with the line labelled $i'$ in $V'_2(D_k)$ contributes a
half-twist $Z_{p(i)_\alpha q(i)_\alpha}$ for all $1\le i\le d$, while
the intersection of $i_\alpha$ and $i'_\alpha$ with $i_\beta$ and $i'_\beta$
(or similarly $i_\gamma$ and $i'_\gamma$) contributes the half-twist 
$Z_{i_\alpha i'_\alpha}$. Therefore, the total contribution of intersection
points is equal to $\prod_{i=1}^d Z_{p(i)_\alpha q(i)_\alpha}\cdot
\left(\prod_{i=1}^n Z_{i_\alpha i'_\alpha}\right)^2$. 

On the other hand, recalling that we are looking for the braid factorization
of a curve in $\CP^2$, the overall relative motions of the 
$2n$ points $1_\alpha,\dots,n'_\alpha$ around each other must amount exactly
to the central element $\Delta_{2n}^2$ in $B_{2n}$~; 
the contribution of the additional nodal intersections is therefore exactly
the difference between the contribution of intersection points and 
$\Delta_{2n}^2$. Moreover, recall
from the discussion at the beginning of this section that the decomposition
of this contribution into a product of positive and negative twists is 
unique up to $m$-equivalence. In order to explicitly compute this
decomposition, we first derive a suitable expression of $\Delta_{2n}^2$.
Viewing the $2n$ points $1_\alpha,1'_\alpha,\dots,n_\alpha,n'_\alpha$
as $n$ groups of two points, it is easy to check that the full twist 
$\Delta_{2n}^2$ can be expressed as 
\begin{equation}\label{eq:2nlines1}
\Delta_{2n}^2=\prod_{i=1}^{n-1}\prod_{j=i+1}^n (Z^2_{i_\alpha j_\alpha}
Z^2_{i_\alpha j'_\alpha} Z^2_{i'_\alpha j_\alpha} Z^2_{i'_\alpha j'_\alpha})
\cdot \prod_{i=1}^n Z^2_{i_\alpha i'_\alpha}.
\end{equation}
Note that the two parts of this expression can be exchanged by Hurwitz
moves. The second part of (\ref{eq:2nlines1}) corresponds exactly to the
contribution of the intersection points with the two other groups of $2n$
lines~; meanwhile, the first $d/2$ factors $Z^2_{i_\alpha j_\alpha}$ 
correspond to the contribution of the points of $\mathcal{I}'_k$ 
(recall the choice of geometric monodromy representation
made above). Therefore, the nodal intersections correspond exactly
to the remaining factors in (\ref{eq:2nlines1}). Inserting these
braids at their respective positions in the factorization, and bringing
the $\check{Z}_{ii'}$ factors back to the beginning of the factorization by
Hurwitz moves, we finally obtain the following result~:

\begin{theo}
Let $X$ be a compact symplectic 4-manifold, and let $f_k:X\to\CP^2$ be
an approximately holomorphic branched covering given by three sections
of $L^{\otimes k}$. Denote by $D_k$ the branch curve of $f_k$, and let
$d=\deg D_k$ and $n=\deg f_k$. Assume that $d\le n(n-1)$. Denote by
$F_k$ the braid factorization corresponding to $D_k$, and assume
that the geometric monodromy representation $\theta:\pi_1(\CP^2-D_k)\to S_n$
is as described at the beginning of \S \ref{ss:braidthm}, i.e.\ 
$\theta(\gamma_i)=(p(i)q(i))$. Then, with the 
notations described in \S \ref{ss:3pt}, the braid factorization 
corresponding to the branch curve $D_{2k}$ of $f_{2k}$ is given up to
$m$-equivalence by the following formula, provided that $k$ is large enough~:
\begin{equation}\label{eq:k2k}
\Delta_{2d+6n}^2=
\mathbf{T}_d\cdot \iota(F_k)\cdot \mathbf{I}^\alpha_{d,n,\theta}
\cdot \mathbf{I}^\beta_{d,n,\theta}\cdot\mathbf{V}^{\alpha\beta}_n
\cdot \mathbf{V}^{\alpha\gamma}_n\cdot \mathbf{V}^{\beta\gamma}_n
\cdot \mathbf{I}^\gamma_{d,n,\theta},
\end{equation}
where $\iota:B_d\hookrightarrow B_{2d+6n}$ is the natural embedding given
by considering a disc containing the $d$ leftmost points
(labelled $1,\dots,d$),

\begin{eqnarray*}
\mathbf{T}_d & = & \prod_{i=1}^d\check{Z}_{ii'}\cdot
\prod_{i=1}^d\hat{Z}_{ii'},
\\
\mathbf{I}^\alpha_{d,n,\theta} &=& 
\prod_{i=1}^{d}\Biggl(
\prod_{j=n}^1\left(\acute{Z}^2_{i'j'_\alpha}
\bigl[\acute{Z}^2_{i'j_\alpha}\bigr]_{j\not\in\{p(i),q(i)\}}\right)
\cdot Z^3_{i'p(i)_\alpha}\cdot Z^3_{i'q(i)_\alpha}\cdot\\
&& \quad Z^3_{i'p(i)_\alpha;(q(i)_\alpha)}\cdot
\prod_{j=n}^1\left(\grave{Z}^2_{i'j'_\alpha}
\bigl[\grave{Z}^2_{i'j_\alpha}\bigr]_{j\not\in\{p(i),q(i)\}}\right)
\cdot\prod_{j=i+1}^d Z^2_{i'j';(\alpha)}\cdot\qquad \\ 
&& \qquad \left[Z^2_{p(i)_\alpha q(i)'_\alpha} 
Z^2_{p'(i)_\alpha q(i)_\alpha} Z^2_{p'(i)_\alpha q'(i)_\alpha}
\right]_{i\equiv 0\ \mathrm{mod}\ 2}\Biggr)\cdot\\
&& \prod_{i=(d/2)+1}^{n(n-1)/2} \left(Z^2_{p(2i)_\alpha q(2i)_\alpha}
Z^2_{p(2i)_\alpha q(2i)'_\alpha} 
Z^2_{p'(2i)_\alpha q(2i)_\alpha} Z^2_{p'(2i)_\alpha q'(2i)_\alpha}\right),
\end{eqnarray*}
\begin{eqnarray*}
\mathbf{I}^\beta_{d,n,\theta} &=& 
\prod_{i=1}^{d}\Biggl(
\prod_{j=n}^1\left(\acute{Z}^2_{i'j'_\beta}
\bigl[\acute{Z}^2_{i'j_\beta}\bigr]_{j\not\in\{p(i),q(i)\}}\right)
\cdot Z^3_{i'p(i)_\beta}\cdot Z^3_{i'q(i)_\beta}\cdot\\
&& \quad Z^3_{i'p(i)_\beta;(q(i)_\beta)}\cdot
\prod_{j=n}^1\left(\grave{Z}^2_{i'j'_\beta}
\bigl[\grave{Z}^2_{i'j_\beta}\bigr]_{j\not\in\{p(i),q(i)\}}\right)
\cdot\prod_{j=i+1}^d Z^2_{i'j';(\beta)}\cdot\qquad \\
&& \qquad \left[Z^2_{p(i)_\beta q(i)'_\beta} 
Z^2_{p'(i)_\beta q(i)_\beta} Z^2_{p'(i)_\beta q'(i)_\beta}
\right]_{i\equiv 0\ \mathrm{mod}\ 2}\Biggr)\cdot\\
&& \prod_{i=(d/2)+1}^{n(n-1)/2} \left(Z^2_{p(2i)_\beta q(2i)_\beta}
Z^2_{p(2i)_\beta q(2i)'_\beta}
Z^2_{p'(2i)_\beta q(2i)_\beta} Z^2_{p'(2i)_\beta q'(2i)_\beta}\right),
\end{eqnarray*}
\begin{eqnarray*}
\mathbf{I}^\gamma_{d,n,\theta} &=& 
\prod_{i=1}^{d}\Biggl(
\prod_{j=n}^1\left(\acute{Z}^2_{i'j'_\gamma}
\bigl[\acute{Z}^2_{i'j_\gamma}\bigr]_{j\not\in\{p(i),q(i)\}}\right)
\cdot Z^3_{i'p(i)_\gamma}\cdot Z^3_{i'q(i)_\gamma}\cdot\\
&& \quad Z^3_{i'p(i)_\gamma;(q(i)_\gamma)}\cdot
\prod_{j=n}^1\left(\grave{Z}^2_{i'j'_\gamma}
\bigl[\grave{Z}^2_{i'j_\gamma}\bigr]_{j\not\in\{p(i),q(i)\}}\right)
\cdot\prod_{j=i+1}^d Z^2_{i'j';(\gamma)}\cdot \qquad \\ 
&& \qquad \left[Z^2_{p(i)_\gamma q(i)'_\gamma} 
Z^2_{p'(i)_\gamma q(i)_\gamma} Z^2_{p'(i)_\gamma q'(i)_\gamma}
\right]_{i\equiv 0\ \mathrm{mod}\ 2}\Biggr)\cdot \\
&& \prod_{i=(d/2)+1}^{n(n-1)/2} \left(Z^2_{p(2i)_\gamma q(2i)_\gamma}
Z^2_{p(2i)_\gamma q(2i)'_\gamma}
Z^2_{p'(2i)_\gamma q(2i)_\gamma} Z^2_{p'(2i)_\gamma q'(2i)_\gamma}\right),
\end{eqnarray*}
\begin{eqnarray*}
\mathbf{V}^{\alpha\beta}_n &=&
\prod_{i=1}^n
\Biggl(\prod_{j=1}^{i-1} \left(\tilde{Z}^2_{i_\alpha j_\beta}
\tilde{Z}^2_{i_\alpha j'_\beta}
\tilde{Z}^2_{i'_\alpha j_\beta}
\tilde{Z}^2_{i'_\alpha j'_\beta}\right) \cdot
\tilde{Z}^3_{i_\alpha i_\beta}\tilde{Z}^3_{i_\alpha i'_\beta}
\tilde{Z}_{i_\alpha i'_\alpha;(i_\beta i'_\beta)}
\tilde{Z}^3_{i'_\alpha i_\beta}\cdot\qquad\\
&& \hspace{4cm} \prod_{j=i+1}^{n} \left(\tilde{Z}^2_{i_\alpha j_\beta}
\tilde{Z}^2_{i_\alpha j'_\beta}
\tilde{Z}^2_{i'_\alpha j_\beta}
\tilde{Z}^2_{i'_\alpha j'_\beta}\right)\Biggr),
\end{eqnarray*}
\begin{eqnarray*}
\mathbf{V}^{\alpha\gamma}_n &=&
\prod_{i=1}^n
\Biggl(\prod_{j=1}^{i-1} \left(\tilde{Z}^2_{i_\alpha j_\gamma}
\tilde{Z}^2_{i_\alpha j'_\gamma}
\tilde{Z}^2_{i'_\alpha j_\gamma}
\tilde{Z}^2_{i'_\alpha j'_\gamma}\right) \cdot
\tilde{Z}^3_{i_\alpha i_\gamma}\tilde{Z}^3_{i_\alpha i'_\gamma}
\tilde{Z}_{i_\alpha i'_\alpha;(i_\gamma i'_\gamma)}
\tilde{Z}^3_{i'_\alpha i_\gamma}\cdot\qquad\\
&& \hspace{4cm} \prod_{j=i+1}^{n} \left(\tilde{Z}^2_{i_\alpha j_\gamma}
\tilde{Z}^2_{i_\alpha j'_\gamma}
\tilde{Z}^2_{i'_\alpha j_\gamma}
\tilde{Z}^2_{i'_\alpha j'_\gamma}\right)\Biggr),
\end{eqnarray*}
\begin{eqnarray*}
\mathbf{V}^{\beta\gamma}_n &=&
\prod_{i=1}^n
\Biggl(\prod_{j=1}^{i-1} \left(\tilde{Z}^2_{i_\beta j_\gamma}
\tilde{Z}^2_{i_\beta j'_\gamma}
\tilde{Z}^2_{i'_\beta j_\gamma}
\tilde{Z}^2_{i'_\beta j'_\gamma}\right) \cdot
\tilde{Z}^3_{i_\beta i_\gamma}\tilde{Z}^3_{i_\beta i'_\gamma}
\tilde{Z}_{i_\beta i'_\beta;(i_\gamma i'_\gamma)}
\tilde{Z}^3_{i'_\beta i_\gamma}\cdot\qquad\\
&& \hspace{4cm} \prod_{j=i+1}^{n} \left(\tilde{Z}^2_{i_\beta j_\gamma}
\tilde{Z}^2_{i_\beta j'_\gamma}
\tilde{Z}^2_{i'_\beta j_\gamma}
\tilde{Z}^2_{i'_\beta j'_\gamma}\right)\Biggr).
\end{eqnarray*}
In these expressions, the notation $[\dots]_{i\equiv 0\ \mathrm{mod}\ 2}$
means that the enclosed factors are only present for even values of $i$~;
the various notations for braids correspond to half-twists along the following
paths~:

\begin{center}
\epsfbox{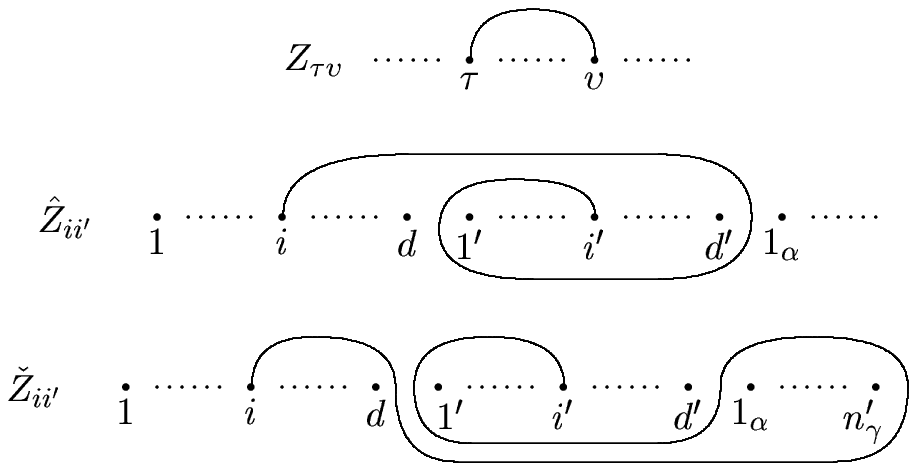}
\end{center}
\begin{center}
\epsfbox{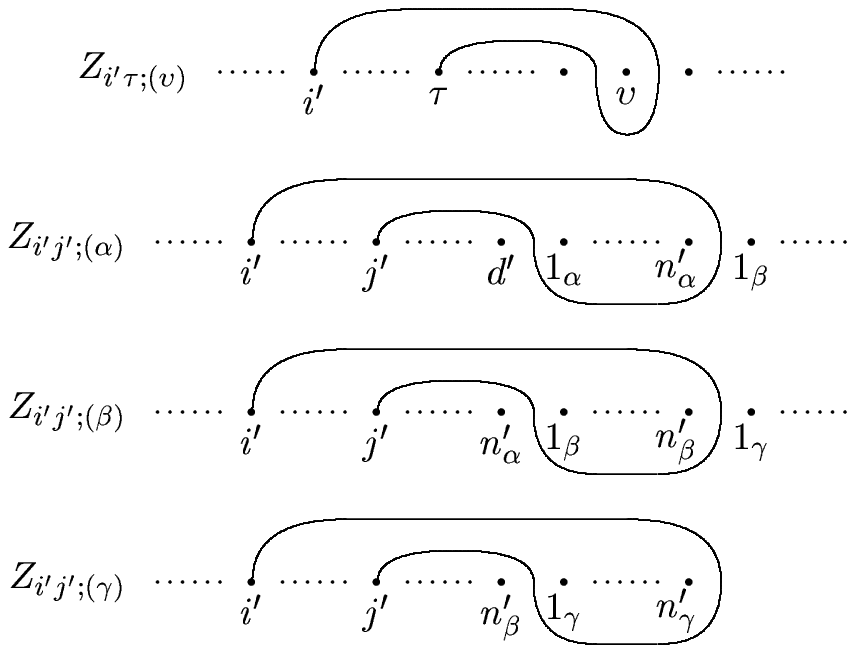}
\end{center}
\begin{center}
\epsfbox{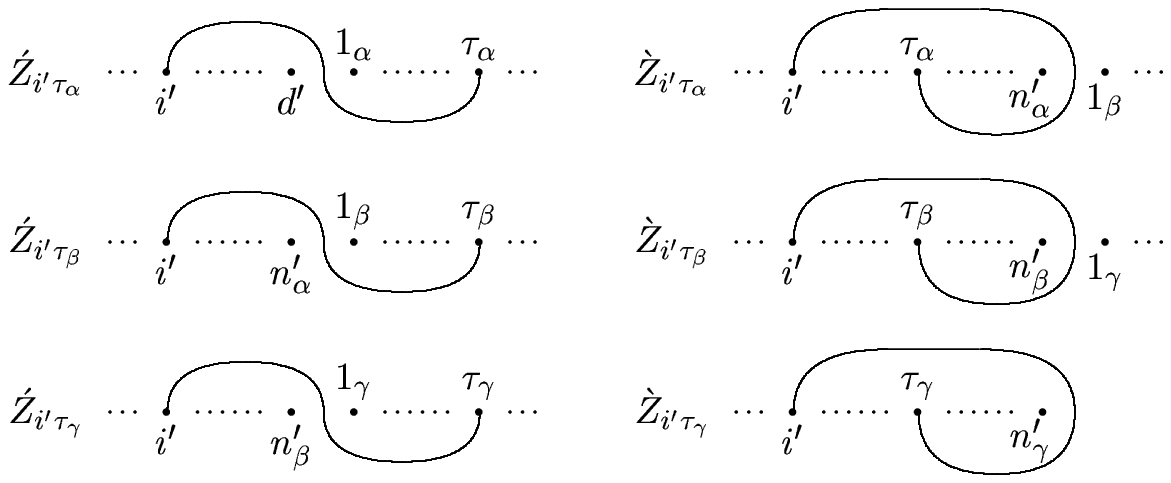}
\end{center}
\begin{center}
\epsfbox{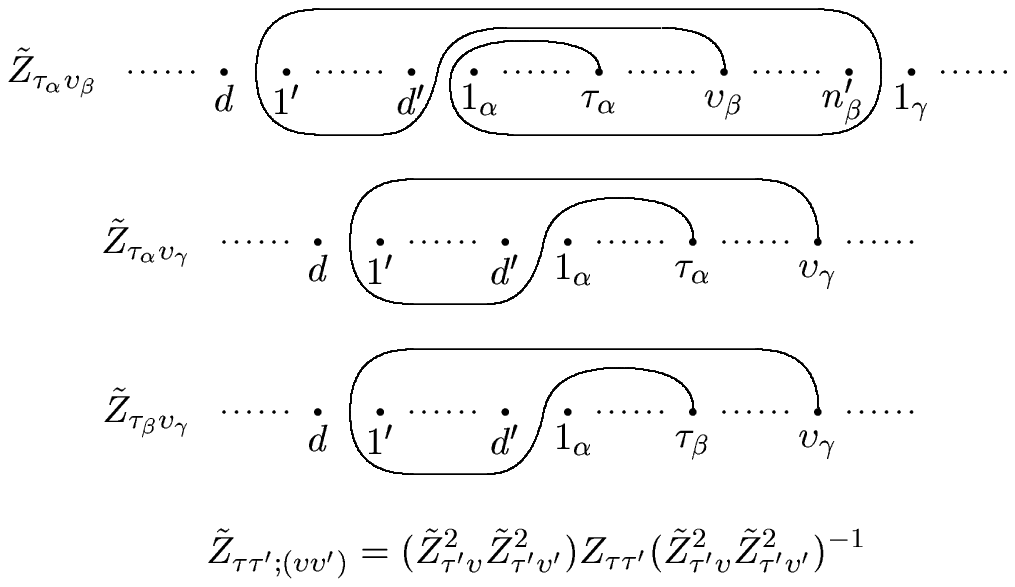}
\end{center}
\end{theo}

\noindent {\bf Remark~:} in the expression (\ref{eq:k2k}) we have made use
of our specific choice of geometric monodromy representation for $f_k$,
which requires the inequality $d\le n(n-1)$ to hold in counterpart for the
relative simplicity of the resulting factors. Also, we have chosen to
insert some of the pure braid factors involving the $2n$ lines
$1_\alpha,\dots,n'_\alpha$ amid the contributions of the intersection
points of these lines with $V'_2(D_k)$, in order to avoid the need for
a rewriting of (\ref{eq:2nlines1}) using Hurwitz moves to isolate these
contributions.

In general, if one wishes to get rid of the assumption made on the structure
of the geometric monodromy representation $\theta$ and to remove the 
constraint $d\le n(n-1)$, the necessary modifications are rather easy and
only involve finding a different expression of $\Delta_{2n}^2$ to 
replace (\ref{eq:2nlines1}). Namely, denote by $(\tau(i)\upsilon(i))\in S_n$
the image by $\theta$ of the $i$-th geometric generator of 
\mbox{$\pi_1(\CP^2-D_k)$} (in the standard situation of (\ref{eq:k2k}) one has 
$\tau(i)=p(i)$ and $\upsilon(i)=q(i)$ but we now want to lift this 
assumption). Then, if we keep our choice of the simplest local geometric
configurations at points of $\mathcal{I}'_k$, the contribution of these
points to the twisting among the lines $1_\alpha,\dots,n'_\alpha$ is given
by the pure braid $\prod_{i=1}^d Z_{\tau(i)_\alpha \upsilon(i)_\alpha}$.
We know that the total contribution of nodal intersections between the
$2n$ lines must be equal to 
$$Q_\alpha=\Bigl(\prod_{i=1}^d Z_{\tau(i)_\alpha\upsilon(i)_\alpha}
\Bigr)^{-1}\cdot\Delta^2_{2n}\cdot
\Bigl(\prod_{i=1}^n Z^2_{i_\alpha i'_\alpha}\Bigr)^{-1}.$$
Since $Q_\alpha$ is a pure braid it can be decomposed into a product of
positive and negative twists involving $1_\alpha,\dots,n'_\alpha$. The
resulting modification of the factors in $\mathbf{I}^\alpha_{d,n,\theta}$
is as follows~: in the first two lines, $p(i)$ and 
$q(i)$ should be replaced by $\tau(i)$ and $\upsilon(i)$ respectively~; 
the third line, consisting only of nodal intersections inserted amid
the other contributions, should be deleted~; and the last line, containing
the main group of nodal intersections, should be replaced by the chosen
factorization of $Q_\alpha$. Similar modifications are also required in
$\mathbf{I}^\beta_{d,n,\theta}$ and $\mathbf{I}^\gamma_{d,n,\theta}$.

As explained previously, the independence of the braid factorization upon
the choice of local configurations and the fact that any two geometric
monodromy representations differ from each other by a global conjugation
imply that the expression obtained for a non-standard choice of $\theta$ 
is m-equivalent to the standard one. In particular, the possible presence
of negative twists in the factorization of $Q_\alpha$ should not be
considered as an indication of the existence of non-removable negative nodes.
\medskip

{\bf Remark~:} when $X$ is a complex projective manifold, braid monodromy
becomes well-defined up to Hurwitz equivalence and global conjugation only,
since no negative nodes may appear in the (holomorphic) branch curve.
However, (\ref{eq:k2k}) only gives the answer up to $m$-equivalence even
in this case. If one looks more closely, the deformation process described in 
\S \ref{ss:folding} can be performed algebraically provided that 
$L^{\otimes k}$ is sufficiently positive, and therefore remains
valid in the complex setting~; in fact, all the braid monodromy computations
described in \S\S \ref{ss:folding}--\ref{ss:assembling} are
valid not only up to $m$-equivalence but also up to Hurwitz equivalence
and conjugation. However, what is not clear from an algebraic point of view
is the exact configuration in which the lines $1_\alpha,\dots,n'_\alpha$ 
are placed by a generic algebraic perturbation performed near the points of 
$\mathcal{I}'_k$. Determining this information now becomes an 
important matter, since our argument to show that all possible configurations
are $m$-equivalent involves cancelling pairs of nodal intersections.

More precisely, provided that $d\le n(n-1)$, by applying formula 
(\ref{eq:k2k}) we obtain a braid factorization without negative twists,
which is $m$-equivalent to the braid factorization describing a generic
algebraic map in degree $2k$, but we don't know for sure whether the 
$m$-equivalence can be realized without creating pairs of nodal
intersections between the $2n$ lines $1_\alpha,\dots,n'_\alpha$
(resp.\ $\beta$, $\gamma$). In fact, the perturbation of $V'_2\circ f_k$
that we perform near the points of $\mathcal{I}'_k$ is isotopic through 
$m$-equivalence to a generic algebraic perturbation of $V'_2\circ
f_k$, which itself would yield the usual algebraic braid monodromy 
invariants as defined by Moishezon and Teicher.

Still, it seems very unlikely that such pair creation
operations are ever needed, and it is reasonable to formulate the following
conjecture~:

\begin{conj} When $X$ is a complex algebraic manifold, the degree doubling
formula $(\ref{eq:k2k})$ is valid up to Hurwitz equivalence and global
conjugation.
\end{conj}

Motivation for this conjecture comes from the following observation. 
Assume that identifying the braid monodromy given by (\ref{eq:k2k}) 
with that of a generic algebraic map requires the creation of pairs of 
nodes. Then, considering only the relative motions of the $2n$ points 
labelled $1_\alpha,\dots,n'_\alpha$ (resp.\ $\beta$, $\gamma$), we 
obtain two factorizations of $\Delta_{2n}^2$ as a product of positive twists 
and half-twists in $B_{2n}$ which are inequivalent in a certain sense. These
two factorizations can be thought of as describing the braid monodromy of 
two symplectic nodal curves in $\CP^2$, both irreducible and of identical 
degree and genus. The braid factorization in $B_{2n}$ arising from 
(\ref{eq:k2k}) is 
easily checked to be that of an algebraic nodal curve. Therefore, the
inequivalence of the two factorizations would be a strong indication
of the possibility of constructing by purely complex algebraic methods a
counterexample to the nodal symplectic isotopy conjecture~; this would
be extremely surprising.

\section{The degree doubling formula for Lefschetz pencils}
\subsection{Braid groups and mapping class groups}\label{ss:mapg}
We now expand on the ideas in \S 5 of \cite{AK} to provide a description
of the relations between the braid monodromy of a branch curve and the
monodromy of the corresponding Lefschetz pencil. 

Recall that the Lefschetz pencils determined by approximately holomorphic
sections of $L^{\otimes k}$ are obtained from the corresponding branched
coverings simply by forgetting one of the three sections, or equivalently by
composing the covering map with the projection $\pi:\CP^2-\{pt\}\to\CP^1$.
In particular the curves making up the pencil are precisely the preimages of
the fibers of $\pi$ by the branched covering, and the base points of the
pencil are the preimages of the pole of the projection $\pi$.

Consider as previously the branched covering $f_k:X\to\CP^2$. Call $n$
its degree and $d$ the degree of its branch curve $D_k$, and let
$\theta:F_d=\pi_1(\C-\{q_1,\dots,q_d\})\to S_n$ be the corresponding 
geometric monodromy representation. The map $\theta$ determines a
simple $n$-fold covering of $\CP^1$ branched at $q_1,\dots,q_d$~;
we will denote this covering by $u:\Sigma_g\to\CP^1$, where $\Sigma_g$ is a 
Riemann surface of genus $g=1-n+(d/2)$. 

It is important for our purposes to observe
that the Riemann surface $\Sigma_g$ naturally comes with $n$ marked
points, corresponding to the base points of the Lefschetz pencil~:
these $n$ points are precisely the preimages by $u$ of the point at
infinity in $\CP^1$. In particular, rather than simply working in the
mapping class group $M_g$ of $\Sigma_g$ in the usual way, we will
consider the mapping class group $M_{g,n}$ of a Riemann surface of
genus $g$ with $n$ boundary components, i.e.\ the set of isotopy classes 
of diffeomorphisms of the complement of $n$ discs centered at the
given points in $\Sigma_g$ which fix each of the $n$ boundary components
(or equivalently, diffeomorphisms of $\Sigma_g$ which fix the $n$ marked
points and whose tangent map at each of these points is the identity).
Describing a Lefschetz pencil by a word in $M_{g,n}$ provides a more
complete picture than the usual description using $M_g$, as it also 
accounts for the relative positions of the base points of the pencil
with respect to the various vanishing cycles.
\medskip

Recall the following construction from \cite{AK}~: let
$\mathcal{C}_n(q_1,\dots,q_d)$ be the (finite) set of
all surjective group homomorphisms $F_d\to S_n$ which map each of the
geometric generators $\gamma_1,\dots,\gamma_d$ of $F_d$ to a transposition 
and map their product $\gamma_1\cdots\gamma_d$ to the identity element in 
$S_n$. Each element of $\mathcal{C}_n(q_1,\dots,q_d)$ determines a
simple $n$-fold covering of $\CP^1$ branched at $q_1,\dots,q_d$.

Denote by
$\mathcal{X}_d$ the space of configurations of $d$ distinct points in
the plane. The set of all simple $n$-fold coverings of $\CP^1$ with $d$ branch
points and such that no branching occurs above the point at infinity
can be thought of as a covering $\tilde{\mathcal{X}}_{d,n}$ above 
$\mathcal{X}_d$,
whose fiber above the configuration $\{q_1,\dots,q_d\}$ identifies with
$\mathcal{C}_n(q_1,\dots,q_d)$. The braid group $B_d$ identifies with the
fundamental group of $\mathcal{X}_d$, and therefore $B_d$ acts on the fiber
$\mathcal{C}_n(q_1,\dots,q_d)$ by deck transformations of the covering
$\tilde{\mathcal{X}}_{d,n}$. 

Define the subgroup $B_d^0(\theta)$ as the set of
all the loops in $\mathcal{X}_d$ whose lift at the point $p_\theta\in
\tilde{\mathcal{X}}_{d,n}$ corresponding to the covering described by 
$\theta$ is a closed loop in $\tilde{\mathcal{X}}_{d,n}$, i.e.\ the set
of all braids which act on $F_d=\pi_1(\C-\{q_1,\dots,q_d\})$ in a manner 
compatible with the covering structure defined by $\theta$.
Denoting by $Q_*$ the action of a braid $Q$ on $F_d$, it is easy to
check that $B_d^0(\theta)$ is the set of all braids $Q$ such that
$\theta\circ Q_*=\theta$.

There exists a natural (tautologically defined) bundle $\mathcal{Y}_{d,n}$
over $\tilde{\mathcal{X}}_{d,n}$ whose fiber is a Riemann surface of genus
$g$. Each of these Riemann surfaces comes naturally as a branched covering
of $\CP^1$, and carries $n$ distinct marked points -- the preimages of the 
point at infinity in $\CP^1$ by the covering.

Given an element $Q$ of $B_d^0(\theta)\subset B_d$, it can be lifted to
$\tilde{\mathcal{X}}_{d,n}$ as a loop based at the point $p_\theta$, and
the monodromy of the fibration $\mathcal{Y}_{d,n}$ around this loop defines
an element of the mapping class group $M_{g,n}$ of a Riemann surface of
genus $g$ with $n$ boundary components, which 
we will call $\theta_*(Q)$.
More intuitively, viewing $Q$ as a compactly supported diffeomorphism of 
the plane preserving $\{q_1,\dots,q_d\}$, the fact that $Q\in B_d^0(\theta)$
means that the diffeomorphism representing $Q$ can be lifted via the
covering $u:\Sigma_g\to\CP^1$ to a diffeomorphism of $\Sigma_g$, whose
class in the mapping class group is $\theta_*(Q)$.\medskip

It is easy to check that the image of the braid monodromy homomorphism
is contained in $B^0_d(\theta)$~: this is because the geometric 
monodromy representation $\theta$ factors through $\pi_1(\CP^2-D_k)$, 
on which the action of the braids arising in the monodromy is clearly trivial. 
Therefore, we can take the image of the braid factorization by the map
$\theta_*$ and obtain a factorization in the mapping class group $M_{g,n}$.
As observed in \cite{AK}, all the factors of degree $\pm 2$ or $3$ in the
factorization lie in the kernel of $\theta_*$~; 
therefore, the only remaining terms are those corresponding
to the tangency points of the branch curve $D_k$, and each of these is
a Dehn twist. 

Recall from \cite{AK} that the image in the mapping class group $M_{g,n}$ 
of a half-twist $Q\in B^0_d(\theta)$ can be constructed as follows. 
Call $\gamma$ the path in $\C$ joining two of the branch points 
(say $q_i$ and $q_j$) which describes the half-twist $Q$ ($\gamma$ 
is the path along which the twisting occurs). Among the $n$ lifts of 
$\gamma$ to $\Sigma_g$, only two hit the branch points of the covering~; 
these two lifts have common end points, and together they define a loop 
$\delta$ in $\Sigma_g$. Equivalently, one may also define $\delta$ as one 
of the two non-trivial lifts of the boundary of a small tubular 
neighborhood of $\gamma$ in $\C$. In any case, one easily checks that
the element $\theta_*(Q)$ in $M_{g,n}$ is a positive Dehn twist along
the loop $\delta$ (see Proposition 4 of \cite{AK}).

As a consequence, one obtains the usual description of the
monodromy of the Lefschetz pencil as a word in the mapping class group
whose factors are positive Dehn twists. However,
as observed by Smith in \cite{smith}, the product of all these Dehn twists
is not the identity element in $M_{g,n}$, because after blowing up the 
pencil at its base points one obtains a Lefschetz fibration in which 
the exceptional sections have the non-trivial normal bundle $O(-1)$. 
Instead, the product of all the factors is equal to $\theta_*(\Delta_d^2)$,
which is itself equal to the product of $n$ positive Dehn twists,
one along a small loop around each of the $n$ base points of the pencil.
\medskip

It follows from the above considerations that we can lift the degree
doubling formula for braid monodromies obtained in \S 3 and obtain a similar
formula for Lefschetz pencils. The task is made even easier by the fact that
we only need to consider the tangency points of the branch curves.

We now introduce the general setup for the degree doubling formula.
To start with, recall that the branch curve $D_{2k}$ is of degree
$\bar{d}=2d+6n$, while the degree of the covering $f_{2k}$ is $4n$.
Recall from \S \ref{ss:3pt} the relation between
the geometric monodromy factorizations $\theta_{2k}:F_{\bar{d}}\to S_{4n}$
and $\theta_k:F_d\to S_n$~: as previously, view
the $4n$ sheets of $f_{2k}$ as four groups of $n$ sheets labelled
$i_a,i_b,i_c,i_d$, $1\le i\le n$, and use the same labelling of the
branch points as in \S 3. With these notations, the transpositions
in $S_{4n}$ corresponding to the geometric generators around
$1,\dots,d,1',\dots,d'$ are directly given by the geometric monodromy 
representation $\theta_k$ associated to $D_k$~: given $1\le r\le d$, 
if $\theta_k$ maps the $r$-th geometric generator to the 
transposition $(ij)$ in $S_n$ then, calling $\gamma_r$ and $\gamma_{r'}$ the 
geometric generators in $F_{\bar d}$ corresponding to $r$ and $r'$, 
one gets $\theta_{2k}(\gamma_r)=\theta_{2k}(\gamma_{r'})=(i_aj_a)$.
Moreover, each of the $n$ copies of $V_2$ connects four sheets to each other,
one in each group of $n$~: the geometric
generators around $i_\alpha$, $i'_\alpha$, $i_\beta$, $i'_\beta$, $i_\gamma$ 
and $i'_\gamma$ are mapped by $\theta_{2k}$ to 
$(i_ai_b)$, $(i_ci_d)$, $(i_ai_c)$, $(i_bi_d)$, $(i_ai_d)$ and $(i_bi_c)$
respectively, for all $1\le i\le n$.

As a consequence, $\theta_{2k}$ determines a $4n$-fold branched covering
$\bar{u}:\Sigma_{\bar g}\to \CP^1$, with $\bar g=2g+n-1$, whose structure
is as follows. First, the preimage of a disc $\mathcal{D}$ containing the 
$d$ points labelled $1,\dots,d$ consists of $3n+1$ components. One of 
these components (the sheets $1_a,\dots,n_a$) is a $n$-fold covering 
identical to the one described by $\theta_k$, i.e.\ it naturally 
identifies with the fiber $\Sigma_g$ of the Lefschetz pencil associated 
to $f_k$, with $n$ small discs removed. These punctures correspond to 
the preimages of a small disc around the point at infinity in the 
covering $u:\Sigma_g\to\CP^1$, i.e.\ they correspond to small discs 
around the base points in $\Sigma_g$. The other $3n$ components of
$\bar{u}^{-1}(\mathcal{D})$, in which no branching occurs, are 
topologically trivial. 

The same picture also describes the preimage of a disc $\mathcal{D}'$
containing the $d$ points labelled $1',\dots,d'$~: there is one non-trivial
component which can be identified with $\Sigma_g$ punctured at its base
points, and the other $3n$ components are just plain discs.

Finally, the preimage by $\bar{u}$ of the cylinder $\CP^1-(\mathcal{D}
\cup\mathcal{D}')$ consists of $n$ components, each of which is a
four-sheeted covering branched at six points, i.e.\ topologically a
sphere with eight punctures. Actually, each of these $n$ components 
may be thought of as the fiber of the Lefschetz pencil corresponding
to the covering $V_2$ (since we restrict ourselves to a cylinder we
get eight punctures). For each $i\in\{1,\dots,n\}$ the corresponding
component of $\bar{u}^{-1}(\CP^1-(\mathcal{D}\cup\mathcal{D}'))$
connects together the non-trivial components of $\bar{u}^{-1}(\mathcal{D})$
and $\bar{u}^{-1}(\mathcal{D}')$ with the trivial components corresponding
to the sheets $i_b$, $i_c$ and $i_d$.

In the end the Riemann surface $\Sigma_{\bar g}$ can be thought of as
two copies of $\Sigma_g$ glued together at the $n$ base points.
This description coincides exactly with the one obtained by Smith in
\cite{smith} via more direct methods.

\subsection{The degree doubling formula for Lefschetz pencils}
\label{ss:mapgthm}
In order to simplify the description of the degree doubling formula
for Lefschetz pencils, we want to slightly modify the setup of \S 3.

First, we want to choose a different picture for $\theta_k$~: recall
that global conjugations in $B_d$ make it possible to choose the most
convenient geometric monodromy representation $\theta_k:F_d\to S_n$.
As a consequence we chose in \S 3 a setup that made the final answer
(\ref{eq:k2k}) relatively easy to express, but as observed in the
remark at the end of \S \ref{ss:braidthm} we could just as well have worked 
with any other choice of $\theta_k$, the only price being a slightly more 
complicated expression for the degree doubling formula. Note that
the change of $\theta_k$ only affects factors of degree $\pm 2$ in
the formula, and therefore the half-twists which are relevant for
our purposes are not affected.

Here we want to choose $\theta_k$ in such a way that the $i$-th geometric
generator $\gamma_i$ is mapped to the transposition $(1,2)$ if 
$i\le d-2(n-1)=2g$, and $\theta_k(\gamma_{d-2j})=\theta(\gamma_{d-2j-1})=
(n-j-1,n-j)$ for all $j\le n-2$. In other words, the transpositions
$\theta_k(\gamma_i)$ correspond to the factorization
$$\textstyle \mathrm{Id}=(1,2)^{2g}\cdot\prod\limits_{i=1}^{n-1} (i,i+1)^2$$
in $S_n$.
Another change that we want to make is in the ordering of the 
$\bar d=2d+6n$ points that appear in the diagrams of \S 3 along
the real axis. Namely, we want to replace the ordering
$1,\dots,d$, $1',\dots,d'$, $1_\alpha,\dots,n'_\gamma$ used in \S 3
by the new ordering $1,\dots,d$, $1_\alpha,\dots,n'_\gamma$, $d',\dots,1'$.
This is done by first moving the $d$ points $1',\dots,d'$ {\em clockwise}
around the points $1_\alpha,\dots,n'_\gamma$ by a half-turn, and then
by rotating a disc containing the $d$ points $1',\dots,d'$ 
{\em counterclockwise} by a half-turn.

Finally, in order to better visualize the positions of the base points of
the pencil (the $4n$ marked points on $\Sigma_{\bar g}$),
we want to move the fiber in which they lie from the point at infinity
in $\CP^1$ back into our picture. We choose to move the base points so
that they correspond to the preimages of a point $b$ on the real axis lying
inbetween the point labelled $d$ and the point labelled $1_\alpha$.
The motion bringing the point at infinity to $b$ is performed along a
vertical line in the upper half-plane (this motion of course affects some
of the braids, but it was chosen in such a way that the resulting changes
are minimal). 

The effect of all these changes is to make the covering $\bar u:\Sigma_{\bar g}
\to \CP^1$ easier to visualize, while simplifying the paths corresponding
to the half-twists in (\ref{eq:k2k}). The picture is the following~:

\begin{center}
\epsfbox{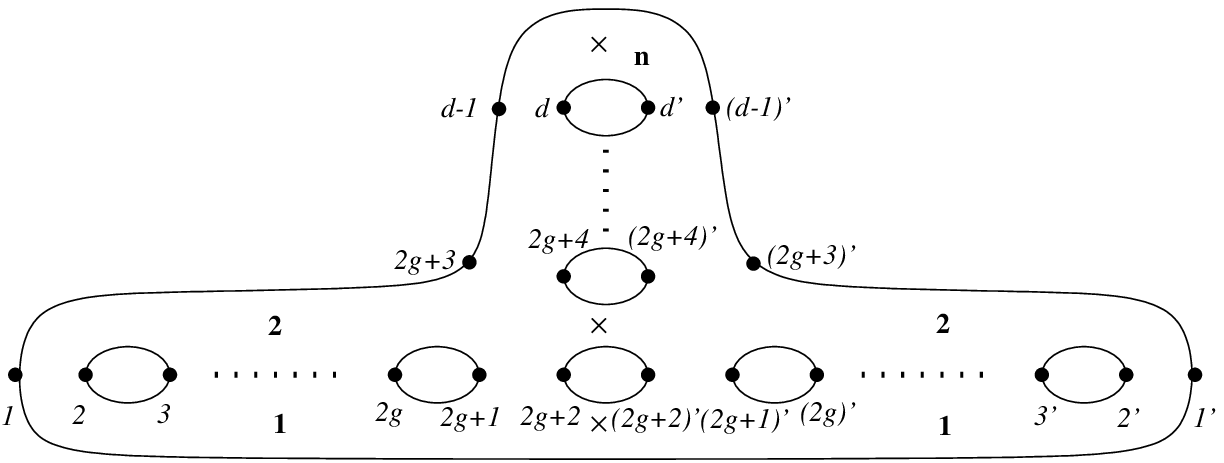}
\end{center}

In this picture, the labels in italics correspond to branch points and those 
in boldface correspond to the sheets of the covering~; for simplicity we have 
omitted the branch points $1_\alpha, \dots,n'_\gamma$, which should be placed 
in the necks joining the two halves, and the $3n$ other sheets which
do not contribute to the topology. When the $3n$ sheets $1_b,\dots,n_d$
are collapsed, the corresponding base points are brought back to the
sheets $1_a,\dots,n_a$ near the branch points $1_\alpha,\dots,n'_\gamma$~;
therefore, on the picture each $\times$ mark corresponds to four base
points.

In order to understand the Lefschetz pencil corresponding to $f_{2k}$,
we need to place the various half-twists appearing in the braid
factorization of $D_{2k}$ on this picture. A first set of half-twists
comes from the braid factorization of $D_k$. These half-twists correspond
exactly to the Dehn twists appearing in monodromy of the Lefschetz pencil
for $f_k$, after a suitable embedding of $M_{g,n}$ into the mapping class
group $M_{\bar g,4n}$. Recall that the braid factorization in $B_d$ 
corresponding to
$D_k$ is embedded into $B_{\bar d}$ by considering a disc $\mathcal{D}$ 
containing the $d$ points labelled $1,\dots,d$. Therefore, the corresponding
embedding of the mapping class group $M_{g,n}$ into the larger mapping 
class group $M_{\bar g,4n}$ is geometrically realized by the embedding
into $\Sigma_{\bar g}$ of the main connected component of 
$\bar{u}^{-1}(\mathcal{D})$, which as we know from \S \ref{ss:mapg}
naturally identifies with the Riemann 
surface $\Sigma_g$ punctured at each of the $n$ base points.
On the above picture of $\Sigma_{\bar g}$ this corresponds to the
left half of the diagram.

Observe that all the other half-twists appearing in the braid factorization
for $D_{2k}$ are completely standard and depend only on $d$ and $n$ rather 
than on the actual topology of the manifold $X$. Therefore, the degree 
doubling formula for Lefschetz pencils is once again a universal formula~:
the word in $M_{\bar g,4n}$ describing the Lefschetz 
pencil in degree $2k$ is obtained by embedding the word describing the 
pencil in degree $k$ via the above-described map from $M_{g,n}$
into $M_{\bar g,4n}$ and adding to it a completely standard set of
Dehn twists which depends only on $g$ and $n$ but not on the actual
topology of the manifold $X$. This observation was already made by
Ivan Smith in \cite{smith}.
\medskip

The extra half-twists appearing in the degree doubling formula for braid
monodromies are $\check{Z}_{ii'}$ and $\hat{Z}_{ii'}$ for $1\le i\le d$,
and $\tilde{Z}_{i_\alpha i'_\alpha;(i_\beta i'_\beta)}$,
$\tilde{Z}_{i_\alpha i'_\alpha;(i_\gamma i'_\gamma)}$ and
$\tilde{Z}_{i_\beta i'_\beta;(i_\gamma i'_\gamma)}$ for $1\le i\le n$,
as described in \S \ref{ss:braidthm} (their total number $2d+3n$ is in 
agreement with an easy calculation of Euler-Poincar\'e characteristics). 
We will now describe the Dehn twists corresponding to these half-twists.

After the global conjugation described above, $\check{Z}_{ii'}$ becomes
a half-twist along the following path~:

\begin{center}
\setlength{\unitlength}{0.25in}
\begin{picture}(13,2)(0,-1)
\put(0,0){\circle*{0.1}}
\put(0,-0.2){\makebox(0,0)[ct]{$1$}}
\multiput(0.5,0)(0.2,0){6}{\circle*{0.02}}
\put(2,0){\circle*{0.1}}
\put(2,-0.2){\makebox(0,0)[ct]{$i$}}
\multiput(2.5,0)(0.2,0){6}{\circle*{0.02}}
\put(4,0){\circle*{0.1}}
\put(4,-0.2){\makebox(0,0)[ct]{$d$}}
\put(5,0.2){\makebox(0,0)[ct]{$\times$}}
\put(5,0.3){\makebox(0,0)[cb]{$b$}}
\put(6,0){\circle*{0.1}}
\put(6,0.2){\makebox(0,0)[cb]{$1_\alpha$}}
\multiput(6.5,0)(0.2,0){6}{\circle*{0.02}}
\put(8,0){\circle*{0.1}}
\put(8.1,0.05){\makebox(0,0)[cb]{$n'_\gamma$}}
\put(9,0){\circle*{0.1}}
\put(9,0.2){\makebox(0,0)[cb]{$d'$}}
\multiput(9.5,0)(0.2,0){6}{\circle*{0.02}}
\put(11,0){\circle*{0.1}}
\put(11,0.2){\makebox(0,0)[cb]{$i'$}}
\multiput(11.5,0)(0.2,0){6}{\circle*{0.02}}
\put(13,0){\circle*{0.1}}
\put(13,0.2){\makebox(0,0)[cb]{$1'$}}
\qbezier(2,0)(2,0.8)(3.25,0.8)
\qbezier(3.25,0.8)(4.5,0.8)(4.5,0)
\qbezier(4.5,0)(4.5,-0.8)(5.8,-0.8)
\put(5.8,-0.8){\line(1,0){4.2}}
\qbezier(10,-0.8)(11,-0.8)(11,0)
\end{picture}
\end{center}

Its lift to the mapping class group $M_{\bar g,4n}$ is a Dehn twist
that we will call $\check\tau_i$, and which can be represented as
follows when $i$ is even and $i\le 2g$~:

\begin{center}
\epsfbox{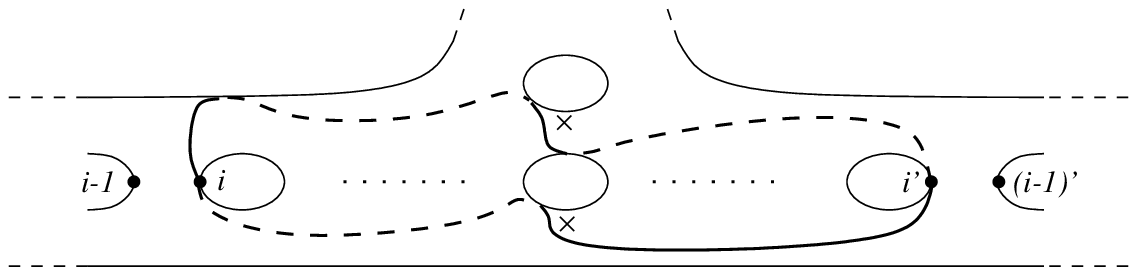}
\end{center}

For $i$ odd and $i\le 2g+1$, the picture describing $\check\tau_i$ becomes
the following~:

\begin{center}
\epsfbox{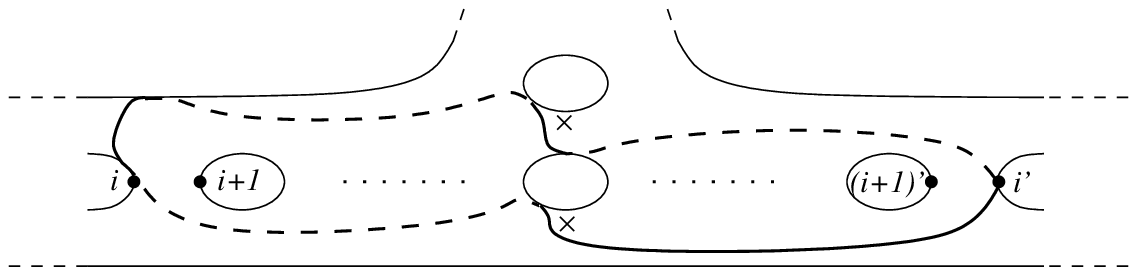}
\end{center}

When $i=1$ the undrawn parts on both sides of the picture are just
discs and the picture can therefore be slightly simplified~; conversely,
when $i=2g+1$ the points labelled $(i+1)$ and $(i+1)'$ are immediately
on both sides of the central neck rather than as pictured.

For $i$ even and $i\ge 2g+2$, $\check\tau_i$ is described by the following
picture (the two necks shown correspond to the sheets numbered $s$ and
$s+1$, where $s=\frac{1}{2}(i-2g)\ge 1$)~:

\begin{center}
\epsfbox{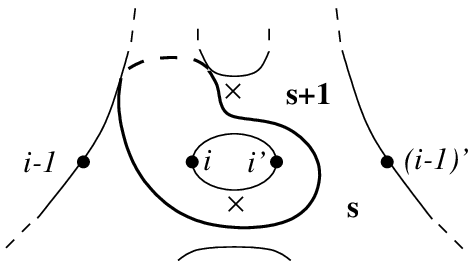}
\end{center}

Finally, when $i$ is odd and $i\ge 2g+3$, the picture describing
$\check\tau_i$ becomes the following (the two necks shown correspond to the
sheets numbered $s$ and $s+1$, where $s=\frac{1}{2}(i+1-2g)\ge 2$)~:

\begin{center}
\epsfbox{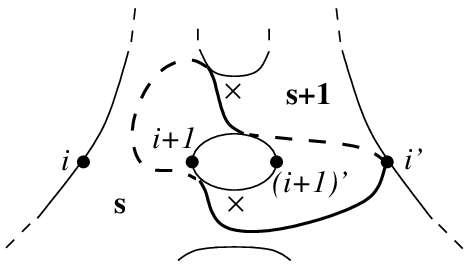}
\end{center}

We now turn to $\hat{Z}_{ii'}$, which after the above-described global
conjugation becomes a half-twist along the following path~:

\begin{center}
\setlength{\unitlength}{0.25in}
\begin{picture}(13,2)(0,-1)
\put(0,0){\circle*{0.1}}
\put(0,-0.2){\makebox(0,0)[ct]{$1$}}
\multiput(0.5,0)(0.2,0){6}{\circle*{0.02}}
\put(2,0){\circle*{0.1}}
\put(2,-0.2){\makebox(0,0)[ct]{$i$}}
\multiput(2.5,0)(0.2,0){6}{\circle*{0.02}}
\put(4,0){\circle*{0.1}}
\put(4,-0.2){\makebox(0,0)[ct]{$d$}}
\put(5,0.2){\makebox(0,0)[ct]{$\times$}}
\put(5,0.3){\makebox(0,0)[cb]{$b$}}
\put(6,0){\circle*{0.1}}
\put(6,-0.2){\makebox(0,0)[ct]{$1_\alpha$}}
\multiput(6.5,0)(0.2,0){6}{\circle*{0.02}}
\put(8,0){\circle*{0.1}}
\put(8,-0.2){\makebox(0,0)[ct]{$n'_\gamma$}}
\put(9,0){\circle*{0.1}}
\put(9,-0.2){\makebox(0,0)[ct]{$d'$}}
\multiput(9.5,0)(0.2,0){6}{\circle*{0.02}}
\put(11,0){\circle*{0.1}}
\put(10.75,-0.2){\makebox(0,0)[ct]{$i'$}}
\multiput(11.5,0)(0.2,0){6}{\circle*{0.02}}
\put(13,0){\circle*{0.1}}
\put(12.8,-0.15){\makebox(0,0)[ct]{$1'$}}
\qbezier(2,0)(2,0.8)(3.25,0.8)
\put(3.25,0.8){\line(1,0){0.5}}
\put(6.25,0.8){\line(1,0){6}}
\qbezier(3.75,0.8)(4.2,0.8)(4.6,0)
\qbezier(4.6,0)(5,-0.8)(5.4,0)
\qbezier(5.4,0)(5.8,0.8)(6.25,0.8)
\qbezier(12.25,0.8)(13.5,0.8)(13.5,0)
\qbezier(13.5,0)(13.5,-0.8)(12.5,-0.8)
\qbezier(12.5,-0.8)(11,-0.8)(11,0)
\end{picture}
\end{center}

This path can be homotoped into the following one, which goes through
the point at infinity in $\CP^1$~:

\begin{center}
\setlength{\unitlength}{0.25in}
\begin{picture}(13,3.2)(0,-1.6)
\put(0,0){\circle*{0.1}}
\put(0,-0.2){\makebox(0,0)[ct]{$1$}}
\multiput(0.5,0)(0.2,0){6}{\circle*{0.02}}
\put(2,0){\circle*{0.1}}
\put(2,-0.2){\makebox(0,0)[ct]{$i$}}
\multiput(2.5,0)(0.2,0){6}{\circle*{0.02}}
\put(4,0){\circle*{0.1}}
\put(4,-0.2){\makebox(0,0)[ct]{$d$}}
\put(5,0.2){\makebox(0,0)[ct]{$\times$}}
\put(5,0.3){\makebox(0,0)[cb]{$b$}}
\put(6,0){\circle*{0.1}}
\put(6,0.2){\makebox(0,0)[cb]{$1_\alpha$}}
\multiput(6.5,0)(0.2,0){6}{\circle*{0.02}}
\put(8,0){\circle*{0.1}}
\put(8.1,0.05){\makebox(0,0)[cb]{$n'_\gamma$}}
\put(9,0){\circle*{0.1}}
\put(9,0.2){\makebox(0,0)[cb]{$d'$}}
\multiput(9.5,0)(0.2,0){6}{\circle*{0.02}}
\put(11,0){\circle*{0.1}}
\put(11,0.2){\makebox(0,0)[cb]{$i'$}}
\multiput(11.5,0)(0.2,0){6}{\circle*{0.02}}
\put(13,0){\circle*{0.1}}
\put(13,0.2){\makebox(0,0)[cb]{$1'$}}
\qbezier(2,0)(2,0.8)(3.25,0.8)
\qbezier(5.9,-1.6)(5.9,-0.8)(6.8,-0.8)
\qbezier(3.25,0.8)(4.5,0.8)(4.5,0)
\qbezier(4.5,0)(4.5,-0.5)(5.1,-0.5)
\qbezier(5.1,-0.5)(5.5,-0.5)(5.5,1.2)
\multiput(5.9,-1.6)(-0.025,0.175){17}{\circle*{0.02}}
\put(6.8,-0.8){\line(1,0){3.2}}
\qbezier(10,-0.8)(11,-0.8)(11,0)
\end{picture}
\end{center}
Therefore, the Dehn twists $\hat{\tau}_i\in M_{\bar g,4n}$ obtained by
lifting $\hat{Z}_{ii'}$ only differ from $\check{\tau}_i$ by a twisting
in each of the necks joining the two halves of $\Sigma_{\bar g}$. As
a result, we get the following pictures (using the same notations as for
$\check{\tau}_i$)~:

\begin{center}
\epsfbox{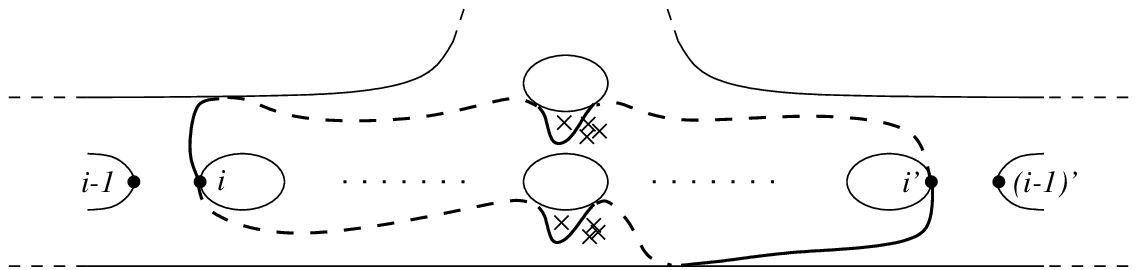}
\end{center}

\begin{center}
\epsfbox{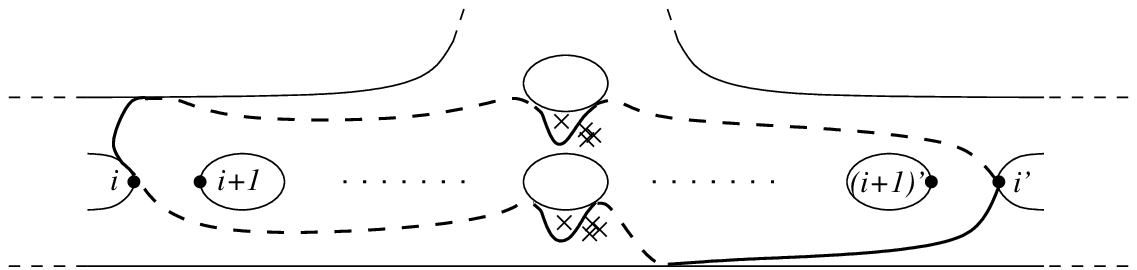}
\end{center}

The first picture corresponds to the case $i$ even, $i\le 2g$~; the second
one to $i$ odd, $i\le 2g+1$. In each of the two necks, the vanishing 
loop circles around the base point corresponding to the sheet labelled
$1_a$ (resp.\ $2_a$), but not around those corresponding to sheets 
$1_b$, $1_c$ and $1_d$ (resp.\ $2_b$, $2_c$, $2_d$).
When $i\ge 2g+2$, the pictures become the
following (the left one is for even $i$, the right one for odd $i$)~:

\begin{center}
\epsfbox{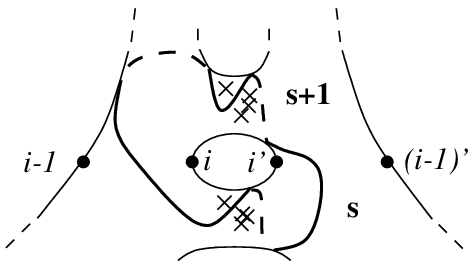}\qquad
\epsfbox{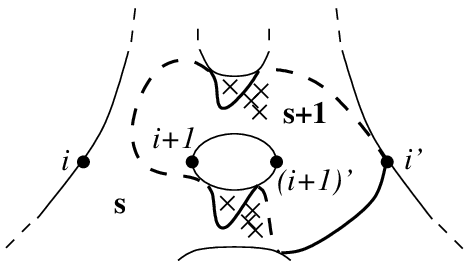}
\end{center}

We now turn to the half-twists 
$\tilde{Z}_{i_\alpha i'_\alpha;(i_\beta i'_\beta)}$,
$\tilde{Z}_{i_\alpha i'_\alpha;(i_\gamma i'_\gamma)}$ and
$\tilde{Z}_{i_\beta i'_\beta;(i_\gamma i'_\gamma)}$ ($1\le i\le n$).
To simplify the diagrams we only represent the relevant points,
i.e. we forget $j_\alpha,j'_\alpha,j_\beta,j'_\beta,j_\gamma,j'_\gamma$
for $j\neq i$ as these points do not play any role. Moreover, we use the 
observation
that, for the purposes of computing the corresponding Dehn twists, we are
allowed to move a path across a branch point if the corresponding sheets
of the covering are distinct. Finally, we further simplify the diagrams
by allowing ourselves to draw paths which go through the point at infinity
in $\CP^1$. With all these simplifications, we get the following diagrams~:

\begin{center}
\setlength{\unitlength}{0.25in}
\begin{picture}(15,3.2)(0,-1.6)
\put(13.5,0){\makebox(0,0)[lc]{
$\tilde{Z}_{i_\alpha i'_\alpha;(i_\beta i'_\beta)}$}}
\put(0,0){\circle*{0.1}}
\put(0,-0.2){\makebox(0,0)[ct]{$1$}}
\multiput(0.5,0)(0.2,0){6}{\circle*{0.02}}
\put(2,0){\circle*{0.1}}
\put(2,-0.2){\makebox(0,0)[ct]{$d$}}
\put(3,0.2){\makebox(0,0)[ct]{$\times$}}
\put(3,0.3){\makebox(0,0)[cb]{$b$}}
\put(4,0){\circle*{0.1}}
\put(4.1,0.1){\makebox(0,0)[cb]{$i_\alpha$}}
\put(5,0){\circle*{0.1}}
\put(4.7,-0.1){\makebox(0,0)[cb]{$i'_\alpha$}}
\put(6,0){\circle*{0.1}}
\put(6,0.2){\makebox(0,0)[cb]{$i_\beta$}}
\put(7,0){\circle*{0.1}}
\put(7,0.2){\makebox(0,0)[cb]{$i'_\beta$}}
\put(8,0){\circle*{0.1}}
\put(8,-0.2){\makebox(0,0)[ct]{$i_\gamma$}}
\put(9,0){\circle*{0.1}}
\put(9,-0.15){\makebox(0,0)[ct]{$i'_\gamma$}}
\put(10,0){\circle*{0.1}}
\put(10,0.2){\makebox(0,0)[cb]{$d'$}}
\multiput(10.5,0)(0.2,0){6}{\circle*{0.02}}
\put(12,0){\circle*{0.1}}
\put(12,0.2){\makebox(0,0)[cb]{$1'$}}
\qbezier(5,0)(5,0.8)(4.2,0.8)
\qbezier(4.2,0.8)(3.5,0.8)(3.5,0)
\qbezier(3.5,0)(3.5,-1)(4.5,-1)
\qbezier(4,0)(4,-0.5)(4.7,-0.5)
\qbezier(4.7,-0.5)(5.5,-0.5)(5.5,0)
\qbezier(5.5,0)(5.5,1.3)(6.5,1.3)
\put(6.5,1.3){\line(1,0){2}}
\qbezier(8.5,1.3)(9.25,1.3)(9.4,1.6)
\qbezier[20](9.4,1.6)(9.7,0)(9.4,-1.6)
\qbezier(9.4,-1.6)(9.3,0.8)(8.8,0.8)
\qbezier(8.8,0.8)(7.5,0.8)(7.5,0)
\qbezier(7.5,0)(7.5,-1)(6.5,-1)
\put(4.5,-1){\line(1,0){2}}
\end{picture}
\end{center}

\begin{center}
\setlength{\unitlength}{0.25in}
\begin{picture}(15,3.2)(0,-1.6)
\put(13.5,0){\makebox(0,0)[lc]{
$\tilde{Z}_{i_\alpha i'_\alpha;(i_\gamma i'_\gamma)}$}}
\put(0,0){\circle*{0.1}}
\put(0,-0.2){\makebox(0,0)[ct]{$1$}}
\multiput(0.5,0)(0.2,0){6}{\circle*{0.02}}
\put(2,0){\circle*{0.1}}
\put(2,-0.2){\makebox(0,0)[ct]{$d$}}
\put(3,0.2){\makebox(0,0)[ct]{$\times$}}
\put(3,0.3){\makebox(0,0)[cb]{$b$}}
\put(4,0){\circle*{0.1}}
\put(4,-0.2){\makebox(0,0)[ct]{$i_\alpha$}}
\put(5,0){\circle*{0.1}}
\put(5,-0.15){\makebox(0,0)[ct]{$i'_\alpha$}}
\put(6,0){\circle*{0.1}}
\put(6,-0.2){\makebox(0,0)[ct]{$i_\beta$}}
\put(7,0){\circle*{0.1}}
\put(6.9,-0.15){\makebox(0,0)[ct]{$i'_\beta$}}
\put(8,0){\circle*{0.1}}
\put(8,-0.2){\makebox(0,0)[ct]{$i_\gamma$}}
\put(9,0){\circle*{0.1}}
\put(9,-0.15){\makebox(0,0)[ct]{$i'_\gamma$}}
\put(10,0){\circle*{0.1}}
\put(10,0.2){\makebox(0,0)[cb]{$d'$}}
\multiput(10.5,0)(0.2,0){6}{\circle*{0.02}}
\put(12,0){\circle*{0.1}}
\put(12,0.2){\makebox(0,0)[cb]{$1'$}}
\qbezier(4,0)(4,1)(5,1)
\put(5,1){\line(1,0){1.5}}
\qbezier(6.5,1)(7.25,1)(7.4,1.6)
\qbezier[20](7.4,1.6)(7.7,0)(7.4,-1.6)
\qbezier(7.4,-1.6)(7.4,0.5)(6.8,0.5)
\qbezier(5,0)(5,0.5)(6,0.5)
\put(6,0.5){\line(1,0){0.8}}
\end{picture}
\end{center}

\begin{center}
\setlength{\unitlength}{0.25in}
\begin{picture}(15,3.6)(0,-1.6)
\put(13.5,0){\makebox(0,0)[lc]{
$\tilde{Z}_{i_\beta i'_\beta;(i_\gamma i'_\gamma)}$}}
\put(0,0){\circle*{0.1}}
\put(0,-0.2){\makebox(0,0)[ct]{$1$}}
\multiput(0.5,0)(0.2,0){6}{\circle*{0.02}}
\put(2,0){\circle*{0.1}}
\put(2,-0.2){\makebox(0,0)[ct]{$d$}}
\put(3,0.2){\makebox(0,0)[ct]{$\times$}}
\put(3,0.3){\makebox(0,0)[cb]{$b$}}
\put(4,0){\circle*{0.1}}
\put(4,-0.2){\makebox(0,0)[ct]{$i_\alpha$}}
\put(5,0){\circle*{0.1}}
\put(5,-0.15){\makebox(0,0)[ct]{$i'_\alpha$}}
\put(6,0){\circle*{0.1}}
\put(6,-0.2){\makebox(0,0)[ct]{$i_\beta$}}
\put(7,0){\circle*{0.1}}
\put(6.9,-0.15){\makebox(0,0)[ct]{$i'_\beta$}}
\put(8,0){\circle*{0.1}}
\put(8,-0.2){\makebox(0,0)[ct]{$i_\gamma$}}
\put(9,0){\circle*{0.1}}
\put(9,-0.15){\makebox(0,0)[ct]{$i'_\gamma$}}
\put(10,0){\circle*{0.1}}
\put(10,0.2){\makebox(0,0)[cb]{$d'$}}
\multiput(10.5,0)(0.2,0){6}{\circle*{0.02}}
\put(12,0){\circle*{0.1}}
\put(12,0.2){\makebox(0,0)[cb]{$1'$}}
\qbezier(6,0)(6,0.8)(6.8,0.8)
\qbezier(6.8,0.8)(7.25,0.8)(7.4,1.6)
\qbezier[20](7.4,1.6)(7.7,0)(7.4,-1.6)
\qbezier(7.4,-1.6)(7.4,0)(7,0)
\end{picture}
\end{center}
It is now clear that the only relevant parts of $\Sigma_{\bar g}$ are
the sheets labelled $i_b$, $i_c$, $i_d$ of the covering, as well as
the part of the sheet labelled $i_a$ that lies inbetween the points
$1,\dots,d$ and $d',\dots,1'$. In particular, the loops we obtain are
entirely located in the $i$-th neck joining the two halves of 
$\Sigma_{\bar g}$~; if we forget about the base points, the
Dehn twists $\tau_{i,\alpha\beta}$, $\tau_{i,\alpha\gamma}$ and 
$\tau_{i,\beta\gamma}$ corresponding to the half-twists 
$\tilde{Z}_{i_\alpha i'_\alpha;(i_\beta i'_\beta)}$,
$\tilde{Z}_{i_\alpha i'_\alpha;(i_\gamma i'_\gamma)}$ and
$\tilde{Z}_{i_\beta i'_\beta;(i_\gamma i'_\gamma)}$ are equal to
each other, and are twists along a loop that simply goes around the
$i$-th neck joining the two halves of $\Sigma_{\bar g}$.

In the presence of the four base points lying in the sheets
$i_a$, $i_b$, $i_c$ and $i_d$ of the covering, we have to be more careful,
but it can be checked that the Dehn twists $\tau_{i,\alpha\beta}$, 
$\tau_{i,\alpha\gamma}$ and $\tau_{i,\beta\gamma}$ are respectively
given by the following diagrams (only the $i$-th neck is shown~; the
base points are labelled $a$, $b$, $c$, and $d$)~:
\begin{center}
\begin{tabular}{ccc}
\epsfbox{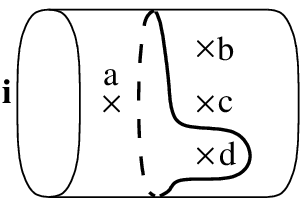} &
\epsfbox{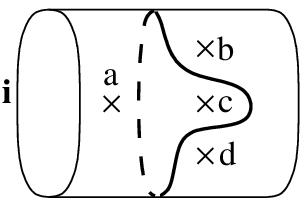} &
\epsfbox{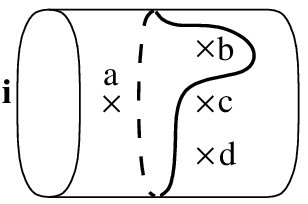}\\
$\tau_{i,\alpha\beta}$ &
$\tau_{i,\alpha\gamma}$ & 
$\tau_{i,\beta\gamma}$
\end{tabular}
\end{center}

Summarizing, we get the following result~:

\begin{theo}
Let $X$ be a compact symplectic 4-manifold, and consider the structure
of symplectic Lefschetz pencil on $X$ given by two sections
of $L^{\otimes k}$. Let $g$ be the genus of the fiber $\Sigma_g$, 
and let $n$ be the number of base points. Let $d=2g-2+2n$, and 
call $\Psi_g$ the word in the mapping class group $M_{g,n}$
describing the monodromy of this pencil.

Let $\bar{g}=2g+n-1$, and view a Riemann surface $\Sigma_{\bar g}$
of genus $\bar{g}$ as obtained by gluing together two copies of
$\Sigma_g$ at the base points. Call $\iota:M_{g,n}\to M_{\bar g,4n}$
the inclusion map discussed above. 

Then, provided that $k$ is large enough and using the notations described
above, the monodromy of the symplectic Lefschetz pencil structure obtained
on $X$ from sections of $L^{\otimes 2k}$ is given by the word $\Psi_{\bar
g}$ in the mapping class group $M_{\bar g,4n}$, where
\begin{equation}\label{eq:slp}
\Psi_{\bar g}=\prod_{i=1}^d \check{\tau}_i\cdot
\prod_{i=1}^d \hat{\tau}_i\cdot \iota(\Psi_g)\cdot
\prod_{i=1}^n \tau_{i,\alpha\beta}\cdot
\prod_{i=1}^n \tau_{i,\alpha\gamma}\cdot
\prod_{i=1}^n \tau_{i,\beta\gamma},\end{equation}
and the Dehn twists $\check{\tau}_i$, $\hat{\tau}_i$,
$\tau_{i,\alpha\beta}$, $\tau_{i,\alpha\gamma}$ and $\tau_{i,\beta\gamma}$
are as described above.
\end{theo}

\noindent {\bf Remark.} One must be aware of the fact that, in the formula
(\ref{eq:slp}), composition products are written from left to right. This
convention, which is the usual one for braid groups, is the opposite of
the usual notation for composition products when working with
diffeomorphisms (the order of the factors then needs to be reversed).

It is also worth observing that the product of the factors in
$\iota(\Psi_g)$ is almost exactly the twist by which $\hat{\tau}_i$ 
differs from $\check{\tau}_i$, the only difference being in the position
of the base points with respect to the vanishing cycle. Therefore, if we
forget about the base points, a sequence of Hurwitz moves in (\ref{eq:slp}) 
yields the following slightly simpler formula (in $M_{\bar g,0}$ instead of
$M_{\bar g,4n}$, and observing that $\tau_{i,\alpha\beta}$, 
$\tau_{i,\alpha\gamma}$ and $\tau_{i,\beta\gamma}$ are equal in $M_{\bar g,0}$):
$$\Psi_{\bar g}=\prod_{i=1}^d \check{\tau}_i\cdot
\iota(\Psi_g)\cdot
\prod_{i=1}^d \check{\tau}_i\cdot
\prod_{i=1}^n \tau_{i,\alpha\beta}^3.$$
It is clear from this expression that the Lefschetz
fibration with total space a blow-up of $X$ and monodromy $\Psi_{\bar g}$
contains many Lagrangian $(-2)$-spheres joining pairs of identical
vanishing cycles among those introduced by the degree doubling procedure; 
however these spheres collapse when the Lefschetz 
fibration is blown down along its exceptional sections, as they intersect 
non-trivially two such sections.
\medskip

The correctness of the formula (\ref{eq:slp}) can be checked easily in 
some simple examples~: for instance, a generic pencil of conics on 
$\CP^2$ has three singular fibers, and can be considered as obtained 
from a pencil of lines by the procedure described above. This corresponds 
to the limit case where $n=1$, $d=0$, $g=0$ and the word $\Psi_g$ is empty.
The three Dehn twists $\tau_{1,\alpha\beta}$, $\tau_{1,\alpha\gamma}$
and $\tau_{1,\beta\gamma}$ in $M_{0,4}$ then coincide with the well-known
picture.

Another simple example that can be considered is the case of a pencil of
curves of degree $(1,1)$ on $\CP^1\times\CP^1$. The generic fiber of this
pencil is a rational curve ($d=2$, $n=2$, $g=0$), and there are two singular
fibers. The corresponding word in $M_{0,2}$ is $\tau\cdot \tau$, where
$\tau$ is a positive Dehn twist along a simple curve separating the two
base points. The degree doubling procedure yields a word in $M_{1,8}$
consisting of $12$ Dehn twists. Forgetting the positions of the
base points, one easily checks that the reduction of this word to 
$M_{1,0}\simeq SL(2,\Z)$ is Hurwitz equivalent to the well-known monodromy
of the elliptic surface $E(1)$, which is exactly what one obtains by blowing
up the eight base points of a pencil of curves of degree $(2,2)$ on 
$\CP^1\times \CP^1$.

\end{document}